\documentclass[12pt]{iopart}
\pdfoutput=1
\usepackage{iopams}
\usepackage{amssymb}
\usepackage{latexsym}

\usepackage{graphicx}
\usepackage{color}

\begin{document}

\bibliographystyle{plain}
\def\debproof{\noindent {\bf Proof.} }
\def\finproof{\hfill {\small $\Box$} \\}
%\renewcommand{\theequation}{\arabic{section}.\arabic{equation}}
%\tableofcontents
\makeatletter % `@' now normal "letter"
\@addtoreset{equation}{section}
\makeatother  % `@' is restored as "non-letter"
\renewcommand\theequation{{\thesection}.{\arabic{equation}}}

\title[]{A Direct Imaging Method for Half-Space Inverse Elastic Scattering Problems}
\author{ Zhiming Chen, Shiqi Zhou}
\address{LSEC, Academy of
	Mathematics and Systems Science, Chinese Academy of Sciences,
	Beijing 100190, China and School of Mathematical Science, University of
        Chinese Academy of Sciences, Beijing 100049, China.}

\begin{abstract}
	We propose a direct imaging method based on the reverse time migration to reconstruct extended
	obstacles in the half space with finite aperture elastic scattering data at a fixed
	frequency. We prove the resolution of the reconstruction method in terms of the
	aperture and the depth of the obstacle embedded in the half space. The resolution
	analysis is studied by virtue of the point spread function and implies that the imaginary 
	part of the cross-correlation imaging function
	always peaks on the upper boundary of the obstacle. Numerical examples
	are included to illustrate the effectiveness of the method. 
\end{abstract}

\maketitle

\newcommand{\eps}{\varepsilon}
\newcommand{\RR}{\mathcal{R}}
\newtheorem{lem}{Lemma}[section]
\newtheorem{prop}{Proposition}[section]
\newtheorem{cor}{Corollary}[section]
\newtheorem{thm}{Theorem}[section]
\newtheorem{rem}{Remark}[section]
\newtheorem{alg}{Algorithm}[section]
\newtheorem{assum}{Assumption}[section]
\newtheorem{definition}{Definition}[section]

\newcounter{RomanNumber}
\newcommand{\MyRoman}[1]{\rm\setcounter{RomanNumber}{#1}\Roman{RomanNumber}}

\newcommand{\bL}{\mathbf{L}}
\newcommand{\bH}{\mathbf{H}}
\newcommand{\bW}{\mathbf{W}}
\newcommand{\bP}{\mathbf{P}}
\newcommand{\bQ}{\mathbf{Q}}
\newcommand{\bp}{\mathbf{p}}
\newcommand{\bq}{\mathbf{q}}
\newcommand{\uL}{u_{_{\rm L}}}
\newcommand{\vL}{v_{_{\rm L}}}
\newcommand{\tuL}{\tilde u_{_{\rm L}}}
\newcommand{\tvL}{\tilde v_{_{\rm L}}}
\newcommand{\fL}{f_{_{\rm L}}}
\newcommand{\gL}{g_{_{\rm L}}}
\newcommand{\bpL}{\bp_{_{\rm L}}}
\newcommand{\bqL}{\bq_{_{\rm L}}}
\newcommand{\tbpL}{\tilde{\bp}_{_{\rm L}}}
\newcommand{\tbqL}{\tilde{\bq}_{_{\rm L}}}
\newcommand{\tbpLf}{\tilde{\bp}_{_{\rm L,1}}}
\newcommand{\tbpLs}{\tilde{\bp}_{_{\rm L,2}}}
\newcommand{\tbqLf}{\tilde{\bq}_{_{\rm L,1}}}
\newcommand{\tbqLs}{\tilde{\bq}_{_{\rm L,2}}}
\newcommand{\bn}{\nu}
\newcommand{\bv}{\mathbf{v}}
\newcommand{\om}{\omega}
\newcommand{\pa}{\partial}
\newcommand{\la}{\langle}
\newcommand{\ra}{\rangle}
\newcommand{\lla}{\la{\hskip -2pt}\la}
\newcommand{\rra}{\ra{\hskip -2pt}\ra}
\newcommand{\jj}{\|{\hskip -0.8pt} |}
\newcommand{\al}{\alpha}
\newcommand{\ze}{\zeta}
\newcommand{\si}{\sigma}
\newcommand{\ep}{\varepsilon}
\newcommand{\na}{\nabla}
\newcommand{\vp}{\varphi}
\newcommand{\ga}{\gamma}
\newcommand{\Ga}{\Gamma}
\newcommand{\Om}{\Omega}
\newcommand{\de}{\delta}
\newcommand{\Th}{\Theta}
\newcommand{\De}{\Delta}
\newcommand{\Lam}{\Lambda}
\newcommand{\lam}{\lambda}
\newcommand{\tri}{\triangle}
\newcommand{\lj}{[{\hskip -2pt} [}
\newcommand{\rj}{]{\hskip -2pt} ]}
\newcommand{\bks}{\backslash}
\newcommand{\diam}{\mathrm{diam}}
\newcommand{\osc}{\mathrm{osc}}
\newcommand{\meas}{\mathrm{meas}}
\newcommand{\dist}{\mathrm{dist}}

\newcommand{\mL}{\mathscr{L}}
\newcommand{\cT}{{\cal T}}
\newcommand{\cM}{{\cal M}}
\newcommand{\cE}{{\cal E}}
\newcommand{\cL}{{\cal L}}
\newcommand{\cF}{{\cal F}}
\newcommand{\cB}{{\cal B}}
\newcommand{\PML}{{\rm PML}}
\newcommand{\FEM}{{\rm FEM}}
\newcommand{\rd}{\,\mathrm{d}}

\renewcommand{\i}{\mathbf{i}}
\renewcommand{\v}{\mathbf{v}}
\renewcommand{\u}{\mathbf{u}}
\renewcommand{\r}{\mathbf{r}}
\newcommand{\gR}{{\mathbb{R}}}
\newcommand{\Z}{{\mathbb{Z}}}
\newcommand{\C}{{\mathbb{C}}}
\newcommand{\I}{{\mathbb{I}}}
\renewcommand{\Re}{\mathrm{Re}\,}
\renewcommand{\Im}{\mathrm{Im}\,}
\renewcommand{\div}{\mathrm{div}}
\newcommand{\curl}{\mathrm{curl}}
\newcommand{\Curl}{\mathbf{curl}}
\newcommand{\pv}{\mathrm{p.v.}}

\newcommand{\Np}{\mathbb{N}_p}
\newcommand{\Ns}{\mathbb{N}_s}
\newcommand{\Tp}{\mathbb{T}_p}
\newcommand{\Ts}{\mathbb{T}_s}
\newcommand{\Na}{\mathbb{N}_\alpha}
\newcommand{\Nb}{\mathbb{N}_\beta}
\newcommand{\Ta}{\mathbb{T}_\alpha}
\newcommand{\Tb}{\mathbb{T}_\beta}
\newcommand{\GG}{\mathcal{G}}

\newcommand{\N}{\mathbb{N}}
\newcommand{\D}{\mathbb{D}}
\newcommand{\T}{\mathbb{T}}
\newcommand{\A}{\mathbb{A}}
\newcommand{\B}{\mathbb{B}}
\newcommand{\G}{\mathbb{G}}
\newcommand{\F}{\mathbb{F}}
\newcommand{\R}{\mathbb{R}}
\newcommand{\W}{\mathbb{W}}
\newcommand{\V}{\mathbb{V}}
\newcommand{\U}{\mathbb{U}}
\newcommand{\J}{\mathbb{J}}
\newcommand{\Zg}{\mathbb{Z}}
\newcommand{\Gtheta}{\mathbb{\Theta}}
\newcommand{\Gphi}{\mathbb{\Phi}}

%%%%%%%%%%%%%%%%%%%%%%%%%%%%%%%%%%%%%%%%%%%%%%%%%%%%%%%%%%%%%%%%%%%%
\newcommand{\be}{\begin{eqnarray}}
\newcommand{\ee}{\end{eqnarray}}
\newcommand{\ben}{\begin{eqnarray*}}
	\newcommand{\een}{\end{eqnarray*}}
\newcommand{\nn}{\nonumber}

\section{Introduction}\label{section1}

Inverse elastic wave scattering problems have considerable interests in diverse application fields including non-destructive testings, medical imaging, and seismic exploration. The purpose of this paper is to propose and study a direct imaging method to find
the shape and location of unknown obstacles embedded in the half-space isotropic and homogeneous elastic medium. We assume the obstacles are far away from the surface of the medium where the sources and receivers are located. The imaging method is
based on the idea of reverse time migration (RTM) and does not require the knowledge of 
physical properties of the obstacles such as penetrable or non-penetrable, and for non-penetrable obstacles, the type of boundary conditions on the boundary of the obstacle.

Let $D\subset\R_+^2=\{(x_1,x_2)^T\in\R^2:x_2>0\}$ be a bounded Lipschitz domain with the unit outer normal $\nu$ to its boundary $\Ga_D$. We
assume the incident wave is emitted by a point source at $x_s$ on the surface $\Ga_0=\{(x_1,x_2)^T\in\R^2:x_2=0\}$, along the polarization direction $q\in\R^2$. Let $\N(x,x_s)$ be the Neumann Green tensor for the half-space elastic scattering problem with free surface condition on $\Ga_0$ (see section 2 below).
The measured data is $u_q(x_r,x_s)=u^s_q(x_r,x_s)+\N(x_r,x_s)q, x_r\in\Ga_0$, where $u_q^s(x,x_s)$ satisfies the following equations
\be
& &\Delta_e u_q^s(x,x_s)+ \rho\,\omega^2u_q^s(x,x_s)= 0 \ \ \ \ \mbox{in }\R_+^2\bks \bar{D},\label{p1}\\
& &u^s_q(x,x_s)=-\N(x,x_s)q \ \ \mbox{on} \ \Ga_D,\ \ \sigma(u_q^s(x,x_s))e_2=0 \ \ \mbox{on} \ \Ga_0,\label{p2}
\ee
where $\Delta_e u:=(\lambda+\mu)\nabla\div u+\mu\Delta u$ is the linear elastic operator with {Lam\'{e}} constants $\lambda$ and $\mu$ satisfying $\lam>0,\mu>0$, $\rho$ is the density, $\omega>0$ is the circular frequency, and
$e_i$ is the unit vector along the $x_i$ axis, $i=1,2$. In the following, we will always assume $\rho=1$. In the boundary condition (\ref{p2}), $\sigma(u)\in\C^{2\times 2}$ is the stress tensor, which relates the strain tensor $\ep(u)\in\C^{2\times 2}$ through the following constitutive law
\ben
\sigma(u) = 2\mu\ep(u) + \lambda\div u\, \I, \ \ \ep(u)=\frac{1}{2}(\na u +(\na u)^T).
\een
Here $\I\in\R^{2\times 2}$ is the identity matrix. 

The equations (\ref{p1})-(\ref{p2}) must be complemented by appropriate boundary conditions at infinity to 
make the problem well-posed. In this paper we shall take the method of limiting absorption principle to define the scattering solution of the problem (\ref{p1})-(\ref{p2}). The limiting absorption principle defines
the scattering solution of (\ref{p1})-(\ref{p2}) as the limit of the solution of the same equations with the complex frequency $\omega(1+\i\eps)$ when $\eps\to 0^+$. The limiting absorption principle for the half-space elastic scattering problems is proved in \cite{Yves1988}, \cite{sini2004} for the scatterer with traction free boundary conditions. The results can be easily extended to cover penetrable scatterers or non-penetrable scatterers with Dirichlet or impedance boundary conditions. We also refer to \cite{nedelec2011}, \cite{Guzina2006} and
the references therein for the study of radiation conditions for half-space elastic scattering problems. 

The RTM method, whose imaging function is defined as the cross-correction between the incident wave field
and the back-propagated wave field using the complex conjugated data, is nowadays widely used in exploration geophysics \cite{claerbout1985imaging, berkhout2012seismic, bleistein2013mathematics}. In \cite{chen2013reverse_acou, chen2013reverse_elec, chen2015reverse_elas, RTMhalf_aco}, the RTM method for reconstructing extended targets using acoustic, electromagnetic and elastic waves at a fixed frequency in the free space is proposed and studied. The resolution
analysis in \cite{chen2013reverse_acou, chen2013reverse_elec, chen2015reverse_elas, RTMhalf_aco} is achieved without using the small inclusion or geometrical optics assumption previously made in the literature (e.g. \cite{ammari2013mathematical, bleistein2013mathematics}). 

In the geophysics literature, the RTM method for elastic scattering data usually consists of three steps:  1) back-propagating the received elastic data on the surface to the medium using the full elastic wave equation; 2) decomposing the back-propagated wave field to obtain the $p$-wave and $s$-wave components by Helmholtz decomposition; and 3) cross-correlating each decomposed mode component with the corresponding mode component of the incident field to output the imaging profile, see e.g.  \cite{sun2001, yan2008, denli2008, chung2012, li}. In this paper, we study the elastic wave RTM method without using the wave field separation, i.e., the imaging condition is defined as the cross-correlation between the back-propagated wave field with the incident wave field \cite{chang1987}. More precisely, we study the following imaging function (see (\ref{cor2}) below):
\ben
\hskip-1cm\hat{I}_d(z)=\Im\sum_{q=e_1,e_2}\int_{\Gamma_0^d}\int_{\Gamma_0^d}\,
[\T_D(x_s,z)^Tq][\T_D(x_r,z)^T\overline{u^s_q(x_r,x_s)}]\,ds(x_r)ds(x_s).
\een
Here $\Ga_0^d=\{x\in\Ga_0: x_1\in (-d,d)\}$, $d>0$, is the interval where the data are collected and $\T_D(x,z)$ is the traction tensor on $\Ga_0$ associated with the Dirichlet Green tensor (see (\ref{DGT2}) below).

Our resolution analysis, which extends the study in \cite{RTMhalf_aco} for the half-space acoustic scattering data, indicates that the imaging function always peaks on the illuminating part of the obstacle. The important Rayleigh surface wave is considered in our resolution analysis which shows that the contribution of the surface wave decays exponentially in the imaging function. The elastic wave RTM method based on the wave separation can be studied using the techniques developed in this paper and will be considered in a forthcoming work.

The layout of the paper is as follows. In section 2 we study the Neumann and Dirichlet Green tensors for the half-space elastic scattering problem by using the method of Fourier transform. In section 3 we study the point spread function defined by the RTM method. In section 4 we study the resolution of our RTM method for locating extended targets. In section 5 the extension of our resolution results to other types of obstacles are briefly considered. In section 6 we report extensive numerical results of our RTM method
for synthesized scattering data. In section 7 we prove a technical result used in the resolution analysis which is of independent interest.

\section{Elastic Green tensors in the half space}

In this section we introduce the elastic Green tensors and study their horizontal asymptotic behavior on the surface $\Ga_0$, which will play a crucial role in the resolution analysis for the RTM method to be proposed in this paper. Throughout the paper, we will assume that for $z\in\mathbb{C}\bks\{0\}$, $z^{1/2}$ is the analytic branch of $\sqrt{z}$ such that $\Im (z^{1/2})\geq0$. This corresponds to the right half real axis as the branch cut in the complex plane. For $z=z_1+\i z_2,z_1,z_2\in\mathbb{R}$, we have
\be \label{convention_1}
z^{1/2}={\rm sgn}(z_2)\sqrt{\frac{|z|+z_1}{2}}+\i\sqrt{\frac{|z|-z_1}{2}},\ \ \forall z\in\C\backslash\bar{\R}_+.
\ee
For $z$ on the upper and lower side of the right half real axis $\R_+=\{z\in\C:\Re z>0,\Im z=0\}$, we take $z^{1/2}$ as the limit of $(z\pm\i\ep)^{1/2}$ as $\ep \to 0^+$.

We start by introducing the Neumann Green tensor $\mathbb{N}(x,y)\in \C^{2\times 2}, y\in\R^2_+,$ which satisfies, for any $q\in\R^2$, 
\be
& & \De_e [\N(x;y)q] + \omega^2 [\N(x,y)q] = -\mathbf{\de}_y(x) q \ \ \mbox{in }\R^2_+ , \label{eq_n1} \\
& & \sigma(\N(x,y)q)e_2 = 0 \ \ \mbox{on } \Ga_0, \label{eq_n2}
\ee
where $\de_y(x)$ is the Dirac source at $y$. We use the method of Fourier transform 
to derive a formula of the Neumann Green tensor which is equivalent to that in \cite{nedelec2011} but is more convenient for our purpose. Let 
\be\label{a1}
\hat \N(\xi,x_2;y_2)= \int_\R\N(x_1,x_2;y) e^{-\i (x_1-y_1)\xi} dx_1,\ \ \forall \xi\in\C,
\ee
be the spectral Neumann Green tensor. Let $\G(x,y)$ be the fundamental solution tensor of the elastic equation \cite{ku63} whose Fourier transform is $\hat{\G}(\xi,x_2;y_2)=\hat{\G}_s(\xi,x_2;y_2)+\hat{\G}_p(\xi,x_2;y_2)$ with
\be
& &\hat{\G}_s(\xi,x_2;y_2)=\frac{\i}{2\omega^2}
\left( \begin{array}{cc}
	\mu_s & -\xi\frac{x_2-y_2}{|x_2-y_2|} \\
	-\xi\frac{x_2-y_2}{|x_2-y_2|} & \frac{\xi^2}{\mu_s}
\end{array} \right)e^{\i\mu_s|x_2-y_2|}, \label{G1}\\
& &\hat{\G}_p(\xi,x_2;y_2)=\frac{\i}{2\omega^2} 
\left( \begin{array}{cc}
	\frac{\xi^2}{\mu_p} & \xi\frac{x_2-y_2}{|x_2-y_2|} \\
	\xi\frac{x_2-y_2}{|x_2-y_2|} & \mu_p
\end{array} \right) e^{\i\mu_p|x_2-y_2|}.\label{G2}
\ee
Here $\mu_\alpha=(k_\alpha^2-\xi^2)^{1/2}$ for $\alpha=s,p$, $k_p=\omega/\sqrt{\lam+2\mu}, k_s=\omega/\sqrt{\mu}$ are the $p$ and $s$ wave numbers. Using the spectral fundamental solution tensor, one can 
write the spectral Neumann Green tensor as
\be\label{NGT}
\hspace{-2cm}\hat \N(\xi,x_2;y_2) = \hat \G(\xi,x_2;y_2)  -\hat \G(\xi,x_2;-y_2) + \frac{\i}{\omega^2\delta(\xi)}\sum_{\alpha,\beta=p,s}\mathbb{A}_{\al\beta}(\xi)e^{\i(\mu_\al x_2+\mu_\beta y_2)}, 
\ee
where $\beta(\xi)=k_s^2-2\xi^2$, $\delta(\xi)=\beta(\xi)^2+4\xi^2\mu_s\mu_p $, and
\ben
&&{\mathbb{A}_{ss}(\xi)} =
\left( \begin{array}{ll}
	\beta^2\mu_s & -4\xi^3\mu_s\mu_p \\
	-\xi\beta^2  & 4\xi^4\mu_p
\end{array} \right),\ \ 
{\mathbb{A}_{sp}(\xi)} =
\left( \begin{array}{ll}
	2\xi^2\beta\mu_s & -2\xi\beta\mu_s\mu_p \\
	-2\xi^3\beta  & 2\xi^2\beta\mu_p
\end{array} \right),\\
&&
{\mathbb{A}_{ps}(\xi)} =
\left( \begin{array}{ll}
	2\xi^2\beta\mu_s & 2\xi^3\beta \\
	2\xi\beta\mu_s\mu_p  & 2\xi^2\beta\mu_p
\end{array} \right),\ \ 
{\mathbb{A}_{pp}(\xi)} =
\left( \begin{array}{ll}
	4\xi^4\mu_s & \xi\beta^2 \\
	4\xi^3\mu_s\mu_p  & \beta^2\mu_p
\end{array} \right).
\een

The desired Neumann Green tensor should be obtained by taking the inverse Fourier transform of the spectral Green tensor $\hat{\N}(\xi,x_2;y_2)$. Unfortunately, one cannot simply take the inverse Fourier transform in the above formula because $\delta(\xi)$ have zeros in the real axis \cite{achenbach1980, Harris2001Linear}.

\begin{lem} \label{lem:2.1}
The Rayleigh equation $\delta(\xi) = 0$ has only two zeros $\pm k_R$, $k_R>k_s$, in the complex plane. 
\end{lem}

\debproof
For the sake of completeness, we include a proof here. By (\ref{convention_1}), It is clear that $\delta(\xi)$ is analytic outside the branch cuts $C_l=\{\xi=\xi_1+\i\xi_2\in\C: \xi_1\in [-k_s,-k_p],\xi_2=0\}$ and 
$C_r=\{\xi=\xi_1+\i\xi_2\in\C: \xi_1\in [k_p,k_s],\xi_2=0\}$. On the branch cuts,
\ben
\delta(\xi)=(k_s^2-2\xi^2)^2+\i\,[4\xi^2(k_s^2-\xi^2)^{1/2}(\xi^2-k_p^2)^{1/2}], \ \ \forall \xi\in C_l\cup C_r.
\een
Thus, $\de(\xi)$ has no zeros in $C_l\cup C_r$ and at least two real zeros $\pm k_R$, $k_R>k_s$, since $\de(\pm k_s)>0$, $\de(\pm\infty)<0$. The upper and lower sides of $C_l,C_r$ are denoted by $C_l^\pm,C_r^\pm$, respectively.

To conclude the proof, we now show that $\delta(\xi)$ has only two zeros in the complex plane by the principle of argument \cite{Ahlfors1979Complex}. Let $\Ga_R$ be a circle with sufficiently large radius $R$. We consider the domain $\mathcal D$ surrounded by the contour $C$ consisting of $\Ga_R$, $\Ga_l$ from $-k_s$ to $-k_p$ along $C_l^+$ and then from $-k_p$ to $-k_s$ along $C_l^-$, and $\Ga_r$ from $k_p$ to $k_s$ along $C_r^+$ and then from $k_s$ to $k_p$ along $C_r^-$. Since $\delta(\xi)$ has no poles in the complex plane,  we know from the principle of argument that the number of zeros 
in $\mathcal D$ is
\be\label{zero}
Z=\frac{1}{2\pi\i}\int_C \frac{\delta'(\xi)}{\delta(\xi)}d\xi.
\ee
It is clear that $\de(\xi)=\de^\pm(\xi)$ for $\xi\in C_r^\pm$, where
\ben
\de^\pm(\xi)=(k_s^2-2\xi^2)^2\mp\i\,[4\xi^2(k_s^2-\xi^2)^{1/2}(\xi^2-k_p^2)^{1/2}\,]:=f_1(\xi)\mp\i f_2(\xi).
\een
Then we have
\ben
\hskip-2cm\int_{\Ga_r} \frac{\delta'(\xi)}{\delta(\xi)}d\xi=\int_{k_p}^{k_s}\left(\frac{{\delta}_{+}' (\xi)}{\delta_{+}(\xi)}-\frac{{\delta}_{-}' (\xi)}{\delta_{-}(\xi)}\right) d\xi
&=&2\i\int_{k_p}^{k_s}\frac{f_1'(\xi) f_2(\xi)-f_1(\xi) f_2'(\xi)}{f_1^2(\xi)+ f_2^2(\xi)} d\xi\\
\hskip-2cm&=&-2\i\arctan \frac{f_2(\xi)}{f_1(\xi)}\Bigg|^{k_s}_{k_p}=0.
\een
Similarly, we have $\int_{\Ga_l}\frac{\delta'(\xi)}{\delta(\xi)}d\xi=0$. Moreover, for $|\xi|$ large, we have $\delta(\xi)=-2(k_p^2+3k_s^2)\xi^2+O(1)$, and consequently 
$\int_{\Ga_R} \frac{\delta'(\xi)}{\delta(\xi)}d\xi=4\pi$ for $R\gg 1$.
This yields $Z=2$ and completes the proof.
\finproof

Let $\mathbb{N}_{\omega(1+\i\ep)}(x,y)$ be the Neumann Green tensor with complex circular frequency $\om(1+\i\ep)$, that is, $\omega$ in (\ref{eq_n1}) is replaced by $\omega(1+\i\ep)$. Let $\N_{\omega(1+\i\ep)}(\xi,x_2;y_2)$ be the corresponding spectral Neumann Green tensor which are obtained by replacing $k_s, k_p$ in (\ref{NGT}) by 
$k_s(1+\i\ep), k_p(1+\i\ep)$, respectively. The Neumann Green tensor $\mathbb{N}(x,y)$ is defined by the limit absorption principle as
\be\label{NGT1}
\hspace{-1.5cm}\N(x,y)=\lim_{\ep\to 0^+} \N_{\om(1+\i\ep)}(x,y)=\lim_{\ep\to 0^+}\frac{1}{2\pi}\int_\R\hat \N_{\om(1+\i\ep)}(\xi,x_2;y_2) e^{\i(x_1-y_1)\xi} d\xi.
\ee
The above limit can be computed by the following lemma on the Cauchy principal value (cf. e.g. \cite[Chapter 4, Theorem 5]{Kuroda}). 

\begin{lem}\label{lem:2.2}
Let $a,b\in\R, a<b$, and $t_0\in (a,b)$. If $\gamma$ is H\"older continuous in $[a,b]$, that is, there exists a constant $\alpha\in (0,1]$ and a constant $C>0$ such that for any $s,t\in [a,b]$, $|\gamma(s)-\gamma(t)|\le C|s-t|^\alpha$, then
\ben
\lim_{z\to t_0,\pm\Im z>0}\int^b_a\frac{\gamma(t)}{t-z}dt={\rm p.v.}\int^b_a\frac{\gamma(t)}{t-t_0}dt\pm\pi\i\ga(t_0),
\een
where ${\rm p.v.}\int^b_a$ denotes the Cauchy principal value of the integral. 
\end{lem}

Lemma \ref{lem:2.2} and (\ref{NGT1}) yield the following representation formula for the Neumann Green function
\ben
\N(x,y)&=&\frac{1}{2\pi}\,{\rm p.v.}\int_{\R}\hat \N(\xi,x_2;y_2) e^{\i(x_1-y_1)\xi} d\xi\\
& &-\frac{1}{2\omega^2}
\left[\sum_{\alpha,\beta=p,s}\frac{\mathbb{A}_{\al\beta}(\xi)}{\de'(\xi)}e^{\i(\mu_\al x_2+\mu_\beta y_2)+\i(x_1-y_1)\xi}\right]^{k_R}_{-k_R},\ \ \forall x,y\in\R^2_+,
\een
where $[f(\xi)]^b_a:=f(b)-f(a)$.

For $x_s\in\Ga_0$, we define $\N(x,x_s),x\in\R^2_+$, as the limit of $\N(x,y)$ when $y\in\R^2_+, y\to x_s$.
It is also easy to check that $\N(x,y)=\N(y,x)^T$ for any $x,y\in\R^2_+$.

In the following we are mostly interested in the Neumann Green tensor $\N(x,y)$ when $x\in\Ga_0, y\in\R^2_+$. In this case, (\ref{NGT}) simplifies to
\be
\hspace{-2cm}
\hat
        \N(\xi,0;y_2)&=&\frac{\i}{\mu\delta(\xi)} \Bigg[ \Bigg(
		\begin{array}{cc}
			2\xi^2\mu_s & -2\xi\mu_s\mu_p\\
			-\xi\beta & \mu_p\beta
		\end{array} \Bigg)e^{\i\mu_p y_2}
		+ \Bigg(
		\begin{array}{cc}
			\mu_s\beta & \xi\beta \\
			2\xi\mu_s\mu_p & 2\xi^2\mu_p
		\end{array} \Bigg)e^{\i\mu_s y_2} \Bigg] \nonumber\\
	  &:=&\frac{1}{\delta(\xi)}(\Np(\xi)e^{\i\mu_p y_2}+\Ns(\xi)e^{\i\mu_s y_2}), \label{d2}
\ee
and consequently, for $x\in\Ga_0, y\in\R^2_+$,
\be\label{c8}\hspace{-2cm}
\N(x,y)=\frac{1}{2\pi}\,{\rm p.v.}\int_{\R}\hat \N(\xi,0;y_2) e^{\i(x_1-y_1)\xi} d\xi+\frac\i 2
\left[\sum_{\al=p,s}\frac{\Na(\xi)}{\de'(\xi)}e^{\i\mu_\al y_2+\i(x_1-y_1)\xi)}\right]^{k_R}_{-k_R}.
\ee

The following representation is useful in studying the horizontal asymptotic behavior of the Neumann Green tensor on $\Ga_0$.

\begin{figure}
	\centering
	\includegraphics[width=0.8\textwidth,height=0.5\textwidth]{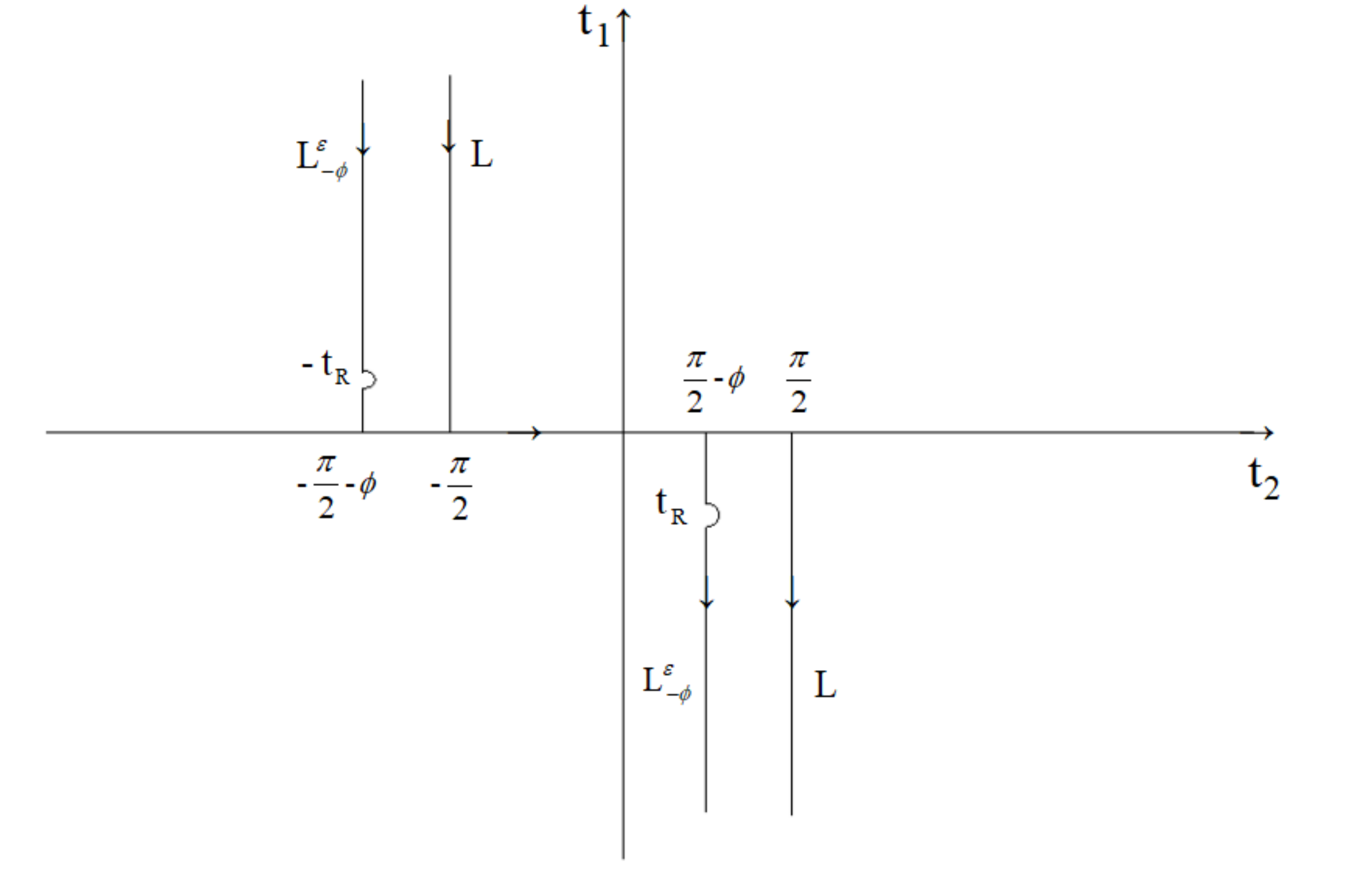}
	\caption{The integral path $L$ and $L^\eps_{-\phi}$.}\label{figure_trans}
\end{figure}

\begin{lem}\label{lem:2.3} Let $x\in\Ga_0, y\in \R^2_+$ and $\phi\in (-\pi/2,\pi/2)$ such that $y_2=|x-y|\cos\phi,
x_1-y_1=|x-y|\sin\phi$. Assume that $x_1\not= y_1$, then we have
\be\hskip-1.5cm
\N(x,y)&=&\frac{1}{2\pi}\int_L\mathbb{N}_0(t)\cos(t+\phi)e^{\i\lam\cos t}dt\pm\i
\left[\sum_{\al=p,s}\frac{\Na(\xi)}{\de'(\xi)}e^{\i\mu_\al y_2+\i(x_1-y_1)\xi}\right]_{\xi=\pm k_R},\label{h3}
\ee
where $\lam=k_s|x-y|$, $L$ is the integral path from $-\pi/2+\i\infty$ to $-\pi/2$, $-\pi/2$ to $\pi/2$, and $\pi/2$ to $\pi/2-\i\infty$ in the complex plan (see Figure \ref{figure_trans}), the sign $\pm$ is taken according to ${\rm sgn}(x_1-y_1)=\pm 1$, and
\be\label{h2}
\mathbb{N}_{0}(t)=\sum_{\al=p,s}k_s\,\frac{\Na(k_s\sin(t+\phi))}{\de(k_s(\sin(t+\phi))}.
\ee
\end{lem}

\debproof
Without loss of generality, we assume $x_1>y_1$ and so ${\rm sgn}(x_1-y_1)=1$. Notice that $\hat{\N}(\xi,0;y_2)=\sum_{\al=p,s}\frac{\Na(\xi)}{\de(\xi)}e^{\i (y_2\mu_\al+(x_1-y_1)\xi)}$. For $\al=p,s$, we use the classical transform $\xi=k_\al\sin t$ to obtain
\ben
\hskip-2cm\frac{1}{2\pi}\,{\rm p.v.}\int_{\R}\hat \N(\xi,0;y_2) e^{\i(x_1-y_1)\xi} d\xi
=\frac 1{2\pi}\,{\rm p.v.}\int_L\sum_{\al=p,s}k_s\,\frac{\Na(k_s\sin t)}{\de(k_s\sin t)}\cos t\, e^{\i \lam\cos (t-\phi)}dt.
\een
Let $L_{-\phi}$ be the integral path which is the shift of $L$ by $-\phi$, then
\be\label{h1}
\hskip-1cm\frac{1}{2\pi}\,{\rm p.v.}\int_{\R}\hat \N(\xi,0;y_2) e^{\i(x_1-y_1)\xi} d\xi
=\frac 1{2\pi}\,{\rm p.v.}\int_{L_{-\phi}}\mathbb{N}_0(t)\cos (t+\phi)\,e^{\i\lam\cos t}dt,
\ee
where $\mathbb{N}_0(t)$ is defined in (\ref{h1}). 

Let $t_R=\pi/2-\i s_R\in L$, $s_R>0$, such that $k_R=k_s\sin t_R$. Thus $k_R$ is the image of $t_R$ under the integral transform $\xi=k_s\sin t$.  For any $\eps>0$, let $L^\eps$ be the integral path from $-\pi/2+\i\infty\to -\pi/2+\i (s_R+\eps)\cup\pa B^+_\eps(-t_R)\to -\pi/2+\i(s_R-\eps)\to -\pi/2
\to\pi/2\to\pi/2-\i(s_R-\eps)\to\pa B^+_\eps(t_R)\to\pi/2-\i(s_R+\eps)\to\pi/2-\i\infty$, where $\pa B^+_\eps(\pm t_R)$ is the right half circle of radius $\eps$ centered at $\pm t_R$ (see Figure \ref{figure_trans}). Let $L^\eps_{-\phi}$ be the shift $L^\eps$ by $-\phi$. Then by the definition of Cauchy principle value, we know that
\ben
\hskip-1.5cm\frac 1{2\pi}\,{\rm p.v.}\int_{L_{-\phi}}\mathbb{N}_0(t)\cos(t+\phi)\,e^{\i \lam\cos t}dt
&=&\lim_{\eps\to 0^+}\frac 1{2\pi}\int_{L^\eps_{-\phi}}\mathbb{N}_0(t)\cos (t+\phi)\,e^{\i \lam\cos t}dt\\
\hskip-1.5cm& &+\frac\i 2\sum_{t'=\pm t_R}{\rm Res}(\mathbb{N}_0(t)\cos (t+\phi)e^{\i \lam\cos t},t').
\een
It is easy to see that the residue
\ben
\frac\i 2\sum_{t'=\pm t_R}{\rm Res}(\mathbb{N}_0(t)\cos (t+\phi)e^{\i \lam\cos t},t')=\frac \i 2\sum_{\xi=\pm k_R}\sum_{\al=p,s}\frac{\Na(\xi)}{\de'(\xi)}e^{\i (y_2\mu_\al+(x_1-y_1)\xi)}.
\een
On the other hand, by Cauchy integral theorem we have
\ben
\frac 1{2\pi}\int_{L^\eps_{-\phi}}\mathbb{N}_0(t)\cos (t+\phi)\,e^{\i \lam\cos t}dt
=\frac 1{2\pi}\int_{L}\mathbb{N}_0(t)\cos (t+\phi)\,e^{\i \lam\cos t}dt.
\een
This completes the proof of the lemma by (\ref{c8}) and (\ref{h1}).
\finproof

This lemma is the starting point of our estimate of the decay behavior of $\N(x,y)$ as $x\to\infty$ on $\Ga_0$.
We recall first the following Van der Corput lemma for the oscillatory integral \cite[P.152]{grafakos}.

\begin{lem}\label{van}
Let $\lam\ge 1$, $f\in C[a,b]$ with absolutely integrable derivative, and $u\in C^k[a,b]$, where $k\ge 1$ and $a<b$. \\
{\rm 1}. If $|u'(t)|\ge 1$ for $t\in (a,b)$ and $u'$ is monotone in $(a,b)$, then 
	\ben
	\left|\int^b_a f(t)e^{\i\lambda u(t)}dt\right|\le
	3\lambda^{-1}\left(|f(b)|+\int^b_a |f'(t)|dt\right).
	\een
{\rm 2}. For $k\geq2$, if $|u^{(k)}(t)|\ge 1$ for $t\in (a,b)$, then 
	\ben
	\left|\int^b_a f(t)e^{\i\lambda u(t)}dt\right|\le
	12k\lambda^{-1/k}\left(|f(b)|+\int^b_a|f'(t)|dt\right).
	\een
\end{lem}

The following lemma, which will be useful in the subsequent analysis, shows that the Van der Corput lemma is still valid when the singular points of the integrand $\phi(t)$ have a gap to the stationary phase points.

\begin{lem}\label{lem:2.5}
Let $\lam\ge 1$ and $f\in C[-\pi/2,\pi/2]$ having absolutely integrable derivative. Then for any $(a,b)\subset (-\pi/2,\pi/2)$, we have
 \be\label{c1}
   \left|\int_a^bf(t)e^{\i\lam\cos t}dt\right| 
   \leq C\lam^{-1/2}\left(|f(0)|+\int_a^b|f'(t)|dt\right),
   \ee
where the constant $C$ is independent of $a,b,\lam$ and the integrand $f$. 
Moreover, let $\kappa\in (0,1)$ and $\phi\in (-\pi/2,\pi/2)$ such that $|\phi|\geq\phi^*>\arcsin \kappa:=\phi_\kappa$, we have
   \be\label{c3}\hspace{-2cm}
   \left|\int_{-\frac\pi 2}^{\frac\pi 2}f(t)(\kappa^2-\sin^2(t+\phi))^{-1/2}e^{\i\lam\cos t}dt\right| 
   \leq C\lam^{-1/2}\left(|f(0)|+\int_{-\frac\pi 2}^{\frac \pi2}|f'(t)|dt\right),
   \ee
   where $C$ depends only on $\phi^*$ and $\kappa$.
\end{lem}
\debproof
The estimate (\ref{c1}) follows directly from Lemma \ref{van} since the interval $(a,b)$ can be divided into several subintervals so that in each subinterval, either $|\sin t|$ or $|\cos t|$ is bounded below by $1/\sqrt 2$. 

Let $g(t)=\kappa^2-\sin^2(t+\phi)$. It is easy to see that $g(t)$ has two zeros $t_1, t_2$ in $(-\pi/2,\pi/2]$, where
$t_1=\phi_\kappa-\phi$ and $t_2=-\phi_\kappa-\phi$ or $t_2=\pi-\phi_\kappa-\phi$ depending on whether $\phi+\phi_\kappa<\pi/2$ or $\phi+\phi_\kappa\ge \pi/2$. Without loss of generality, we assume the later case and thus $t_2=\pi-\phi_\kappa-\phi$. 

Let $\eps_0=\min(\frac{\phi^*-\phi_\kappa}{2},\frac{\phi_\kappa}{2})>0$. Obviously, $t_1-\eps_0\le -\pi/2, t_1+\eps_0\le -(\phi^*-\phi_\kappa)/2$ and $t_2-\eps_0\ge(\phi^*-\phi_\kappa)/2$. We divide $(-\pi/2,\pi/2)$ into five intervals:
$I_1=(-\pi/2,t_1-\eps_0), I_2=(t_1-\eps_0,t_1+\eps_0), I_3=(t_1+\eps_0,t_2-\eps_0), I_4=(t_2-\eps_0,t_2+\min(t_2+\eps_0,\pi/2))$ and $I_5=(\min(t_2+\eps_0,\pi/2),\pi/2)$.
By (\ref{c1}) we have
\be\label{c4}
\hspace{-2cm}  \left|\int_{I_1\cup I_3\cup I_5}f(t)(\kappa^2-\sin^2(t+\phi))^{-1/2}e^{\i\lam\cos t}dt\right| 
  \leq C\lam^{-1/2}\left(|f(0)|+\int_{-d}^0|f'(t)|dt\right),
\ee
where the constant $C$ depends only on $\phi^*$ and $\kappa$. 

Now we estimate the integral in $I_2, I_4$.  We first observe that $|\sin t|\ge \sin((\phi^*-\phi_\kappa)/2)$ in $I_2\cup I_4$. Moreover, $|g'(t)|=|\sin(2(t+\phi))|\ge \min(\sin\phi_\kappa,\sin(\phi^*+\phi_\kappa))$ in $I_2\cup I_4$. Let $\de\in (0,\eps_0)$ be sufficiently small. Since $g(t_j)=0, j=1,2$, by the mean value theorem, we have
\ben
\hspace{-1cm}|g(t)|\ge \min(\sin\phi_\kappa,\sin(\phi^*+\phi_\kappa))\de,\ \ \forall \de\le |t-t_j|\le \eps_0,j=1,2.
\een
By integration by parts we then obtain
\ben
\hskip-1cm& &\left|\int_{t_1-\eps_0}^{t_1-\de}f(t)g(t)^{-1/2}e^{\i\lam\cos t}dt\right| \le C\delta^{-1/2}\lam^{-1}\left(|f(0)|+\int_{-\frac\pi 2}^{\frac \pi 2}|f'(t)|dt\right).
\een
Similarly, 
\ben
\hspace{-1.cm}\left|\int_{t_1+\de}^{t_1+\eps_0}f(t)g(t)^{-1/2}e^{\i\lam\cos t}dt\right| 
\le C\delta^{-1/2}\lam^{-1}\left(|f(0)|+\int_{-\frac\pi 2}^{\frac \pi 2}|f'(t)|dt\right).
\een
Finally, 
\ben
\hspace{-1.cm}\left|\int_{t_1-\delta}^{t_1+\de}f(t)g(t)^{-1/2}e^{\i\lam\cos t}dt\right| 
&\leq&C\max_{t\in(-\pi/2,\pi/2)}|f(t)|\int_{-\delta}^{\de}|\kappa -\sin(\phi_\kappa+t)|^{-1/2}dt\\
\hspace{-1.cm}&\leq&C\de^{1/2}\max_{t\in(-\pi/2,\pi/2)}|f(t)|.
\een
In conclusion, by taking $\delta=\lam^{-1}$, we obtain
\ben
\hspace{-1.5cm}\left|\int_{I_2}f(t)(\kappa^2-\sin^2(t+\phi))^{-1/2}e^{\i\lam\cos t}dt\right| 
\leq C\lam^{-1/2}\left(|f(0)|+\int_{-\frac\pi 2}^{\frac \pi 2}|f'(t)|dt\right).
\een
The integral in $I_4$ can be estimated similarly. This completes the proof by combining with the estimate in (\ref{c4}).
\finproof

The following lemma collects some facts about the Rayleigh function $\de(\xi)$.

\begin{lem}\label{delta}
Let $d_R=(k_R-k_s)/2$, where $k_R$ is the wave number of the Rayleigh surface wave. There exist constants $C_1,C_2$ depending only on $\kappa$ such that $|\delta^{(k)}(\xi)|\le C_2k_s^{4-k}, k=0,1,2$, $|\delta(\xi)|\ge C_1k_s^4$ for any $|\xi|\le k_R+d_R$, and $|\de(\xi)|\ge 2k_s^2(|\xi|^2-k_R^2)$ for any $|\xi|\ge k_R$. Moreover, let $\de(\xi)=\de_1(\xi)(\xi^2-k_R^2)$ for $\xi\in\R$, then $|\de_1(\xi)|\ge C_1k_s^2$ for any $k_R-d_R\le |\xi|\le k_R+d_R$.
\end{lem}

\debproof
By definition, we know that for $|\xi|\ge k_s$, $\de(\xi)=k_s^4f({\xi^2}/{k_s^2})$, where
\ben
f(t)=(2t-1)^2-4t\sqrt{t-1}\sqrt{t-\kappa^2},\ \ \forall t\ge 1.
\een
It is easy to see that $f'(t)\le -2$ for $t\ge 1$, $f(1)=1$ and $f((2-\kappa^2)/(1-\kappa^2))<0$. Thus $k_R^2\le\frac{2-\kappa^2}{1-\kappa^2}k_s^2$. The estimates for $\de(\xi)$ follows easily. 

Next by the mean value theorem, 
\ben 
\min_{k_R-d_R\le|\xi|\le k_R+d_R}|\de_1(\xi)|\ge\min_{k_R-d_R\le|\xi|\le k_R+d_R}\frac{|\de'(\xi)|}{|\xi|+k_R}\ge C_1k_s^2,
\een
where we have used the fact $f'(t)\le -2$ for $t\ge 1$. This completes the proof.
\finproof

\begin{lem}\label{lem:2.7}
Let $\phi\in (0,\pi/2)$ and $H$ be the hyperbola $\{\xi=\xi_1+\i\xi_2\in\C:(\xi_1/(k_s\cos\phi))^2-(\xi_2/(k_s\sin\phi))^2=1\}$. Let $f(\xi)$ be analytic on $H$. Then there exists a constant $C$ depending only $\kappa$ such that
\ben
\hskip-2.2cm& &\left|\int_{L\bks [-\pi/2,\pi/2]}f(k_s\sin(t+\phi))e^{\i\lam\cos t}dt\right|+\left|\int_{L\bks [-\pi/2,\pi/2]}f(k_s\sin(t+\phi))\cos te^{\i\lam\cos t}dt\right|\\
\hskip-2cm&\le&C\lam^{-1}\max_{\xi\in H}(|f(\xi)|+k_s|f'(\xi)|).
\een 
\end{lem}
\debproof
Notice that for $t=-\pi/2+\i s$, $s>0$, $k_s\sin(t+\phi)=-\cosh(s)\cos\phi+\i\sinh(s)\sin\phi\in H$. Thus
\ben
& &\int_{-\pi/2}^{-\pi/2+\i\infty}f(k_s\sin(t+\phi))e^{\i\lam\cos t}dt\\
&=&\i\int^\infty_0f(-\cosh(s)\cos\phi+\i\sinh(s)\sin\phi)e^{-\lam\sinh(s)}ds.
\een
This implies the estimate on the integral path $-\pi/2+\i\infty\to -\pi/2$ by integration by parts. The estimate on the path $\pi/2\to\pi/2-\i\infty$ is similar. Thus the proof of the first term in the lemma follows. The second term can be proved similarly. Here we omit the details. This completes the proof.
\finproof

The following theorem is the first main result in this section.

\begin{thm}\label{thm:2.1}
Let $x\in\Gamma_0$, $y\in\R_+^2$ satisfy $|x_1-y_1|/|x-y|\ge(1+\kappa)/2$ and $k_sy_2\ge 1$. There exists a constant $C$ depending only on $\kappa$ such that
	\ben
	|\N(x,y)|+k_s^{-1}|\na_y\N(x,y)|\leq \frac{C}{\mu}\left(\frac{k_sy_2}{(k_s|x-y|)^{3/2}}+e^{-\sqrt{k_R^2-k_s^2}y_2}\right).
	\een
\end{thm}

\debproof We only prove the estimate for $\N(x,y)$. The estimate of $|\na_y\N(x,y)|$ can be proved similarly. The starting point is (\ref{h3}) in Lemma \ref{lem:2.3}. Without loss of generality, we assume $x_1>y_1$ and thus $\phi\in (0,\pi/2)$ which satisfies $\phi\ge\phi^*=\arcsin (1+\kappa)/2>\phi_\kappa$. It is easy to see from Lemma \ref{delta} that the second term in (\ref{h3}) is bounded by $C\mu^{-1}e^{-\sqrt{k_R^2-k_s^2}y_2}$. 

For the first term in (\ref{h3}) we first note that
\ben
\hskip-2cm\frac 1{2\pi}\int_L\mathbb{N}_0(t)\cos(t+\phi)e^{\i\lam\cos t}dt=
\frac 1{2\pi}\int_L\sum_{\al=p,s} k_s\frac{\Na(k_s\sin(t+\phi))}{\de(k_s\sin(t+\phi))}\cos(t+\phi)e^{\i\lam\cos t}dt.
\een
We shall only estimate the term including $[\Np(k_s\sin(t+\phi))]_{22}=\mu^{-1}(\beta\mu_p)(k_s\sin(t+\phi))$. The other terms can be proved similarly. Thus denote by
\ben
\hskip-2cmg(t)=k_s\frac{[\Np(k_s\sin(t+\phi))]_{22}}{\delta(k_s\sin(t+\phi))}:=f(t)(\kappa^2-\sin^2(t+\phi))^{1/2},\ \ f(t)=\frac{k_s^2}{\mu}\frac{\beta(k_s\sin(t+\phi))}{\de(k_s\sin(t+\phi))}.
\een
Note first by integration by parts that
\ben
\hskip-2cm\int_{L}g(t)\cos(t+\phi)e^{\i\lam\cos t}dt&=&\cos\phi\int_Lg(t)\cos te^{\i\lam\cos t}dt-\sin\phi\int_{L}g(t)\sin te^{\i\lam\cos t}dt\\
&=&\cos\phi\int_Lg(t)\cos te^{\i\lam\cos t}dt-\frac{\sin\phi}{\i\lam}\int_{L}g'(t)e^{\i\lam\cos t}dt\\
&=&{\rm I}_1+{\rm I}_2.
\een
By Lemma \ref{delta} and (\ref{c1}) in Lemma \ref{lem:2.5}, we know that
\ben
\hskip-1cm\left|\int_{-\frac\pi 2}^{\frac \pi 2}g(t)\cos te^{\i\lam\cos t}dt\right|\le C\lam^{-1/2}\left(|g(0)|+\int^{\frac\pi 2}_{-\frac\pi 2}|(g(t)\cos t)'|dt\right)\le C\mu^{-1}\lam^{-1/2}.
\een
Since $g(t)=k_s([\Np(k_s\sin(t+\phi))]_{22}\de^{-1}(k_s\sin(t+\phi))$, $|[\Np(\xi)]_{22}|\le C|\xi|^3$, $|[\Np'(\xi)]_{22}|\le C|\xi|^2$, and $\de(\xi)\ge Ck_s^2|\xi|^2$ on the hyperbola $H$, we conclude from Lemma \ref{lem:2.7} that
\ben
\left|\int_{L\bks [-\frac\pi 2,\frac\pi 2]}g(t)\cos te^{\i\lam\cos t}dt\right|\le C\mu^{-1}\lam^{-1}.
\een
Thus $|{\rm I}_1|\le C\mu^{-1}\lam^{-1/2}\cos\phi$. Similarly, we can obtain $|{\rm I}_2|\le C\mu^{-1}\lam^{-3/2}$. Indeed, the only difference is that since $g'(t)=f'(t)(\kappa^2-\sin^2(t+\phi))^{1/2}-f(t)(\kappa^2-\sin^2(t+\phi))^{-1/2}\sin(t+\phi)\cos(t+\phi)$, one has to use both (\ref{c1}) and (\ref{c3}) in Lemma \ref{lem:2.5} to obtain
\ben
\left|\int_{-\frac\pi 2}^{\frac \pi 2}g'(t)e^{\i\lam\cos t}dt\right|\le  C\mu^{-1}\lam^{-1/2}.
\een
Thus $|{\rm I}_1+{\rm I}_2|\le C\mu^{-1}\lam^{-1/2}\cos\phi+C\mu^{-1}\lam^{-3/2}\le C\mu^{-1}(k_sy_2)/(k_s|x-y|)^{3/2}$, where have used the condition $k_sy_2\ge 1$. This completes the proof.
\finproof

Now we introduce the Dirichlet Green tensor $\D(x,y), y\in\R^2_+$ which satisfies \cite{arens1999}
\be
& & \De_e [\D(x,y)q] + \omega^2 [\D(x,y)q] = -\mathbf{\de}_y(x)q \ \ \mbox{in } \R^2_+ , \label{eq_d1} \\
& &  \D(x,y)q = 0 \ \ \mbox{on } \Ga_0. \label{eq_d2}
\ee 
The spectral Dirichlet Green tensor $\hat{\D}(\xi,x_2;y_2)$ is defined similar to the spectral Neumann Green tensor in (\ref{a1})
and it follows that  
\be
\hspace{-2cm}\hat \D(\xi,x_2;y_2) = \hat \G(\xi,x_2;y_2)  -\hat \G(\xi,x_2;-y_2) + \frac{\i}{\omega^2 \gamma(\xi)}\sum_{\al,\beta=s,p}\mathbb{B}_{\al\beta}(\xi)e^{\i(x_2\mu_\alpha+y_2\mu_\beta)},\label{DGT}
\ee
where $\gamma(\xi)=\xi^2+\mu_s\mu_p$, $\mathbb{B}_{sp}(\xi)=-\mathbb{B}_{ss}(\xi), \mathbb{B}_{ps}(\xi)=-\mathbb{B}_{pp}(\xi)$, and
\ben
&&{\mathbb{B}_{ss}(\xi)} =
\left( \begin{array}{ll}
	\xi^2\mu_s & -\xi\mu_s\mu_p \\
	-\xi^3  & \xi^2\mu_p
\end{array} \right),\ \ \ \ \ \
{\mathbb{B}_{pp}(\xi)} =
\left( \begin{array}{ll}
	\xi^2\mu_s & \xi^3 \\
	\xi\mu_s\mu_p  & \xi^2\mu_p
\end{array} \right).
\een
The Dirichlet Green tensor $\D(x,y)$ is obtained as
the limit of $\D_{\om(1+\i\ep)}(x,y)$ when $\ep\to 0^+$, where $\D_{\om(1+\i\ep)}(x,y)$ is Dirichlet Green tensor with
the complex circular frequency $\om(1+\i\ep)$, that is $\om$ in (\ref{eq_d1}) is replaced by $\om(1+\i\ep)$.
The corresponding spectral Dirichle Green tensor $\hat\D_{\om(1+\i\ep)}(x,y)$ are obtained by replacing $k_s, k_p$ in (\ref{DGT}) by $k_s(1+\i\ep), k_p(1+\i\ep)$, respectively. Thus 
\ben
\hspace{-1cm}\D(x,y)=\lim_{\ep\to 0^+} \D_{\om(1+\i\ep)}(x,y)=\lim_{\ep\to 0^+}\frac{1}{2\pi}\int_\R\hat \D_{\om(1+\i\ep)}(\xi,x_2;y_2) e^{\i(x_1-y_1)\xi} d\xi.
\een
We have the following representation of Dirichlet Green tensor
\be\label{DGT1}
\hskip-2cm\D(x,y)=\G(x,y)-\G(x,y')+\frac{\i}{2\pi\omega^2}\int_{\R}
\sum_{\al,\beta=s,p}\frac{\mathbb{B}_{\al\beta}(\xi)}{\ga(\xi)}e^{\i(\mu_\alpha x_2+\mu_\beta y_2)+\i(x_1-y_1)\xi}d\xi.
\ee
It is easy to check that $\D(x,y)=\D(y,x)^T$ for $x,y\in\R^2_+$.

For $x\in\Ga_0, y\in\R^2_+$, let $\T_D(x,y)\in\C^{2\times 2}$ denote the traction tensor of $\D(x,y)$ in the direction $e_2$ with respect to $x$, that is,
$\T_D(x,y)q=\sigma(\D(x,y)q)e_2, \forall q\in\R^2$. By  (\ref{DGT1}) we have
\be\label{DGT2}
\T_D(x,y)=\frac{1}{2\pi}\int_{\R}\hat \T_D(\xi,0;y_2) e^{\i(x_1-y_1)\xi}d\xi,\ \ \ \ \forall x\in\Ga_0,
\ee
where 
\be\hspace{-1cm} 
\hat\T_D(\xi,0;y_2)&=&\frac 1{\gamma(\xi)}\left[\left(   \begin{array}{cc}
	\xi^2 & -\xi\mu_p \\
	-\xi\mu_s & \mu_p\mu_s
\end{array} \right)e^{\i\mu_p y_2}+
\left(   \begin{array}{cc}
	\mu_s\mu_p & \xi\mu_p \\
	\xi\mu_s & \xi^2
\end{array}\right)e^{\i\mu_s y_2}\right]\nonumber\\
&:=&\Tp(\xi)e^{\i\mu_p y_2}+\Ts(\xi)e^{\i\mu_s y_2}.\label{d1}
\ee

The following theorem on the decay behavior of the traction tensor on $\Ga_0$ improves \cite[Lemma 2.2]{arens1999} in the sense that we provide exact dependence on $y_2$ in the numerator. It can be proved by the same (and simpler) argument as that in the proof of Theorem \ref{thm:2.1}. We omit the details.

\begin{thm}\label{thm:2.2}
	Let $x\in\Gamma_0$, $y\in\R_+^2$ satisfy $|x_1-y_1|/|x-y|\ge (1+\kappa)/2$ and $k_s y_2\ge 1$. There exists a constant $C$ depending only on $\kappa$ such that
\ben
|\T_D(x,y)|+k_s^{-1}|\na_y\T_D(x,y)|\leq C\frac{k_s^2 y_2}{(k_s|x-y|)^{3/2}}.
\een
\end{thm}

\section{The point spread function}

In this section we introduce the point spread function for imaging a point source embedded in the 
half-space elastic medium, which extends the study in \cite{RTMhalf_aco} for acoustic waves.   
Let $\N(x,y)$
be the Neumann Green tensor which is the data collected on the surface $\Ga_0^d=\{(x_1,x_2)^T\in\Ga_0: x_1\in (-d,d)\}$ with a point source $y\in\R^2_+$, where $d>0$ is the aperture. The finite aperture point spread function $\J_d(x,y)$, $x,y\in\R^2_+$, is a $\C^{2\times 2}$ matrix, which is
the back-propagated field with $\N(x,y)\chi_{(-d,d)}$ as the Dirichlet boundary condition, where $\chi_{(-d,d)}$ is the characteristic function of the interval $(-d,d)$. More precisely, $\J_d(x,y)e_j, j=1,2,$ is the scattering solution of the following problem
\ben
& &\De_e[\J_d(x,y)e_j]+\om^2[\J_d(x,y)e_j]=0\ \ \mbox{in }\R^2_+,\\
& &\J_d(x,y)e_j=[\overline{\N(x,y)}e_j]\chi_{(-d,d)}\ \ \mbox{on }\Ga_0.
\een
By the integral representation formula, for any $z,y\in\R^2_+$,
\ben
\hskip-1.5cm[\J_d(z,y)]_{ij}=e_i\cdot[\J_d(z,y)e_j]=\int_{\Ga_0^d}\T_D(x,z)e_i\cdot\overline{\N(x,y)}e_jds(x),\ \ i,j=1,2,
\een
or, more concisely,
\be\label{jd}
\J_d(z,y)=\int_{\Ga_0^d}\T_D(x,z)^T\overline{\N(x,y)}ds(x).
\ee
By Theorem \ref{thm:2.1} and Theorem \ref{thm:2.2}, we know that the integral converges as $d\to\infty$. Thus we can define the half-space elastic point spread function $\J(x,y)\in \C^{2\times 2}$, $x,y\in\R^2_+$, as
\be\label{j}
\J(z,y)=\int_{\Ga_0}\T_D(x,z)^T\overline{\N(x,y)}ds(x).
\ee
By the limiting absorption principle, we know that
\ben
\J(z,y)=\lim_{\ep\to 0^+}\int_{\Ga_0} \T_D^{\,\omega(1+\i\eps)}(x,z)^T\,
\overline{\N_{\om(1+\i\ep)}(x,y)}ds(x),
\een
where $\T_D^{\,\omega(1+\i\eps)}(x,z)q=\sigma(\D_{\om(1+\i\omega)}(x,z)q)e_2,\forall q\in\R^2$.
By using Parserval identity, Lemma \ref{lem:2.2}, (\ref{d2}) and (\ref{d1}), we obtain
\be
\J(z,y)&=&\frac{1}{2\pi}\sum_{\al,\beta=p,s}{\rm p.v.}\int_{\R}\frac{{\Ta}(\xi)^T \overline{\Nb(\xi)}}{\overline{\delta(\xi)}} e^{\i (\mu_\alpha z_2-\overline{\mu}_\beta y_2)+\i(y_1-z_1)\xi} d\xi \nn \\
& &-\frac\i 2\sum_{\al,\beta=p,s}\left[\frac{{\Ta}(\xi)^T \overline{\Nb(\xi)}}{\overline{\delta'(\xi)}} e^{\i (\mu_\alpha z_2-\overline{\mu}_\beta y_2)+\i(y_1-z_1)\xi}\right]^{k_R}_{-k_R}. \label{d3}
\ee
To proceed, we define
\be
\F(z,y)&=&\frac{1}{2\pi}\int^{k_p}_{-k_p} \  \frac{{\Tp}(\xi)^T \overline{\Np(\xi)}}{\overline{\delta(\xi)}} e^{\i \mu_p (z_2- y_2)+\i(y_1-z_1)\xi} d\xi\nn \\
&&+\frac{1}{2\pi}\int^{k_s}_{-k_s} \  \frac{{\Ts}(\xi)^T \overline{\Ns(\xi)}}{\overline{\delta(\xi)}} e^{\i \mu_s (z_2- y_2)+\i(y_1-z_1)\xi} d\xi. \label{d4}
\ee

Let $\Omega$ be the imaging domain and $h=\dist(\Omega,\Gamma_0)$ be the distance between $\Omega$ and $\Gamma_0$. We assume  there exist constants $0<c_1<1,c_2>0$ such that
\be\label{d0}
|x_1|\leq c_1 d , \ \ |x-y|\leq c_2 h ,\ \ \ \ \forall x,y \in \Omega.
\ee
We remark that this assumption is rather mild in practical applications.
The aim of this section is to show that for $z,y\in\Om$, $\F(z,y)$ is the main contribution in $\J_d(z,y)$. Moreover, $\F(z,y)$ decays as $|z-y|\to\infty$ and the imaginary part of the function $|\Im\F_{ii}(z,y)|, i=1,2$, peaks when $z=y$.

We start with the following lemma.

\begin{lem} \label{lem:3.1}
Let $k_s h\geq 1$ and $d\gg h$. For any $z,y\in\Omega$, we have
\ben
& &|\J(z,y)-\J_d(z,y)|+k_s^{-1}|\nabla_y(\J(z,y)-\J_d(z,y))| \\
&\leq&\frac{C}{\mu} \left[\left(\frac{h}{d}\right)^{2}+(k_s h)^{1/2}e^{-k_s h\sqrt{\kappa_R^2-1}}\left(\frac{h}{d}\right)^{1/2}\right],
\een
where the constant C depends only on $\kappa$.
\end{lem}
\debproof
By Theorem \ref{thm:2.1} and Theorem \ref{thm:2.2},  when $k_s h\geq 1$ and $d\gg h$, we have
\ben
& &\left| \int_{d}^{\infty} \left[\T_D(x,z)^T\overline{\N(x,y)}\right]_{x_2=0}dx_1
\right| \\
&\leq&
\frac{C}{\mu}\int_{d}^{\infty}\frac{k_s^{1/2} z_2}{|x-z|^{3/2}}\left(
\frac{k_s^{-1/2} y_2}{|x-y|^{3/2}}+e^{-\sqrt{k_R^2-k_s^2}y_2}\right) dx_1\\
&\leq&
\frac{C}{\mu}\int_{(1-c_1)d/h}^{\infty}\left(\frac{1}{(1+t^2)^{3/2}}+\frac{(k_s h)^{1/2}}{(1+t^2)^{3/4}} e^{-\sqrt{k_R^2-k_s^2}h}\right)  dt\\
&\leq&\frac{C}{\mu} \left[\left(\frac{h}{d}\right)^{2}+\frac{(k_s h)^{1/2}}{ e^{\sqrt{k_R^2-k_s^2}h}}\left(\frac{h}{d}\right)^{1/2}\right].
\een
Here we have used the first inequality in (\ref{d0}). Similarly, we can prove that the estimate for the integral in $(-\infty,-d)$. This shows the estimate for $\J(z,y)-\J_d(z,y)$. The estimate for $\nabla_y(\J(z,y)-\J_d(z,y))$ can be proved similarly.
\finproof

The following lemma shows the second term on the right-hand side of (\ref{d3}) is small.

\begin{lem}\label{lem:3.2}
There exists a constant $C$ depending only on $\kappa$ such that for any $z,y\in\Om$,
\ben
\left|\sum_{\al,\beta=p,s}\left[\frac{{\Ta}(\xi)^T \overline{\Nb(\xi)}}{\overline{\delta'(\xi)}} e^{\i (\mu_\alpha z_2-\overline{\mu}_\beta y_2)+\i(y_1-z_1)\xi}\right]^{k_R}_{-k_R}\right|\le \frac C\mu e^{-\sqrt{k_R^2-k_s^2}h}.
\een
\end{lem}

\debproof
We know that from (\ref{d1}), (\ref{d2}) that for $\al=p,s$, $|\Ta(\pm k_R)|\le Ck_R^2/k_s^2\le C$, $|\Na(\pm k_R)|\le Ck_R^3$. The lemma now follows easily by using Lemma \ref{delta}. 
\finproof

\begin{lem}\label{lem:3.3}
Let $g\in C^1(\R)\cap L^1(\R)$. There exists a constant $C$ depending only on $\kappa$ such that for any $z,y\in\Om$, 
\ben
\hskip-2cm\left|{\rm p.v.}\int_{|\xi|>k_s}\frac{g(\xi)}{\delta(\xi)}d\xi\right|\leq Ck_s^{-4}\int_{|\xi|>k_s}|g(\xi)|d\xi+
Ck_s^{-3}\max_{\xi\in(k_R-d_R,k_R+d_R)}(|g(\xi)|+k_s|g'(\xi)|).
\een
where $d_R =(k_R-k_s)/2$.
\end{lem}
\debproof
Without loss of generality, we prove the lemma for the integral in $(k_s,\infty)$. We write $\de(\xi)=(\xi^2-k_R^2)\de_1(\xi)$ as in Lemma \ref{delta}, where $\de_1(\xi)\not=0$ for $\xi>k_s$. By the definition of
the Cauchy principle value
\be\label{l4}
\hskip-2cm\pv\int^\infty_{k_s}\frac{g(\xi)}{\delta(\xi)}d\xi&=&\int_{k_s}^{k_R-d_R}\frac{g(\xi)}{\delta(\xi)}d\xi+
\int^\infty_{k_R+d_R}\frac{g(\xi)}{\delta(\xi)}d\xi\nn\\
\hskip-2cm&&+\int^{k_R+d_R}_{k_R-d_R}\frac{g(\xi)((\xi+k_R)\de_1(\xi))^{-1}-g(k_R)(2k_R\de_1(k_R))^{-1}}{(\xi-k_R)}d\xi.
\ee
By Lemma \ref{delta} we obtain easily
\ben
\left|\int_{k_s}^{k_R-d_R}\frac{g(\xi)}{\delta(\xi)}d\xi+
\int^\infty_{k_R+d_R}\frac{g(\xi)}{\delta(\xi)}d\xi\right|\le Ck_s^{-4}\int^\infty_{k_s}|g(\xi)|d\xi.
\een
The last term in (\ref{l4}) can be proved by using the mean value theorem and the bounds for $\de(\xi),\de_1(\xi)$ in Lemma \ref{delta}. This completes the proof.
\finproof

\begin{lem}\label{lem:3.4}
Let $k_sh\ge 1$. There exists a constant $C$ depending only on $\kappa$ such that for any $z,y\in\Om$, 
\ben
\left|\sum_{\al,\beta=p,s}{\rm p.v.}\int_{|\xi|>k_s}\frac{{\Ta}(\xi)^T \overline{\Nb(\xi)}}{\overline{\delta(\xi)}} e^{\i (\mu_\alpha z_2-\overline{\mu}_\beta y_2)+\i(y_1-z_1)\xi} d\xi\right|\le \frac C\mu(k_sh)^{-1}.
\een
\end{lem}
\debproof
For $\al,\beta=p,s$, we denote $g_{\al\beta}(\xi)=\Ta(\xi)^T\overline{\Nb(\xi)}e^{\i (\mu_\alpha z_2-\overline{\mu}_\beta y_2)+\i(y_1-z_1)\xi}$. By using Lemma \ref{lem:3.3}, we obtain easily
\ben
\hskip-2cm\left|\pv\int_{|\xi|>k_s}\frac{g_{\al\beta}(\xi)}{\overline{\de(\xi)}}d\xi\right|
&\le&\frac {C}{k_s^6\mu}\int^\infty_{k_s}|\xi|^5e^{-\sqrt{\xi^2-k_s^2}(y_2+z_2)}d\xi+\frac C\mu(k_sh) e^{-\sqrt{(k_R-d_R)^2-k_s^2}(y_2+z_2)}\\
&\le&\frac C\mu\int^\infty_1t^5e^{-\sqrt{t^2-1}(y_2+z_2)}dt+\frac C\mu (k_sh) e^{-\sqrt{(k_R-d_R)^2-k_s^2}(y_2+z_2)}\\
&\le&\frac C\mu (k_sh)^{-1}+\frac C\mu (k_sh) e^{-\sqrt{(k_R-d_R)^2-k_s^2}(y_2+z_2)},
\een
where we have used that $y_2,z_2\ge h$ and $d_R=(k_R-k_s)/2\ge C_1k_s$ for some constant $C_1>0$ depending only on $\kappa$.
This completes the proof as the second term decays exponentially in $k_sh$.
\finproof

\begin{lem}\label{lem:3.5}
Let $\phi(t)=\sqrt{1-t^2}-\tau\sqrt{\kappa^2-t^2}+\nu t$, where $\kappa\in (0,1), \tau\ge\tau_0>0, \nu\in\R$.
There exists a constant $C$ depending only on $\kappa, \tau_0$ but independent of $\nu$ such that for any $\lam>0$ and $f\in C[0,\kappa]$ with absolutely integrable derivative,
\ben\hskip-1cm
\left|\int^\kappa_{-\kappa}f(t)e^{\i\lam\phi(t)}dt\right|+\left|\int^\kappa_{-\kappa}f(t)e^{-\i\lam\phi(t)}dt\right|\leq C\lambda^{-1/4}\left(|f(0)|+\int_{-\kappa}^{\kappa}|f'(t)|dt\right).
\een
\end{lem}
\debproof
We only prove the estimate for the first integral in the interval $(0,\kappa)$. The other cases can be proved similarly. It is easy to check that for $t\in (0,\kappa), m\ge 2$, the $m$-th derivative $\phi^{(m)}(t)=\tau\kappa^{-(m-1)}\psi_m(t/\kappa)-\psi_m(t)$, where
\ben\hskip-1cm
\psi_2(t)=(1-t^2)^{-3/2},\ \ \psi_3(t)=3t(1-t^2)^{-5/2},\ \ 
\psi_4(t)=3(1+4t^2)(1-t^2)^{-7/2}.
\een
Obviously, $\psi_m(t),m\ge 2$, are increasing functions in $(0,\kappa)$. 

We first consider the case when $\tau\ge \kappa^2$. This implies $\tau\kappa^{-3}\ge\kappa^{-1}$ and thus
\ben
\phi^{(4)}(t)\ge(\kappa^{-1}-1)\psi_4(t)\ge 3(\kappa^{-1}-1).
\een
By using the Van der Corput Lemma \ref{van}, we have
\be\label{k2}
\left|\int^\kappa_{0}f(t)e^{\i\lam\phi(t)}dt\right|\leq C\lambda^{-1/4}\left(|f(0)|+\int_{-\kappa}^{\kappa}|f'(t)|dt\right).
\ee

Next we consider the case when $\tau<\kappa^2$. Now $\phi''(t)$ has only one zero in $(0,\kappa)$ at $t=t_2$ and either $\phi'''(t)\ge 0$ in $(0,\kappa)$ when $\kappa^3\le\tau<\kappa$ or $\phi'''(t)$ has only one zero in $(0,\kappa)$ at $t_3$ when $\tau<\kappa^3$, where
\ben
t_2^2=\kappa^2-\frac{1-\kappa^2}{(\tau\kappa^2)^{-2/3}-1},\ \ t_3^2=\kappa^2-\frac{1-\kappa^2}{(\tau\kappa^2)^{-2/5}-1}.
\een

When $\kappa^3\le\tau<\kappa$, $\phi''(t)$ is increasing in $(0,\kappa)$. Thus for sufficiently small $\de>0$,
\be\label{k3}
\hskip-1cm|\phi''(t)|\ge \min(|\phi''(t_2+\de)|,|\phi''(t_2-\de)|),\ \ \forall t\in (0,t_2-\de)\cup( t_2+\de,\kappa).
\ee
On the other hand, when $\tau<\kappa^3$, we have $t_3<t_2$ and $\phi'''(t)\ge 0$ for $t\ge t_3$ and $\phi'''(t)\le 0$ for $t\le t_3$. Therefore, $\phi''(t)$ is increasing in $(t_3,\kappa)$ and decreasing in $(0, t_3)$. Thus
\be\label{k4}
 \hskip-2cm|\phi''(t)|\ge \min(|\phi''(t_2+\de)|,|\phi''(t_2-\de)|,|\phi''(0)|),\ \ \forall t\in (0,t_2-\de)\cup(t_2+\de,\kappa).
\ee

To estimate the lower bound of $|\phi''(t_2\pm\de)|$, we observe that 
since $\tau\kappa^2<\kappa^4$, $t_2^2\ge\kappa^2-(1-\kappa^2)/(\kappa^{-8/3}-1)$, and thus $|\phi'''(t_2)|\ge c_0\tau\ge c_0\tau_0$ for some constant $c_0$ depending only on $\kappa$. 
Moreover, for any $t\in [t_2-\de,t_2+\de]$, $|\phi'''(t)-\phi'''(t_2)|\le C_1\de$ for some constant $C_1$ depending only on $\kappa$. Thus, if $\de\le c_0\tau_0/(2C_1)$, $|\phi'''(t)|\ge c_0\tau_0/2$ in $[t_2-\de,t_2+\de]$. This implies by using the mean value theorem that $|\phi''(t_2\pm\de)|\ge (c_0\tau_0/2)\de$. Notice that $|\phi''(0)|=1-\tau\kappa^{-1}\ge 1-\kappa$, from (\ref{k3})-(\ref{k4}) we conclude that for sufficiently small $\de>0$,
\be\label{k5}
|\phi''(t)|\ge (c_0\tau_0/2)\de,\ \ \forall t\in (0,t_2-\de)\cup( t_2+\de,\kappa).
\ee
Now we split the integral
\ben
\hskip-2cm\int^{\kappa}_0f(t)e^{\i\lam\phi(t)}dt&=&\int^{t_2-\de}_0f(t)e^{\i\lam\phi(t)}dt+\int^{t_2+\de}_{t_2-\de}f(t)e^{\i\lam\phi(t)}dt+\int_{t_2+\de}^\kappa f(t)e^{\i\lam\phi(t)}dt\\
&:=&{\rm II}_1+{\rm II}_2+{\rm II}_3.
\een
From (\ref{k5}), by Van der Corput Lemma \ref{van}, we have 
\ben
|{\rm II}_1+{\rm II}_3| \le C(\lam\de)^{-1/2}\left(|f(0)|+\int^\kappa_0|f'(t)|dt\right).
\een
It is obvious that $|{\rm II}_2|\le 2\de\max_{t\in (0,\kappa)}|f(t)|$. This yields after taking $\de=\lam^{-1/3}$, 
\ben
\left|\int^\kappa_{0}f(t)e^{\i\lam\phi(t)}dt\right|\leq C\lam^{-1/3}\left(|f(0)|+\int_{-\kappa}^{\kappa}|f'(t)|dt\right).
\een
This completes the proof by noticing (\ref{k2}).
\finproof

The following theorem is the main result of this section.
\begin{thm}\label{thm:3.1}
Let $k_s h\ge 1$. There exists a constant $C$ depending only on $\kappa$ such that for any $z,y\in\Om$,
\ben
|\J(z,y)-\F(z,y)|+k_s^{-1}|\na_y(\J(z,y)-\F(z,y))|\leq \frac{C}{\mu}(k_s h)^{-1/4}.
\een
\end{thm}

\debproof
By Lemma \ref{lem:3.2}, Lemma \ref{lem:3.4} and the definitions of $\J(z,y),\F(z,y)$ in (\ref{d3})-(\ref{d4}), we know that we are left to estimate
\ben
& &\frac 1{2\pi}\sum_{\al,\beta=p,s}\int_{-k_s}^{k_s}\frac{{\Ta}(\xi)^T \overline{\Nb(\xi)}}{\overline{\delta(\xi)}} e^{\i (\mu_\alpha z_2-\overline{\mu}_\beta y_2)+\i(y_1-z_1)\xi} d\xi-\F(z,y)\\
\hskip-1cm&=&\frac {1}{2\pi}\sum_{\stackrel{\al,\beta=p,s}{_{(\al,\beta)\not= (s,s)}}}\int_{(-k_s,k_s)\backslash[-k_p,k_p]}\frac{{\Ta}(\xi)^T \overline{\Nb(\xi)}}{\overline{\delta(\xi)}} e^{\i (\mu_\alpha z_2-\overline{\mu}_\beta y_2)+\i(y_1-z_1)\xi} d\xi\\
\hskip-1cm&+&\frac 1{2\pi}\int_{-k_p}^{k_p}\left[\frac{\Tp(\xi)\overline{\Ns(\xi)}}{\overline{\de(\xi)}}e^{\i(\mu_py_2-\mu_s z_2)}+\frac{\Ts(\xi)\overline{\Np(\xi)}}{\overline{\de(\xi)}}e^{\i(\mu_sy_2-\mu_p z_2)}\right]e^{\i(y_1-z_1)\xi}d\xi\\
\hskip-1cm&:=&{\rm II}_1+{\rm II}_2.
\een
For $k_p<|\xi|<k_s$, we know that $|\de(\xi)|\ge Ck_s^4$ by Lemma \ref{delta}, and for $\al,\beta=p,s$, $|\Ta(\xi)|\le C, |\Nb(\xi)|\le C\mu^{-1}k_s^2$. This implies
\ben
|{\rm II_1}|\le \frac{C}{k_s\mu}\int^{k_s}_{k_p}e^{-\sqrt{\xi^2-k_p^2}h}d\xi\le\frac C\mu (k_sh)^{-1}.
\een
For the term ${\rm II}_2$ we use Lemma \ref{lem:3.5}. The first term in ${\rm II}_2$ can be reduced to the integral in Lemma \ref{lem:3.5} by setting
\ben
f(t)=k_s\frac{\Tp(k_st)\overline{\Ns(k_st)}}{\overline{\de(k_st)}},\ \ \ \lam=k_sz_2,\tau=\frac{y_2}{z_2},\nu=\frac{y_1-z_1}{z_2}.
\een
By the assumption (\ref{d0}), it is then straightforward by using Lemma \ref{lem:3.5} to see that
\ben
\left|\int_{-k_p}^{k_p}\frac{\Tp(\xi)\overline{\Ns(\xi)}}{\overline{\de(\xi)}}e^{\i(\mu_py_2-\mu_s z_2)}d\xi\right|\le\frac{C}{\mu}(k_sh)^{-1/4}.
\een
The second integral in ${\rm II}_2$ can be estimated similarly. This completes the proof.
\finproof

The following theorem shows that $\F(z,y)$ has the similar behavior as the imaginary part of the elastic fundamental solution $\Im\G(z,y)$.

\begin{thm} \label{thm:3.2}
For any $z,y\in \R_+^2$, $\F(z,y)^T=\F(z,y)$. When $z=y$, $\Im [\F(z,y)]_{12} = \Im [\F(z,y)]_{21} =0$ and
	\be\label{d6}
-\Im [\F(z,y)]_{ii}\geq \frac{1}{4(\lambda+2\mu)} \ , \ i =1 ,2.
	\ee
When $z\neq y$,
	\be\label{d7}
	|\F(z,y)|&\le \frac{C}{\mu}\left(\frac 1{(k_s|z-y|)^{1/2}}+\frac 1{k_s|z-y|}\right),
	\ee
	where constant $C$ depends only on $\kappa$.
\end{thm}

\debproof
Substitute (\ref{d1}) and (\ref{d2}) into (\ref{d4}), we obtain
\be   
\hskip-1.5cm\F(z,y)&=&-\frac{1}{2\pi}\int_{-k_p}^{k_p} \frac{\i k_s^2\mu_s}{\mu\gamma(\xi)\delta(\xi)}
\Bigg(
\begin{array}{cc}
	\xi^2 & -\xi\mu_p \\
	-\xi\mu_p & \mu_p^2
\end{array}\Bigg)e^{\i\mu_p (z_2-y_2) +\i\xi(y_1-z_1)}d\xi \nn\\
\hskip-1.5cm& &-\frac{1}{2\pi}\int_{-k_p}^{k_p} \frac{\i k_s^2\mu_p}{\mu\gamma(\xi)\delta(\xi)}
\Bigg(
\begin{array}{cc}
	\mu_s^2 & \xi\mu_s \\
	\xi\mu_s & \xi^2
\end{array}		\Bigg)e^{\i\mu_s (z_2-y_2) +\i\xi(y_1-z_1)}d\xi \nn\\ 
\hskip-1.5cm& &
-\frac{1}{2\pi}\int_{(-k_s,k_s)\bks[-k_p,k_p]} \frac{\i(k_s^2-4\xi^2)\mu_p}{\mu\gamma(\xi)\overline{\delta(\xi)}}
\Bigg(
\begin{array}{cc}
	\mu_s^2 & \xi\mu_s \\
	\xi\mu_s & \xi^2
\end{array}		\Bigg)e^{\i\mu_s (z_2-y_2) +\i\xi(y_1-z_1)}d\xi \nn\\
\hskip-1.5cm&:=&{\rm III}_1+{\rm III}_2+{\rm III}_3. \label{d8}
\ee
It is easy to show that $\Im [\F(z,y)]_{12} = \Im [\F(z,y)]_{21} =0$ when $z=y$. 

Now we show the inequality (\ref{d6}) for the case of $i=j=1$. The other case is similar. 
Notice that for $\xi\in (-k_p,k_p)$, $\delta(\xi)\le k_s^4$ and $\mu_p\le\mu_s$. Then, if $z=y$,
\ben
\hskip-1.5cm-\Im ({\rm III}_1+{\rm III}_2)=\frac{1}{2\pi\mu}\int^{k_p}_{-k_p}\frac{k_s^2\mu_s}{\de(\xi)}d\xi\geq\frac{1}{2\pi\mu}\int_{-k_p}^{k_p} \frac{\mu_p}{k_s^2}d\xi = \frac{1}{4(\lambda+2\mu)}.
\een
If $\xi\in(-k_s,k_s)\bks[-k_p,k_p]$, $\mu_p=\i\sqrt{\xi^2-k_p^2}$, we have
\ben
\hspace{-1.5cm}
-{\rm III}_3=\frac{1}{2\pi\mu}\int_{(-k_s,k_s)\bks(-k_p,k_p)} \frac{\mu_s^2\sqrt{\xi^2-k_p^2}(k_s^2-4\xi^2)}{(\xi^2+\i\mu_s\sqrt{\xi^2-k_p^2})(\beta^2-\i4\xi^2\mu_s\sqrt{\xi^2-k_p^2})} d\xi.
\een
A simple computation shows that $\Im[(\xi^2+\i\mu_s\sqrt{\xi^2-k_p^2})(\beta^2-\i4\xi^2\mu_s\sqrt{\xi^2-k_p^2})]=k_s^2\mu_s\sqrt{\xi^2-k_p^2}(k_s^2-4\xi^2)$. It is then clear that $-\Im({\rm III}_3)\ge 0$. This shows $-\Im[\F(z,y)]_{11}\ge 1/[4(\lam+2\mu)]$ when $z=y$.

For $z\neq y$, we denote $y-z=|y-z|(\cos\phi,\sin\phi)^T$ for some $0\le\phi\le \pi$. Then it is easy to see that
\ben
{\rm III}_1=\frac{1}{\mu}\int_{0}^{\pi} A(\theta,\kappa) e^{\i k_s |z-y| \cos(\theta-\phi)}d\theta,
\een
for some function $A(\theta,\kappa)$. By Van der Corput Lemma \ref{van}, we can show easily
\ben
|{\rm III}_1|\le \frac C\mu\left(\frac 1{(k_s|z-y|)^{1/2}}+\frac 1{k_s|z-y|}\right).
\een
The estimate for ${\rm III}_2+{\rm III}_3$ can be proved similarly. 
This completes the proof.
\finproof

\section{The reverse time migration algorithm}

We start by introducing some notation. For any Lipschitz domain $\mathcal{D}\subset \R^2$ with boundary $\Ga_\mathcal{D}$, let $\|u\|_{H^1(\mathcal{D})}=(\|\na \phi\|_{L^2(\mathcal{D})}^2+d_\mathcal{D}^{-2}\|\phi\|_{L^2(\mathcal{D})}^2)^{1/2}$ be the weighted $H^1(\mathcal{D})$ norm
and
$\|v\|_{H^{1/2}(\Ga_\mathcal{D})}=(d_\mathcal{D}^{-1}\|v\|_{L^2(\Ga_\mathcal{D})}^2+|v|_{\frac 12,\Ga_\mathcal{D}}^2)^{1/2}$ be the weighted $H^{1/2}(\Ga_\mathcal{D})$ norm,
where $d_\mathcal{D}$ is the diameter of $\mathcal{D}$ and
\ben
|v|_{\frac 12,\Ga_\mathcal{D}}=\left(\int_{\Ga_\mathcal{D}}\int_{\Ga_\mathcal{D}}\frac{|v(x)-v(y)|^2}{|x-y|^2}ds(x)ds(y)\right)^{1/2}.
\een
By the scaling argument and trace theorem we know that there exists a constant $C>0$ independent of $d_\mathcal{D}$ such that for any $\phi\in C^1(\bar{\mathcal{D}})^2$ \cite[corollary 3.1]{RTMhalf_aco},
\be\label{q0}
\|\phi\|_{H^{1/2}(\Ga_\mathcal{D})}+\|\sigma(\phi)\nu\|_{H^{-1/2}(\Ga_\mathcal{D})}\le C\max_{x\in \bar\mathcal{D}}(|\phi(x)|+d_\mathcal{D}|\na\phi(x)|).
\ee
In this paper, for any Sobolev space $X$, we still denote $X$ the vector valued space $X^2$ or tensor valued space $X^{2\times 2}$. The norms of $X, X^2, X^{2\times 2}$ are all denoted by $\|\cdot\|_X$.

\begin{lem}\label{lem:4.1}
Let $k_s h\geq 1, d\gg h$, there exists a constant $C$ depending only on $\kappa$ but independent of $k_s, h, d, d_D$ such that for any $z\in\Om$, $j=1,2$,
\ben
\hskip-2cm& &\|\F(z,\cdot)e_j\|_{H^{1/2}(\Ga_D)}+\|\sigma(\F(z,\cdot)e_j)\nu\|_{H^{-1/2}(\Ga_D)}\le\frac C\mu(1+k_sd_D),\\\hspace{-2cm}
\hskip-2cm& &\|\R_d(z,\cdot)e_j\|_{H^{1/2}(\Gamma_D)}+\|\sigma(\R_d(z,\cdot)e_j)\nu\|_{H^{-1/2}(\Gamma_D)} \le
\frac{C}{\mu}(1+k_sd_D)\left[\left(\frac hd\right)^{2}+(k_sh)^{-1/4}\right],
\een	
where $\R_d(z,\cdot)=\J_d(z,\cdot)-\F(z,\cdot)$.
\end{lem}

\debproof
The first estimate follows easily from (\ref{q0}) and the definition of $\F(z,\cdot)$ in (\ref{d4}). The second estimate follows from (\ref{q0}), Lemma \ref{lem:3.1} and Theorem \ref{thm:3.1}. This completes the proof.
\finproof

Now we briefly recall the classical argument of limiting absorption principle (see e.g. \cite{leis, wilcox1975, Yves1988}) to define the scattering solution for the exterior elastic scattering problem in the half space:
\be
\Delta_e u + \omega^2 u =0 \ \ \mbox{\rm in } \R^2_+\bks \bar{D}, \label{pp1}\\
u= g \ \ \mbox{\rm on } \Ga_D, \ \ \ \ \sigma(u)e_2=0 \ \ \mbox{\rm on } \Ga_0,  \label{pp2}
\ee
where $g \in H^{1/2}(\Ga_D)$. Let $\eps>0$ and $u_\eps$ be the solution of the problem
\be
\Delta_e u_\eps + [\omega(1+\i\eps)]^2 u_\eps =0 \ \ \mbox{\rm in } \R^2_+\bks \bar{D}, \label{pp3}\\
u_\eps= g \ \ \mbox{\rm on } \Ga_D, \ \ \ \ \sigma(u_\eps)e_2=0 \ \ \mbox{\rm on } \Ga_0.  \label{pp4}
\ee
By the Lax-Milgram lemma, the problem (\ref{pp3})-(\ref{pp4}) has a unique solution $u_\eps\in H^1(\R^2_+\bks\bar D)$. Let $\mathcal{D}(\De_e)=\{v\in H^1(\R^2_+\bks\bar D): \De_e v\in L^2(\R^2_+\bks\bar D), v=0\ \ \mbox{on }\Ga_D, \sigma(v)e_2=0\ \ \mbox{on }\Ga_0\}$ as the domain of the operator $-\De_e$, it is shown in \cite{Yves1988} that if $\om^2$ is not the eigenvalue for $-\De_e$ in the domain $\mathcal D(\De_e)$, $u_\eps$ converges to some function $u$ satisfying (\ref{pp1})-(\ref{pp2}) in $H^{1,-s}(\R^2_+\bks\bar D)$, $s>1/2$, where the weighted Sobolev space $H^{1,s}(\R^2_+\bks\bar D),s \in \R$, is defined as the set of functions in $L^{2,s}(\R^2_+\bks\bar D)=\{v \in L^2_{\rm loc}(\R^2_+\bks\bar D): (1+|x|^2)^{s/2}v \in L^2(\R^2_+\bks\bar D) \}$ whose first derivatives are also in $L^{2,s}(\R^2_+\bks\bar D)$. The norm $\| v \|_{ H^{1,s}(\R^2_+\bks\bar D)} = (\| v \|^2_{ L^{2,s} (\R^2_+\bks\bar D)} + \| \nabla v \|^2_{ L^{2,s}(\R^2_+\bks\bar D)})^{1/2}$, where $\| v \|_{ L^{2,s}(\mathcal D)} = (\int_{\mathcal D}(1+|x|^2)^{s}|v|^2 dx )^{1/2}$. The absence of the positive eigenvalue for the operator $-\De_e$ is proved in \cite{sini2004} in the domain $\mathcal D'(\De_e)=\{v\in H^1(\R^2_+\bks\bar D), \De_e v\in L^2(\R^2_+\bks\bar D), \sigma(v)\nu=0\ \mbox{on }\Ga_D, \sigma(v)e_2=0\ \mbox{on }\Ga_0\}$. One can easily extend the argument in \cite{sini2004} to show the absence of the positive eigenvalue for
$-\De_e$ also in the domain $\mathcal D(\De_e)$ and thus obtain the following theorem for the forward scattering problem.

\begin{thm} \label{thm:4.1}
Let $g \in H^{1/2}(\Ga_D)$. The half-space elastic scattering problem (\ref{pp1})-(\ref{pp2})
admits a unique solution $u\in H^{1}_{\rm loc}(\R^2_+ \backslash \bar D)$. Moreover, for any bounded open set $\mathcal O\subset \R^2_+\bks\bar D$ there exists a constant $C>0$ such that
$\|u\|_{H^{1}(\mathcal O)}\le C\|g\|_{H^{1/2}(\Ga_D)}$.
\end{thm}

For the sake of convenience, we introduce the following notation: for any $u,v\in H^1(\R^2\bks\bar D)$ such that $\De_e u, \De_e v\in L^2(\R^2\bks\bar D)$,
\be\label{g1}
\GG(u,v)=\int_{\Ga_D} [u(x)\cdot \sigma(v(x))\nu- \sigma(u(x))\nu\cdot v(x)]ds(x).
\ee
Using this notation, the integral representation formula for the solution of the half-space elastic scattering problem reads: 
\be\label{g2}
u(y)\cdot q=\GG(u(\cdot),\N(\cdot,y)q), \ \ \forall y\in\R^2_+\bks\bar D,\ \ \forall q\in\R^2.
\ee

Now we introduce the RTM algorithm for the half-space inverse elastic scattering problem. Assume that there are $N_s\ge 1$ sources and $N_r\ge 1$ receivers uniformly distributed on $\Gamma^d_0$. 

For any $q\in\R^2$, let $u^i_q$ be the incident field which satisfies
\ben
\Delta_e u_q^i(x,x_s) + \omega^2 u_q^i(x,x_s) =0 \ \ \mbox{\rm in } \R^2_+,\ \ u^i_q(x,x_s)=q\de_{x_s}(x)\ \ \mbox{on }\Ga_0.
\een
By the integral representation formula, $u^i_q(x,x_s)=\T_D(x_s,x)^Tq$. The following algorithm extends the algorithm in \cite{RTMhalf_aco, Zhang2007} for acoustic waves.

\begin{alg}{\sc (RTM algorithm for half-space elastic scattering data)}\label{alg_rtm}\\
Given the data $u_q^s(x_r,x_s)$ which is the measurement of the scattered field at $x_r$ when the source is emitted at $x_s$ along the polarized direction $q=e_1, e_2$, $s=1,\dots, N_s$, $r=1,\dots,N_r$. 

$1^\circ$ Back-propagation: Compute $v_q(x,x_s)$ as the scattering solution of the following half-space elastic scattering problem:
\ben
\Delta_e v_q(x,x_s) + \omega^2 v_q(x,x_s) =0 \ \ \mbox{\rm in } \R^2_+, \\
v_q(x,x_s)=\frac{|\Ga_0^d|}{N_r}\sum_{r=1}^{N_r}\overline{u_q^s(x_r,x_s)}\delta_{x_r}(x) \ \ \mbox{\rm on } \Ga_0.
\een
$2^\circ$ Cross-correlation: For each $z\in\Om$, compute the imaging function
\be\label{cor1} 
I_d(z)=\Im\sum_{q=e_1,e_2}\left\{\frac{|\Gamma_0^d|}{N_s}\sum^{N_s}_{s=1} u^i_q(z,x_s)\cdot v_q(z,x_s)\right\}. 
\ee
\end{alg}
By the integral representation formula, we know that
\ben
v_q(x,x_s)\cdot e_j=\frac{|\Ga_0^d|}{N_r}\sum_{r=1}^{N_r}\T_D(x_r,x)e_j\cdot\overline{u_q^s(x_r,x_s)},
\een
which yields
\be\label{cor}
\hskip-1cmI_d(z)=\Im\sum_{q=e_1,e_2}\left\{\frac{|\Gamma_0^d|^2}{N_sN_r}\sum^{N_s}_{s=1}\sum^{N_r}_{r=1}
[\T_D(x_s,z)^Tq]\cdot[\T_D(x_r,z)^T\overline{u^s_q(x_r,x_s)}]\right\}.
\ee
This is the formula will be used in our numerical examples in section 6.
By letting $N_s,N_r\to\infty$, we know that (\ref{cor1}) can be viewed as an approximation of the following continuous integral:
\be\hspace{-1.5cm}
\hat{I}_d(z)=\Im\sum_{q=e_1,e_2}\int_{\Gamma_0^d}\int_{\Gamma_0^d}\,
[\T_D(x_s,z)^Tq]\cdot[\T_D(x_r,z)^T\overline{u^s_q(x_r,x_s)}]\,ds(x_r)ds(x_s).\label{cor2}
\ee

The following theorem which extends \cite[Theorem 4.1]{RTMhalf_aco} for acoustic waves will be proved in the Appendix of this paper. It shows that the difference between the half-space scattering solution and the full space scattering solution is small when the scatterer is far away from the boundary $\Ga_0$.

\begin{thm}\label{thm:4.2}
Let $g\in H^{1/2}(\Ga_D)$ and $u_1,u_2$ be the scattering solution of following problems:
	\be\label{e1}
	\Delta_e u_1 + \omega^2 u_1=0 \ \ \mbox{\rm in } \R^2_+\bks \bar{D},\ \  u_1= g \ \ \mbox{\rm on } \Ga_D,\ \ \sigma(u_1)e_2=0 \ \ \mbox{\rm on } \Ga_0,\\
	\label{e2}
	\Delta_e u_2 + \omega^2 u_2=0 \ \ \mbox{\rm in }\R^2\bks \bar{D},\ \ u_2 = g \ \ \mbox{\rm on } \Ga_D.
	\ee
	Then there exits a constant C depending only on $\kappa$ but independent of $k_s, h,d_D$ such that
	\ben\hskip-2cm
	\|\sigma(u_1-u_2)\nu\|_{H^{-1/2}(\Gamma_D)}
	\le\frac{C}{\mu}(1+\|T_1\|)(1+\|T_2\|)(1+k_s d_D)^2(k_sh)^{-1/2}\|g\|_{ H^{1/2}(\Ga_D)}.
	\een
Here $T_1, T_2:H^{1/2}(\Ga_D)\to H^{-1/2}(\Ga_D)$ are the Dirichlet to Neumann mapping associated with the elastic scattering problem (\ref{e1}) and (\ref{e2}), respectively. $\|T_1\|, \|T_2\|$ denote their operator norms.
\end{thm}

We remark that the well-posedness of the full space elastic scattering problem (\ref{e2}) under the so-called Sommerfeld-Kupradze radiation condition is well known (cf. e.g. \cite{ku63}). It is equivalent to the solution defined by the limiting absorption principle \cite{leis, cxz2016}.

The following theorem, which relates the imaging function $\hat I_d(z)$ to the point spread function in section 3, shows the resolution of our RTM imaging algorithm for the half-space inverse elastic scattering problems.

\begin{thm}\label{thm:4.3}
	For any $z\in\Omega$, let $\U(z,x)\in\C^{2\times2}$ such that $\U(z,x)e_j$, $j=1,2$, is the scattering solution of the problem:
	\ben
	\hskip-2cm\Delta_e [\U(z,x)e_j]+ \omega^2[\U(z,x)e_j]= 0 \ \ \mbox{in }\R^2\bks \bar{D},\ \ \ \ 
	\U(z,x)e_j= -\overline{\F(z,x)}e_j \ \ \mbox{on }\Ga_D.  
	\een
	Then, we have
	\ben\hspace{-2cm}
	\hat{I}_d(z)=\Im\sum_{j=1}^2\int_{\Gamma_D}[\sigma(\U(z,x)e_j+\overline{\F(z,x)}e_j)\nu]\cdot [\overline{\F(z,x)}e_j]ds(x)+R_d(z),
	\een
	where $|R_d(z)|\leq C\mu^{-2}(1+\|T_1\|)(1+\|T_2\|)(1+k_s d_D)^3\left[\left(\frac hd\right)^{2}+(k_sh)^{-1/4}\right]$ for some constant $C$ depending only on $\kappa$ but independent of $k_s,k_p, h, d, d_D$.
\end{thm}
\debproof
From (\ref{cor2}) we know that
\be\label{g5}
\hat I_d(z)=\Im\sum_{q=e_1,e_2}\int_{\Ga_0^d}[\T_D(x_s,z)^Tq]\cdot\hat v_q(z,x_s)ds(x_s),
\ee
where for $j=1,2$,
\ben
\hat v_q(z,x_s)\cdot e_j=\int_{\Ga_0^d}\T_D(x_r,z)e_j\cdot\overline{u^s_q(x_r,x_s)}ds(x_r).
\een
By (\ref{g2}) we know that $u^s_q(x_r,x_s)\cdot e_i=\GG(u^s_q(\cdot,x_s),\N(\cdot,x_r)e_i), i=1,2$, and thus
\ben
\hat v_q(z,x_s)\cdot e_j=\GG(\overline{u^s_q(\cdot,x_s)},\left[\int_{\Ga_0^d}\sum^2_{i=1}[\T_D(x_r,z)]_{ij}\overline{\N(\cdot,x_r)}e_ids(x_r)\right]\,).
\een
By using the reciprocity relation $\N(x,x_r)=\N(x_r,x)^T$ and the definition of $\J_d(\cdot,\cdot)$ in (\ref{jd}), we obtain
\be\label{g6}
\int_{\Ga_0^d}\sum^2_{i=1}[\T_D(x_r,z)]_{ij}\overline{\N(\cdot,x_r)}e_ids(x_r)=\J_d(z,x)^Te_j.
\ee
This implies $\hat v_q(z,x_s)e_j=\GG(\overline{u^s_q(\cdot,x_s)},\J_d(z,\cdot)^Te_j)$. Substitute it into (\ref{g5}) we have
\be\label{g3}
\hat I_d(z)=\Im\sum_{j=1}^2\GG(\W(z,\cdot)e_j,\J_d(z,\cdot)^Te_j),
\ee
where $\W(z,x)\in \C^{2\times 2}$ is the tensor defined by
\ben
\W(z,x)e_j=\int_{\Ga_0^d}\sum^2_{k=1}[\T_D(x_s,z)]_{kj}\overline{u^s_{e_k}(x,x_s)}ds(x_s),\ \ j=1,2.
\een
Notice that $\overline{\W(z,x)}e_j$ can be viewed as the weighted superposition of $u^s_{e_k}(x,x_s)$ and thus it satisfies
\be\label{g7}
\hskip-2cm\De_e[\overline{\W(z,x)}e_j]+\om^2[\overline{\W(z,x)}e_j]=0\ \ \mbox{in }\R^2_+\bks\bar D,\ \ \ \ \sigma(\overline{\W(z,x)}e_j)e_2=0\ \ \mbox{on }\Ga_0.
\ee
On the boundary of the obstacle $\Gamma_D$, by using (\ref{g6}) we have
\be
\hskip-1cm\overline{\W(z,x)}e_j=-\int_{\Ga_0^d}\sum^2_{k=1}[\overline{\T_D(x_s,z)}]_{kj}\N(x,x_s)e_kds(x_s)
=-\overline{\J_d(z,x)}^Te_j.\label{g8}
\ee
Now define the tensor $\W_d(z,x)\in \C^{2\times 2}$ such that $\W_d(z,x)e_j$, $j=1,2$, is the scattering solution of the problem
\be
\hskip-2cm& & \Delta_e [\W_d(z,x)e_j]+ \omega^2 [\W_d(z,x)e_j]= 0 \ \ \ \ \mbox{in }\R^2_+\bks \bar{D},\label{g9}\\
\hskip-2cm& &\W_d(z,x)e_j= -\overline{\F(z,x)}e_j \ \ \mbox{on }\Ga_D,\ \ \ \ 
\sigma(\W_d(z,x)e_j)e_2=0 \ \ \mbox{on }\Ga_0.\label{g10}
\ee
By (\ref{g3}) we deduce
\be
\hat I_d(z)&=&\Im\sum^2_{j=1}\GG(\W(z,\cdot)e_j,J_d(z,\cdot)^Te_j-\F(z,\cdot)e_j)\nn\\
& &+\Im\sum^2_{j=1}\GG(\W(z,\cdot)e_j-\overline{\W_d(z,\cdot)}e_j,\F(z,\cdot)e_j)\nn\\
& &+\Im\sum^2_{j=1}\GG(\overline{\W_d(z,\cdot)}e_j-\overline{\U(z,\cdot)}e_j,\F(z,\cdot)e_j)\nn\\
& &+\Im\sum^2_{j=1}\GG(\overline{\U(z,\cdot)}e_j,\F(z,\cdot)e_j):={\rm VI}_1+\cdots+{\rm VI}_4.\label{g11}
\ee
Recall that $\F(z,y)^T=\F(z,y)$ by Theorem \ref{thm:3.2}. By Lemma \ref{lem:4.1}, 
\ben
\|\J_d(z,\cdot)e_j\|_{H^{1/2}(\Ga_D)}+\|\sigma(\J_d(z,\cdot)e_j)\nu\|_{H^{-1/2}(\Ga_D)}\le\frac C\mu (1+k_sd_D).
\een
This implies, by (\ref{g7})-(\ref{g8}) and Lemma \ref{lem:4.1}, that
\ben
|{\rm VI}_1|&\le&\sum_{j=1}^2\Big(\|\W(z,\cdot)e_j\|_{H^{1/2}(\Ga_D)}\|\sigma(\J_d(z,\cdot)^Te_j-\F(z,\cdot)e_j)e_2\|_{H^{-1/2}(\Ga_D)}\\
& &+\|\sigma(\W(z,\cdot)e_j)e_2\|_{H^{-1/2}(\Ga_D)}\|\J_d(z,\cdot)^Te_j-\F(z,\cdot)e_j\|_{H^{1/2}(\Ga_D)}\Big)\\
&\le&\frac C{\mu^2}(1+\|T_1\|)(1+k_sd_D)^2\left[\left(\frac hd\right)^{2}+(k_sh)^{-1/4}\right].
\een
From (\ref{g7})-(\ref{g8}) and (\ref{g9})-(\ref{g10}), we obtain by using Lemma \ref{lem:4.1} that
\ben
|{\rm VI}_2|\le\frac C{\mu^2}(1+\|T_1\|)(1+k_sd_D)^2\left[\left(\frac hd\right)^{2}+(k_sh)^{-1/4}\right].
\een
To estimate the third term, we use Theoremn \ref{thm:4.2} and Lemma {\ref{lem:4.1} to obtain
\ben
|{\rm VI}_3|\le\frac C{\mu^2}(1+\|T_1\|)(1+\|T_2\|)(1+k_sd_D)^3(k_sh)^{-1/2}.
\een
Finally, since $\U(z,x)e_j=-\overline{\F(z,x)}e_j$ on $\Ga_D$,
\ben
\hskip-1.5cm{\rm IV}_4&=&\Im\sum^2_{j=1}\int_{\Ga_D}(\overline{\U(z,x)}e_j\cdot\sigma(\F(z,x)e_j)\nu-\sigma(\overline{\U(z,x)}e_j)\nu\cdot\F(z,x)e_j)ds(x)\\
\hskip-1.5cm&=&-\Im\sum^2_{j=1}\int_{\Ga_D}\sigma(\overline{\U(z,x)}e_j+\F(z,x)e_j)\nu\cdot\F(z,x)e_jds(x)\\
\hskip-1.5cm&=&\Im\sum^2_{j=1}\int_{\Ga_D}\sigma(\U(z,x)e_j+\overline{\F(z,x)}e_j)\nu\cdot\overline{\F(z,x)}e_jds(x).
\een
The theorem follows now from (\ref{g11}).
\finproof

By (\ref{d8}) we know that for any fixed $z\in\Om$ and some functions $A_j(\xi), B_j(\xi)$, $j=1,2$,
\ben
\hskip-2cm\F(z,x)e_j&=&\int_{-k_p}^{k_p}A_j(\xi)\left(\begin{array}{c}
\hskip-6pt-\xi \hskip-6pt \\
\hskip-6pt \mu_p \hskip-6pt
\end{array}\right)e^{\i(z-x)\cdot(-\xi,\mu_p)^T}d\xi+\int_{-k_s}^{k_s}B_j(\xi)\left(\begin{array}{c}
\hskip-6pt\mu_s \hskip-6pt\\
\hskip-6pt\xi \hskip-6pt
\end{array}\right)e^{\i(z-x)\cdot(-\xi,\mu_s)^T}d\xi\\
\hskip-2cm&=&\int^\pi_0\left[\tilde A_j(\theta)\tau(\theta)e^{\i k_p(z-x)\cdot\tau(\theta)}+\tilde B_j(\theta)\tau(\theta)^\perp e^{\i k_s(z-x)\cdot\tau(\theta)}\right]d\theta,
\een
where $\tilde A_j(\theta)=k_pA_j(k_p\cos\theta)\sin\theta, \tilde B_j(\theta)=k_sB_j(k_s\cos\theta)\sin\theta$, $\tau(\theta)=(-\cos\theta,\sin\theta)^T$ and $\tau(\theta)^\perp=(\sin\theta,\cos\theta)^T$.
Thus $\overline{\F(z,x)}e_j$ is the weighted superposition of planar $p$ and $s$ waves and thus satisfies the elastic wave equation. Therefore, $\U(z,x)e_j$ can be viewed as the scattering solution of the elastic equation with the
incident wave $\overline{\F(z,x)}e_j$. By Theorem \ref{thm:3.2} we know that $\overline{\F(z,x)}$ decays as $|x-z|$ becomes large. Thus the imaging function $\hat{I}_d(z)$ becomes small when $z$ moves away from the boundary $\Ga_D$ if $k_s h \gg 1$ and $d\gg h$.

To understand the behavior of the imaging function when $z$ is close to the boundary of the scatterer, we introduce the concept of the scattering coefficient for incident plane waves.

\begin{definition}\label{scarr_con}
For any unit vector $\tau\in \R^2$, let $u^i_p =\tau e^{\i k_p x\cdot \tau}$, $u^i_s= \tau^\perp e^{\i k_s x\cdot \tau}$ be the incident planar $p$ and $s$ wave.  Let $u^s_\alpha (x) = u^s_\alpha(x;\tau), \al=p,s$, be the corresponding scattering solution of the elastic wave equation:
	\ben
	\De_e u^s_\alpha + \om^2u^s_\alpha = 0\ \ \mbox{in } \R^2\bks\bar{D}, \ \ \ \ 
	u^s_\alpha =-u^i_\alpha \ \ \mbox{on } \Ga_D.
	\een
	The scattering coefficient $R_\al(x;\tau)$, $x\in\Ga_D$, is defined by the relation
	\ben
	\sigma(u^s_\alpha(x)+u^i_\alpha(x))\nu= \i k_\alpha R_\alpha(x;\tau)e^{\i k_\alpha x\cdot \tau}  \ \ \ \mbox{on }\Ga_D.
	\een
Here for $\tau=(\tau_1,\tau_2)^T\in\R^2$, $\tau^\perp=(\tau_2,-\tau_1)^T$.
\end{definition}

With this definition we deduce from Theorem \ref{thm:4.3} that for any $z\in\Ga_D$,
\ben
\hskip-1.5cm\hat I_d(z)&\approx&\Im\sum^2_{j=1}\int_{\Ga_D}\left[
\int^\pi_0\overline{\tilde A_j(\theta)}\i k_pR_p(x;\tau(\theta))e^{\i k_p(x-z)\cdot\tau(\theta)}d\theta\right]\cdot\overline{\F(z,x)}e_j ds(x)\\
\hskip-1.5cm& &+\Im\sum^2_{j=1}\int_{\Ga_D}\left[
\int^\pi_0\overline{\tilde B_j(\theta)}\i k_sR_s(x;\tau(\theta))e^{\i k_s(x-z)\cdot\tau(\theta)}d\theta\right]\cdot\overline{\F(z,x)}e_j ds(x).
\een

In the case of Kirchhoff high-frequency approximation, the scattering coefficient is approximately zero in the shadow region of the obstacle: 
\ben
R_\alpha(x;\tau)\approx 0\ \ \mbox{if } x \in \Ga_D^{+}(\tau)=\{x\in \Ga_D, \nu(x)\cdot \tau>0\},\ \ \al=p,s.
\een
Let $x(s)$, $0<s<L$, be the arc length parametrization of the boundary $\Ga_D$ and $x_{\pm}(\theta)$ be the points on $\Ga_D$ such that $\nu(x_\pm(\theta))=\pm\tau(\theta)$. By using the method of the stationary phase and the above Kirchhoff approximation we can obtain as in \cite{RTMhalf_aco} that
\ben
\hskip-2.5cm\hat I_d(z)&\approx&\Im\sum^2_{j=1}\sqrt{2\pi k_p}
\int^\pi_0\overline{\tilde A_j(\theta)}e^{\i k_p(x_-(\theta)-z)\cdot\tau(\theta)+\i\frac\pi 4}\,\frac{R_p(x_-(\theta);\tau(\theta))\cdot\overline{\F(z,x_-(\theta))}e_j}{\sqrt{\kappa(x_-(\theta))}}d\theta\\
\hskip-2.5cm& &+\Im\sum^2_{j=1}\sqrt{2\pi k_s}
\int^\pi_0\overline{\tilde B_j(\theta)}e^{\i k_s(x_-(\theta)-z)\cdot\tau(\theta)+\i\frac\pi 4\,}\frac{R_s(x_-(\theta);\tau(\theta))\cdot\overline{\F(z,x_-(\theta))}e_j}{\sqrt{\kappa(x_-(\theta))}}d\theta.
\een
Here $\kappa(x)$ is the curvature of $\Ga_D$. 
Now for $z$ in the part of $\Ga_D$ which is back to $\Ga_0$, i.e., $\nu(z)\cdot\tau(\theta)>0$ for any $\theta\in (0,\pi)$, we know that $z$ and $x_-(\theta)$ are far away and thus $\hat{I}_d(z)\approx0$. This indicates that one cannot image the back part of the obstacle with only the data collected on $\Ga_0$. This is confirmed in our numerical examples in section 6.

\section{Extensions}

In this section we consider the reconstruction of non-penetrable obstacles with the impedance boundary condition and penetrable obstacles in the half space by the RTM algorithm \ref{alg_rtm}. For non-penetrable obstacles with the impedance boundary condition on the obstacle, the measured data $u_q(x_r,x_s)=u_q^s(x_r,x_s)+\N(x_r,x_s)q$, $q=e_1, e_2$, where $u^s_q(x,x_s)$ is the scattering solution of the following problem:
\ben
\hskip-1cm\Delta_e u^s_q(x,x_s) + \omega^2 u^s_q(x,x_s) =0\ \ \mbox{\rm in } \R^2_+\bks \bar{D}, \\
\hskip-1cm \sigma(u^s_q(x,x_s))\nu+\i\eta(x)u^s_q(x,x_s)=-[\sigma(\N(x,x_s)q)\nu+\i\eta(x)\N(x,x_s)q]\ \ \mbox{\rm on } \Ga_D, \\ 
\hskip-1cm\sigma(u^s_q(x,x_s))e_2=0 \ \ \mbox{\rm on } \Ga_0,
\een
where $\eta\in L^\infty(\Ga_D)$ and $\eta\ge 0$ on $\Ga_D$. By modifying the argument in the proof of Theorem \ref{thm:4.3}, we can show the following theorem whose proof is omitted.
\begin{thm}\label{thm:5.1}
	For any $z\in\Omega$, let $\U(z,x)\in\C^{2\times2}$ such that $\U(z,x)e_j$, $j=1,2$, is the scattering solution of the problem:
	\ben
	\hskip-1cm& & \Delta_e[\U(z,x)e_j]+ \omega^2[\U(z,x) e_j]= 0 \ \ \mbox{\rm in }\R^2\bks \bar{D},\\
	\hskip-1cm& &\sigma(\U(z,x)e_j)\nu+\i\eta(x)[\U(z,x)e_j]= -[\sigma(\overline{\F(z,x)}e_j)\nu+\i\eta(x)\overline{\F(z,x)}e_j ]\ \ \mbox{\rm on} \ \Ga_D.
	\een
Then the imaging function (\ref{cor2}) for the half-space elastic scattering data $u^s_q(x_r,x_s)$ of the non-penetrable obstacle with the impedance boundary condition satisfies
	\ben\hspace{-2.cm}
	\hat{I}_d(z)&=&-\Im\sum_{j=1}^2\int_{\Gamma_D} [\U(z,x)e_j+\overline{\F(z,x)}e_j]\cdot[\sigma(\overline{\F(z,x)}e_j)\nu+\i\eta(x)\overline{\F(z,x)}e_j]ds(x)\\
\hspace{-2.cm}	& &+R_d(z),\ \ \forall z\in\Om,
	\een
where $|R_d(z)|\leq C\mu^{-2}(1+k_s d_D)^3\left[\left(\frac hd\right)^{2}+(k_sh)^{-1/4}\right]$ for some constant $C$ depending only on $\kappa$ but independent of $k_s,k_p, h, d, d_D$.
\end{thm}

For the penetrable obstacle, the measured data is $u_q(x_r.x_s)=u_q^s(x_r,x_s)+\N(x_r,x_s)q$, $q=e_1,e_2$, where $u^s_q(x,x_s)$ is the scattering solution of the following problem:
\ben
\Delta_e u^s_q(x,x_s) + \omega^2n(x) u^s_q(x,x_s) =-\om^2(n(x)-1)\N(x,x_s)q\ \ \mbox{\rm in } \R^2_+, \\
\sigma(u^s_q(x,x_s))e_2=0 \ \ \mbox{\rm on }\Ga_0, 
\een
where $n(x)\in L^{\infty}({\R^2_+})$ is a positive function which is equal to 1 outside $D$. By modifying the argument in Theorem \ref{thm:4.3}, the following theorem can be proved.

\begin{thm}\label{resolution2}
	For any $z\in\Omega$, let $\U(z,x)\in\C^{2\times2}$ such that $\U(z,x)e_j$, $j=1,2$, is the scattering solution of the problem:
	\ben
	& & \Delta_e [\U(z,x)e_j] + \omega^2n(x)[\U(z,x)e_j]= -\omega^2(n(x)-1)\overline{\F(z,x)}e_j \ \ \mbox{\rm in }\R^2.
	\een
	Then the imaging function (\ref{cor2}) for the half-space elastic scattering data $u^s_q(x_r,x_s)$ of the penetrable obstacle satisfies
	\ben\hspace{-2cm}
	\hat{I}_d(z)=\Im\sum_{j=1}^2 \int_{D}\omega^2(n(x)-1)[(\U(z,x)e_j+\overline{\F(z,x)}e_j)\cdot\overline{\F(z,x)}e_j]dx+R_d(z),
	\een
	where $|R_d(z)|\leq C\mu^{-2}(1+k_s d_D)^3\left[\left(\frac hd\right)^{2}+(k_sh)^{-1/4}\right]$ for some constant $C$ depending only on $\kappa$ but independent of $k_s,k_p, h, d, d_D$.
\end{thm}

\section{Numerical examples}

In this section we present several numerical examples to show the effectiveness of our
RTM algorithm. To synthesize the scattering data we compute the solution $u^s_q(x,x_s)$ of
the scattering problems by representing the ansatz solution as the single layer potential
with the Neumann Green tensor $\N(x,y)$ as the kernel and discretizing the integral equation by
standard Nystr\"{o}m methods \cite{colton-kress}. The boundary integral equations on $\Ga_D$ are solved on
a uniform mesh over the boundary with ten points per probe wavelength. The sources
and receivers are both equally placed on the surface $\Ga^d_0$. In all our numerical examples we choose $h = 10, d = 50$ and {Lam\'{e}} constant $\lambda=1/2$, $\mu=1/4$. The boundaries
of the obstacles used in our numerical experiments are parameterized as follows, 
\ben
\mbox{Circle:}\ \ \ \ &&x_1=\rho\cos(\theta),\ \ x_2=\rho\sin(\theta);\ \ \\
\mbox{Kite:}\ \ \ \ &&x_1=\cos(\theta) + 0.65\cos(2\theta) - 0.65,\ \ x_2=1.5 \sin (\theta);\ \ \\
\mbox{$p$-leaf:}\ \ \ \ &&r(\theta)=1+0.2\cos(p\theta); \\
\mbox{peanut:}\ \ \ \ &&x_1 = \cos \theta + 0:2 \cos 3\theta; x_2 = \sin \theta + 0:2 \sin 3\theta; \\
\mbox{square:}\ \ \ \ &&x_1 = \cos3 \theta + \cos \theta; x_2 = \sin3\theta + \sin \theta.
\een
where
$\theta\in[0,2\pi]$. The numerical imaging function is (\ref{cor}) in section 4.

In the following by Dirichlet, Neumann or impedance obstacle we mean the non-penetrable obstacle that
satisfies Dirichlet, Neumann or impedance boundary condition on the boundary of the obstacle.

\bigskip
\textbf{Example 1}.
We consider imaging of a Dirichlet, a Neumann, an impedance, and a penetrable obstacle. The imaging domain $\Om=(-2,2) \times (8, 12)$ with the sampling grid $201 \times 201$. We set $N_s = N_r = 401$. The angular frequency $\om=2\pi$.
 \begin{figure}
 	\centering
 	\includegraphics[width=0.24\textwidth]{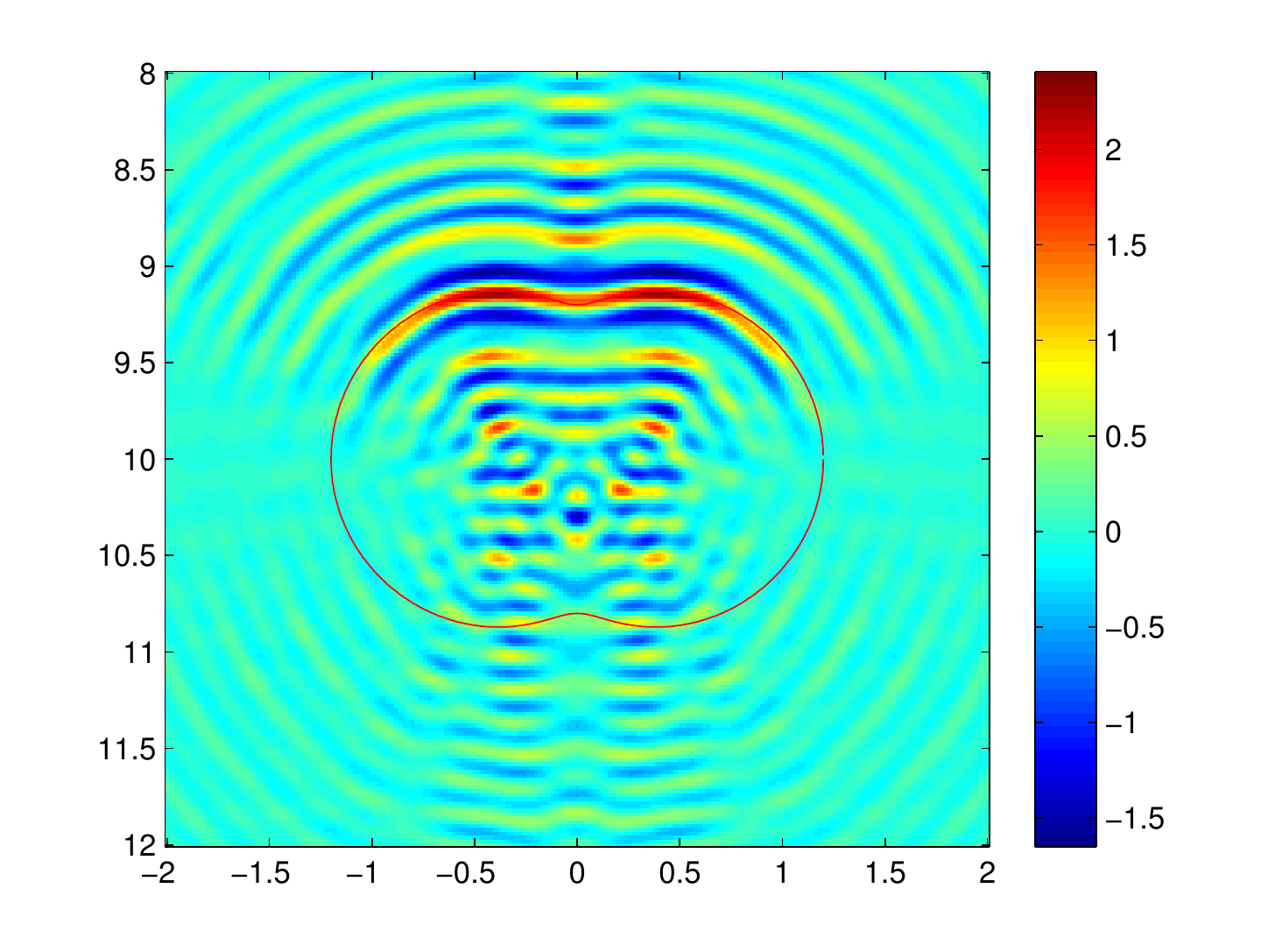}
 	\includegraphics[width=0.24\textwidth]{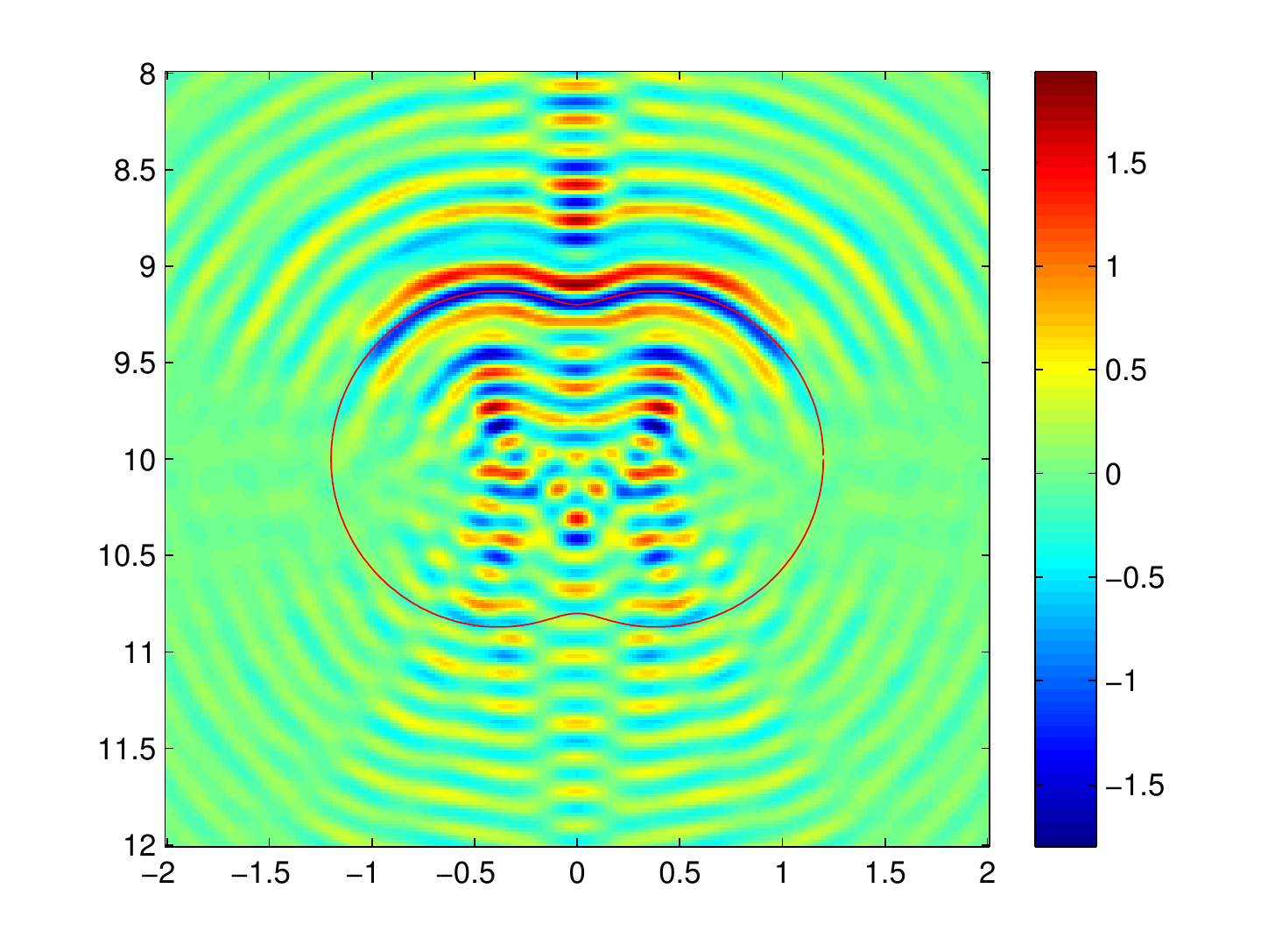}
 	\includegraphics[width=0.24\textwidth]{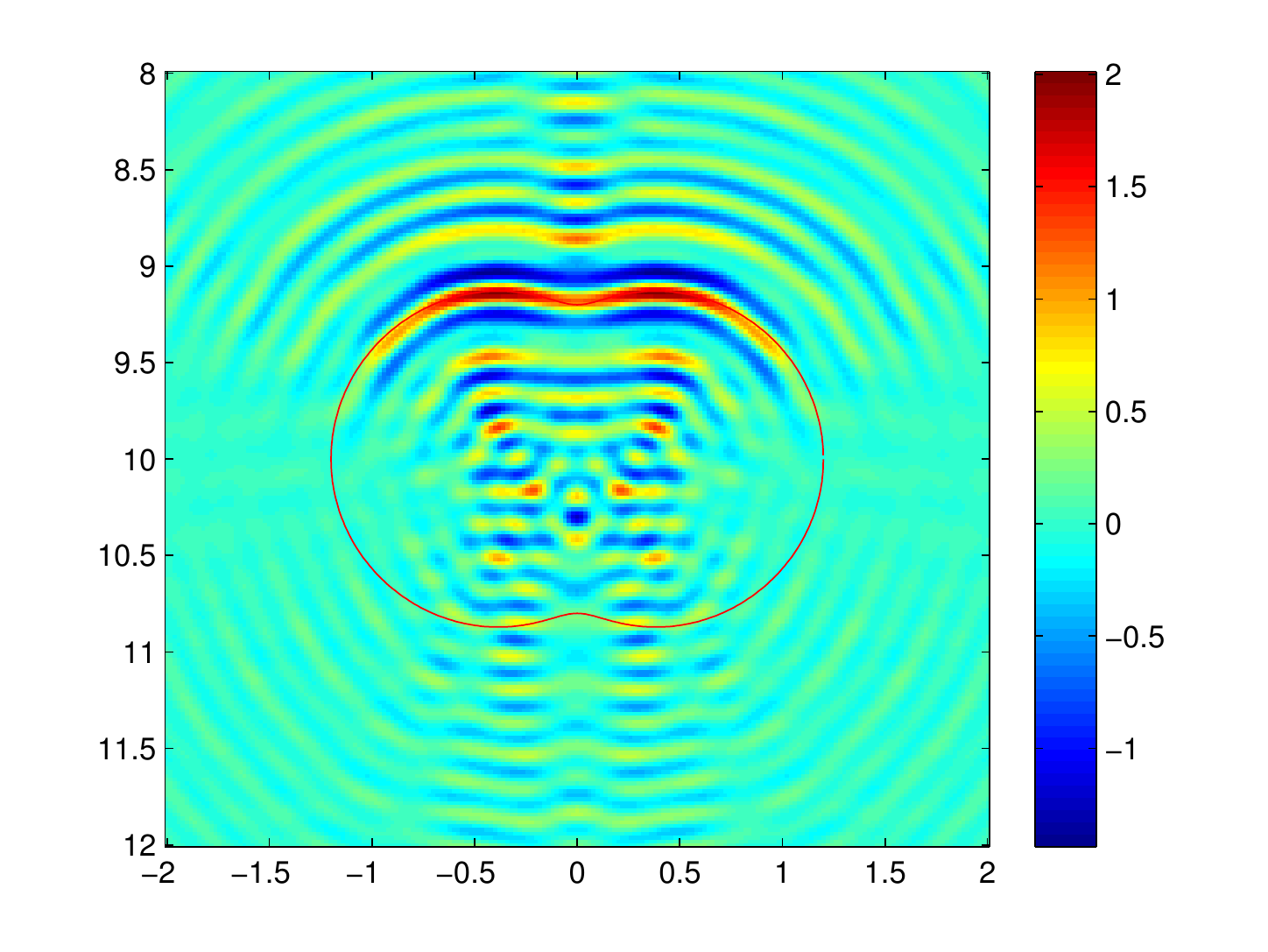}
 	\includegraphics[width=0.24\textwidth]{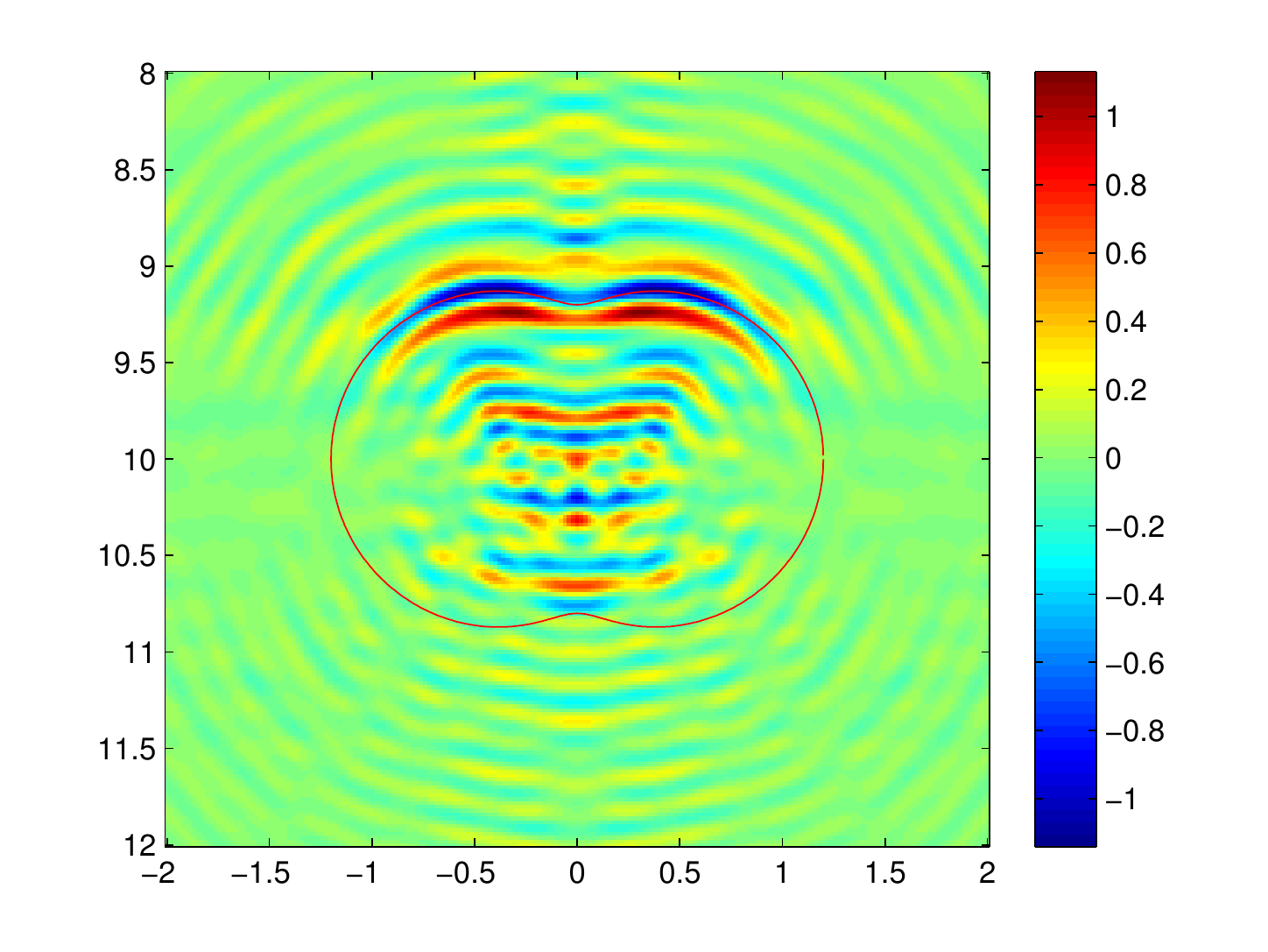}
 	\caption{Example 1: From left to right: imaging results of a Dirichlet, a Neumann, an impedance with $\eta(x)=1$, and a penetrable obstacle with diffractive index $n(x)=0.25$.} \label{figure_1}
 \end{figure}
 
 The imaging results are shown in Figure \ref{figure_1}. It demonstrates clearly that our RTM
 algorithm can effectively image the upper boundary illuminated by the sources and
 receivers distributed along the boundary $\Ga_0$ for non-penetrable obstacles. The imaging
 values decrease on the shadow part of the obstacles and at the points away from the
 boundary of the obstacle.

\bigskip
\textbf{Example 2}. We consider the imaging of Dirichlet obstacles with different shapes including a circle, a peanut, a $p$-leaf and a rounded square. The imaging domain $\Om = (-2, 2) \times (8, 12)$ with the sampling grid $201 \times 201$. We set $N_s = N_r = 401$. The angular frequency $\om = 3\pi,4\pi$ for the single frequency and $\om=\pi\times[2:0.5:8]$ for the test of multiple frequencies.

\begin{figure}
	\centering
	\includegraphics[width=0.32\textwidth]{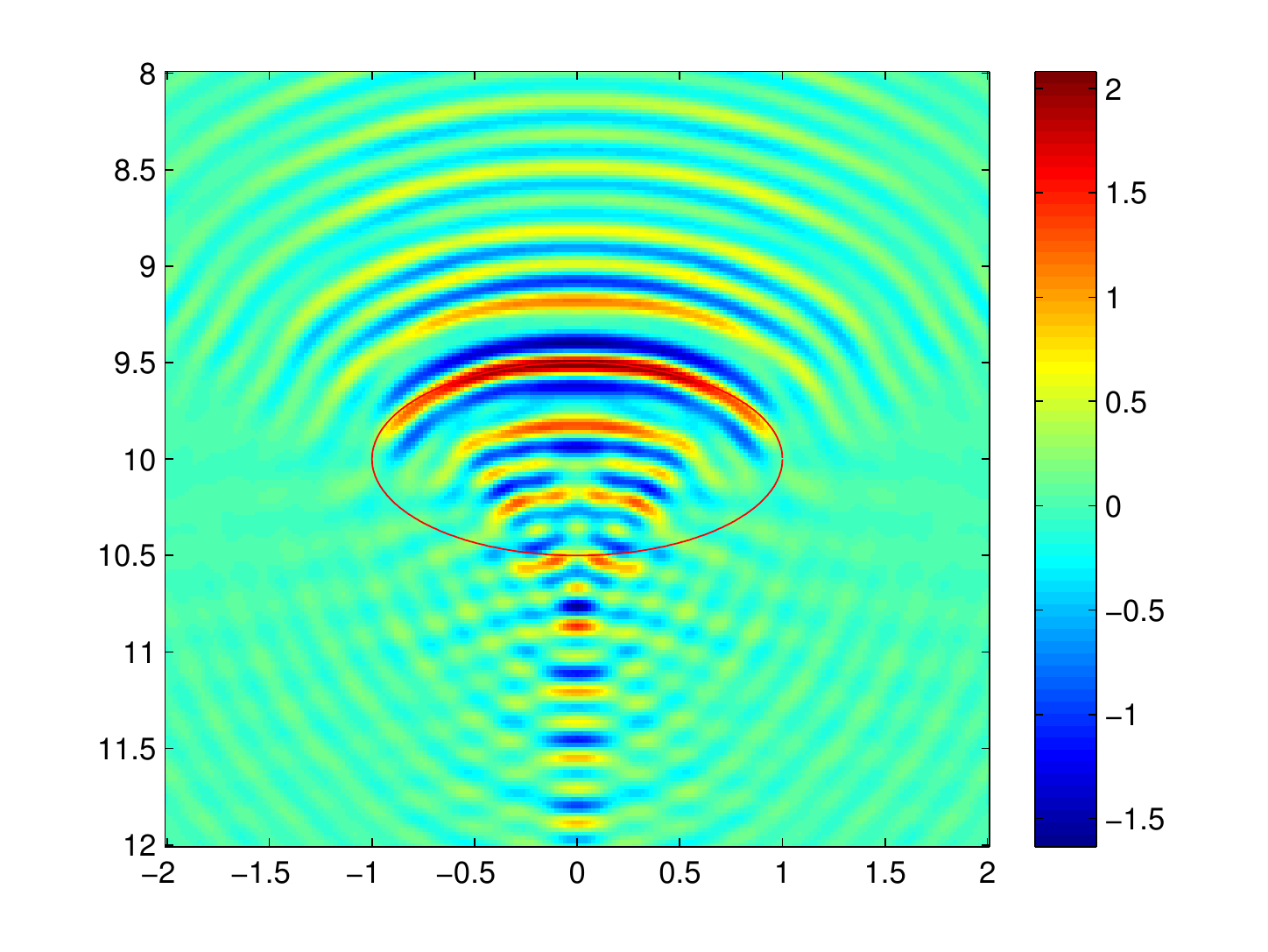}
	\includegraphics[width=0.32\textwidth]{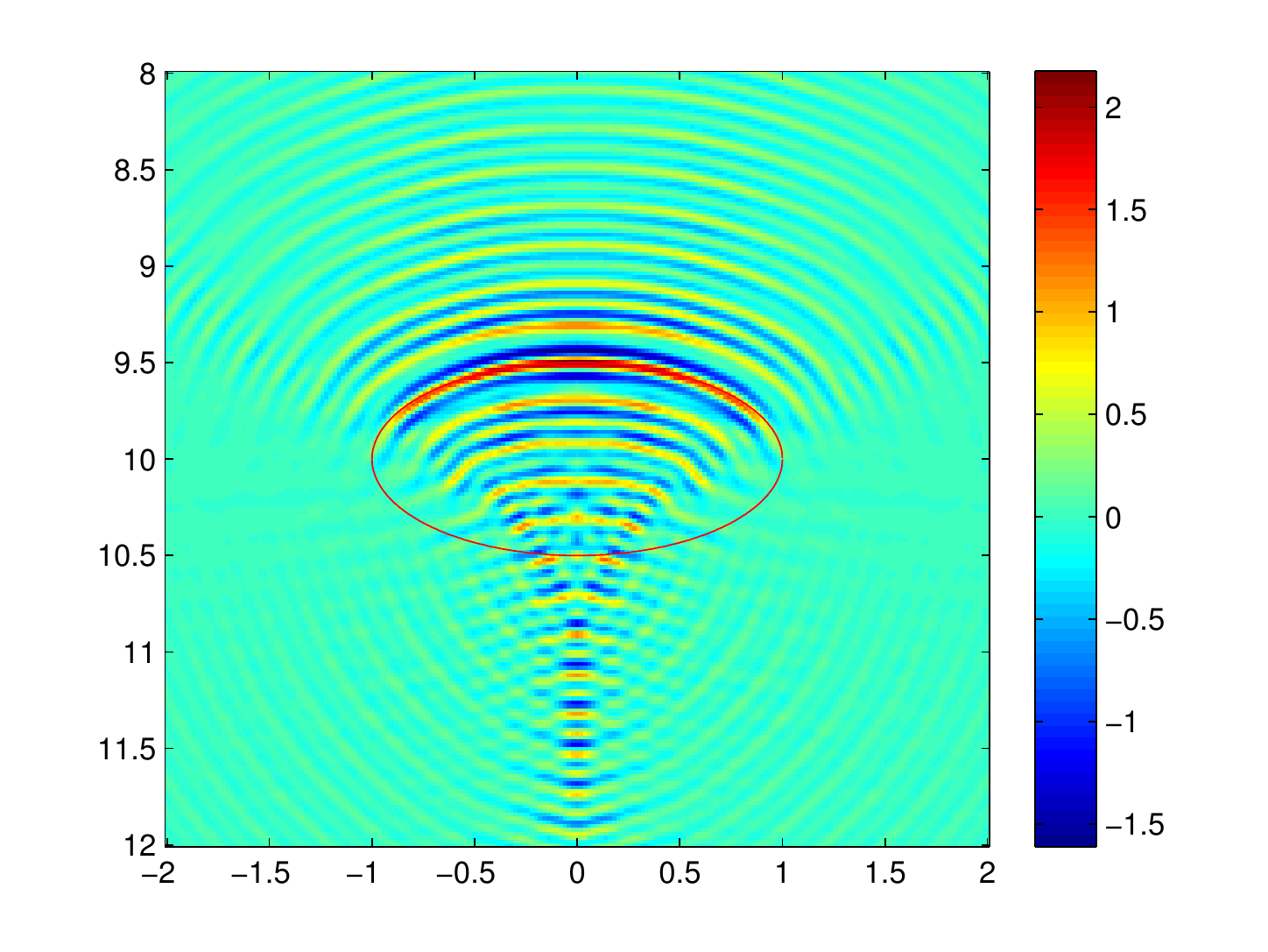}
	\includegraphics[width=0.32\textwidth]{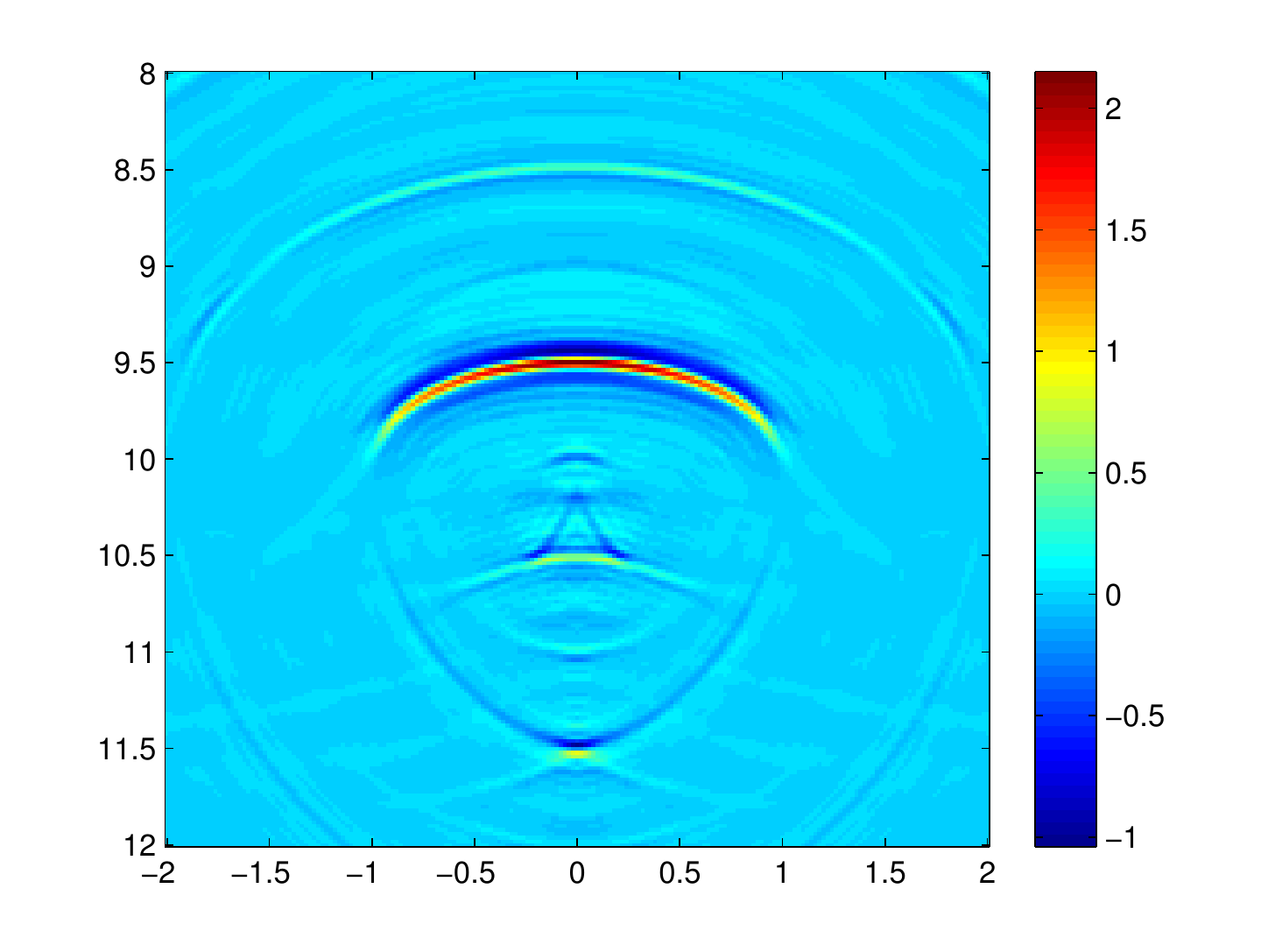}
	\includegraphics[width=0.32\textwidth]{./graphic/peanut_3pi-eps-converted-to.pdf}
	\includegraphics[width=0.32\textwidth]{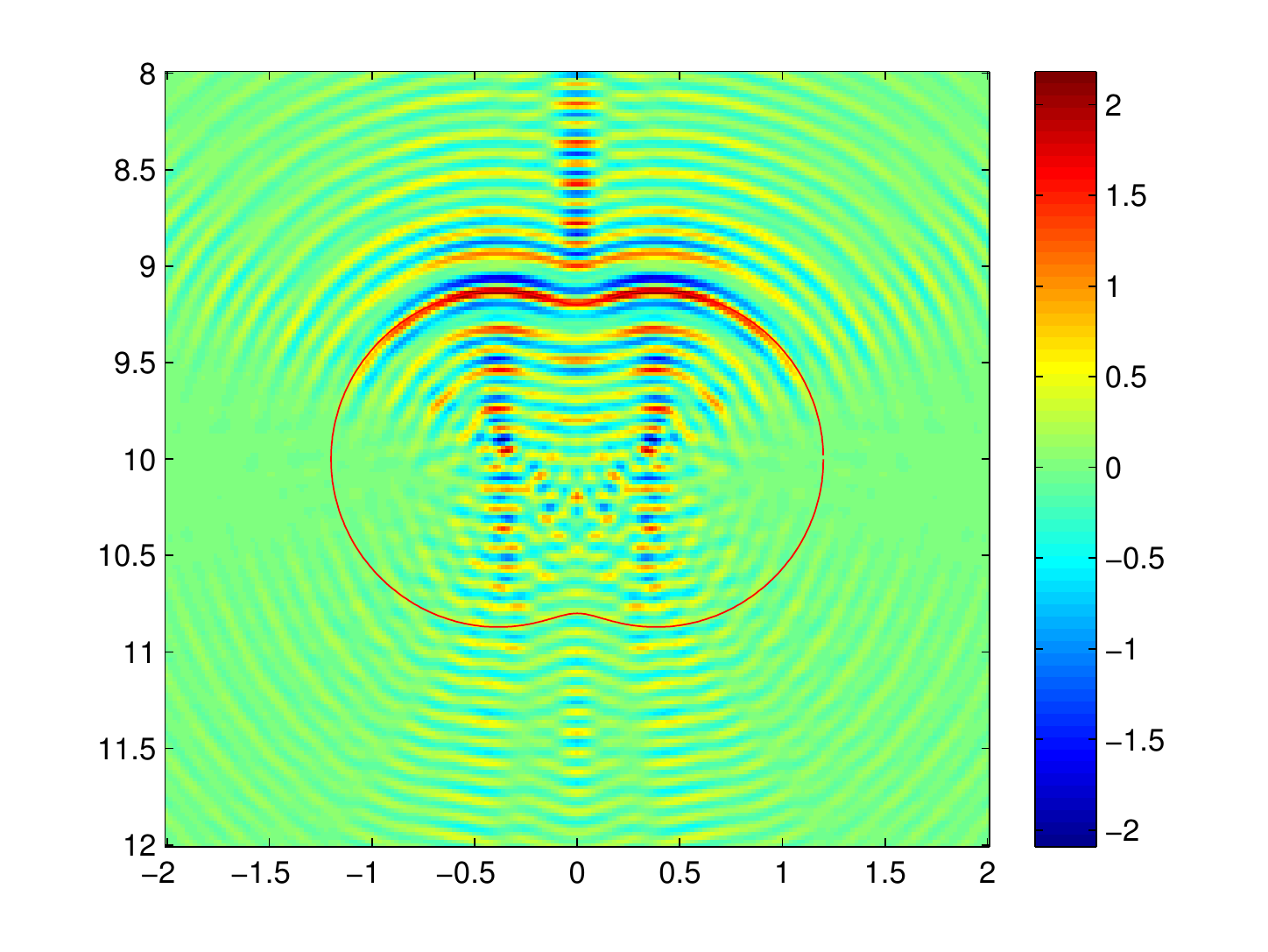}
	\includegraphics[width=0.32\textwidth]{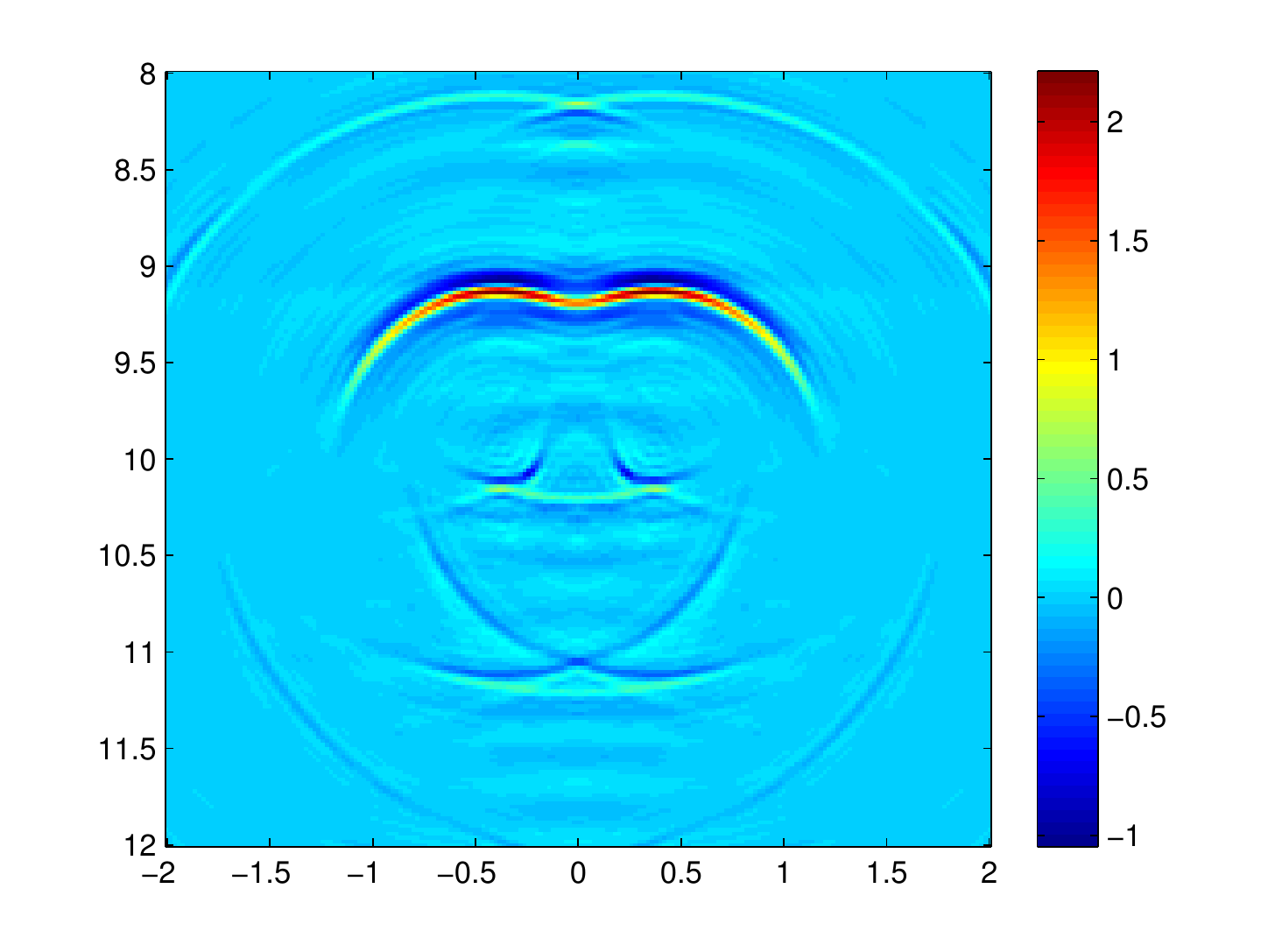}
	\includegraphics[width=0.32\textwidth]{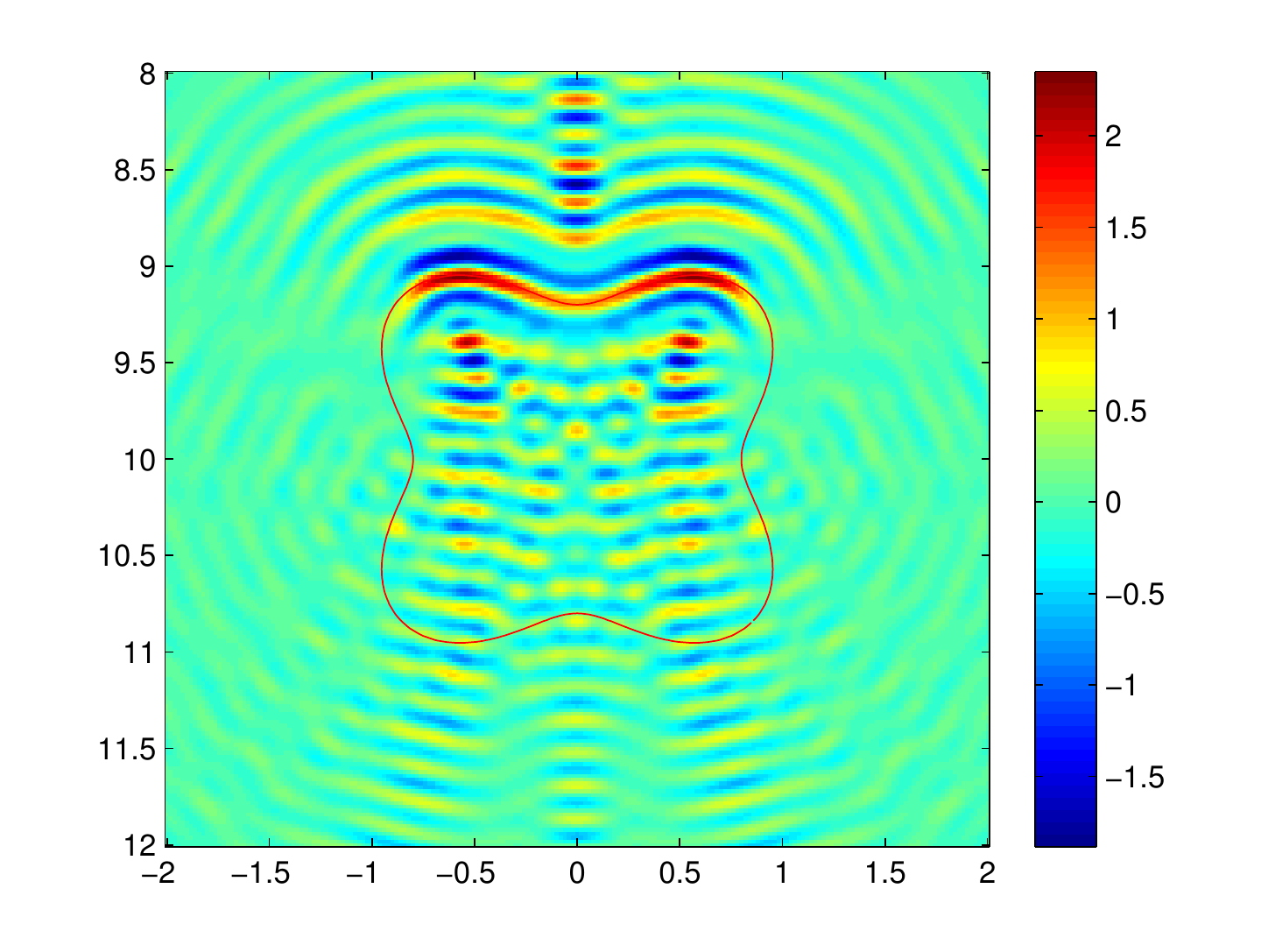}
	\includegraphics[width=0.32\textwidth]{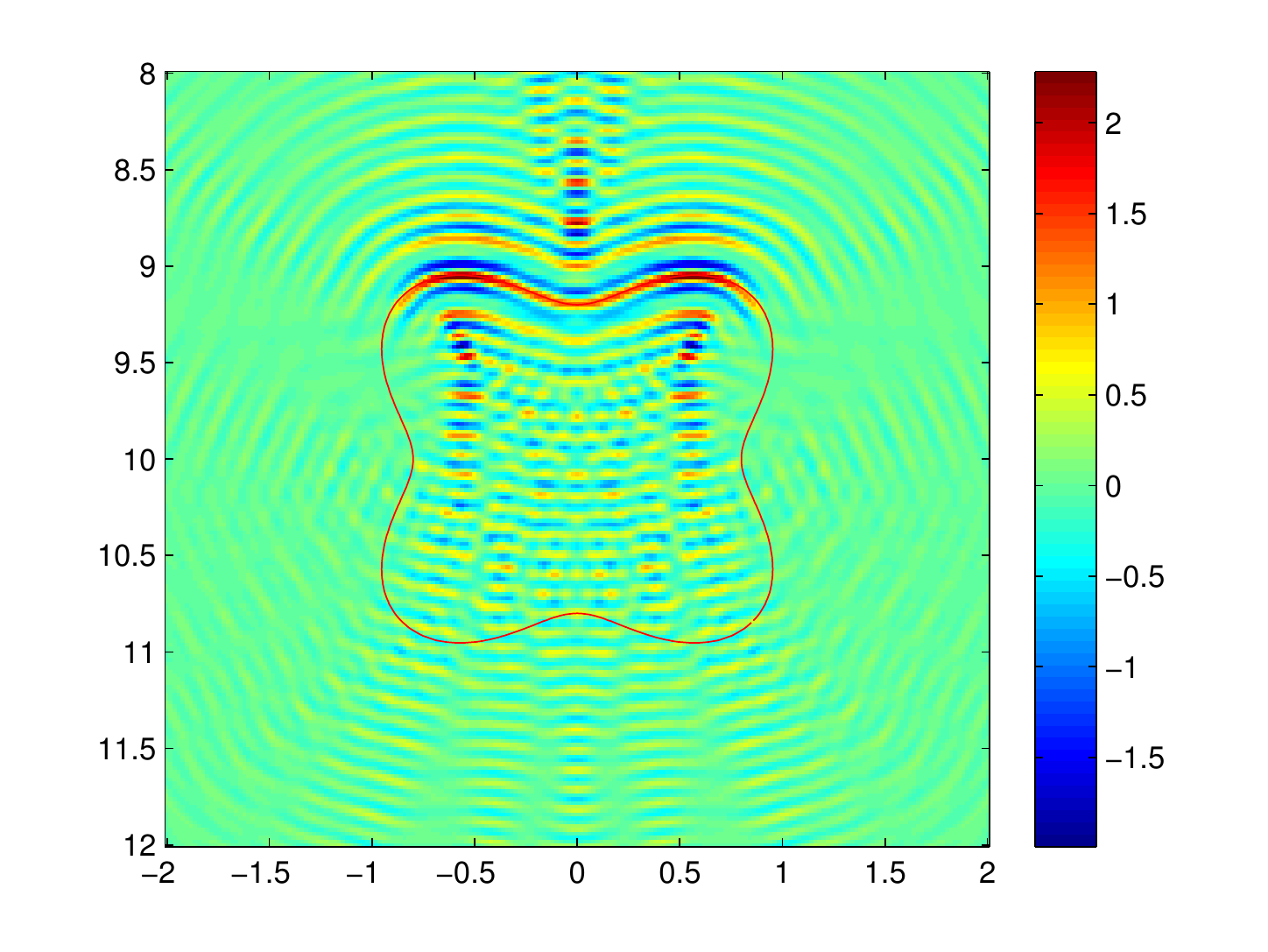}
	\includegraphics[width=0.32\textwidth]{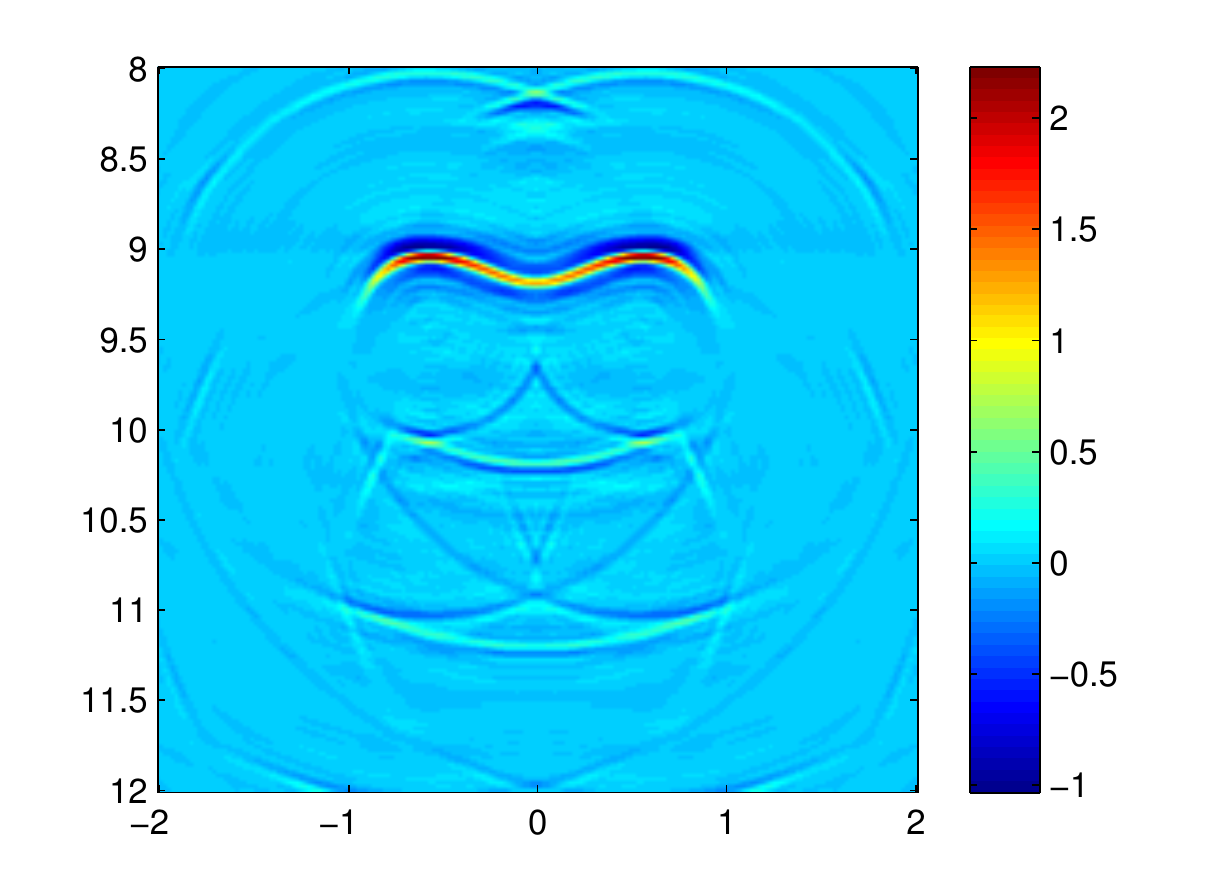}
	\includegraphics[width=0.32\textwidth]{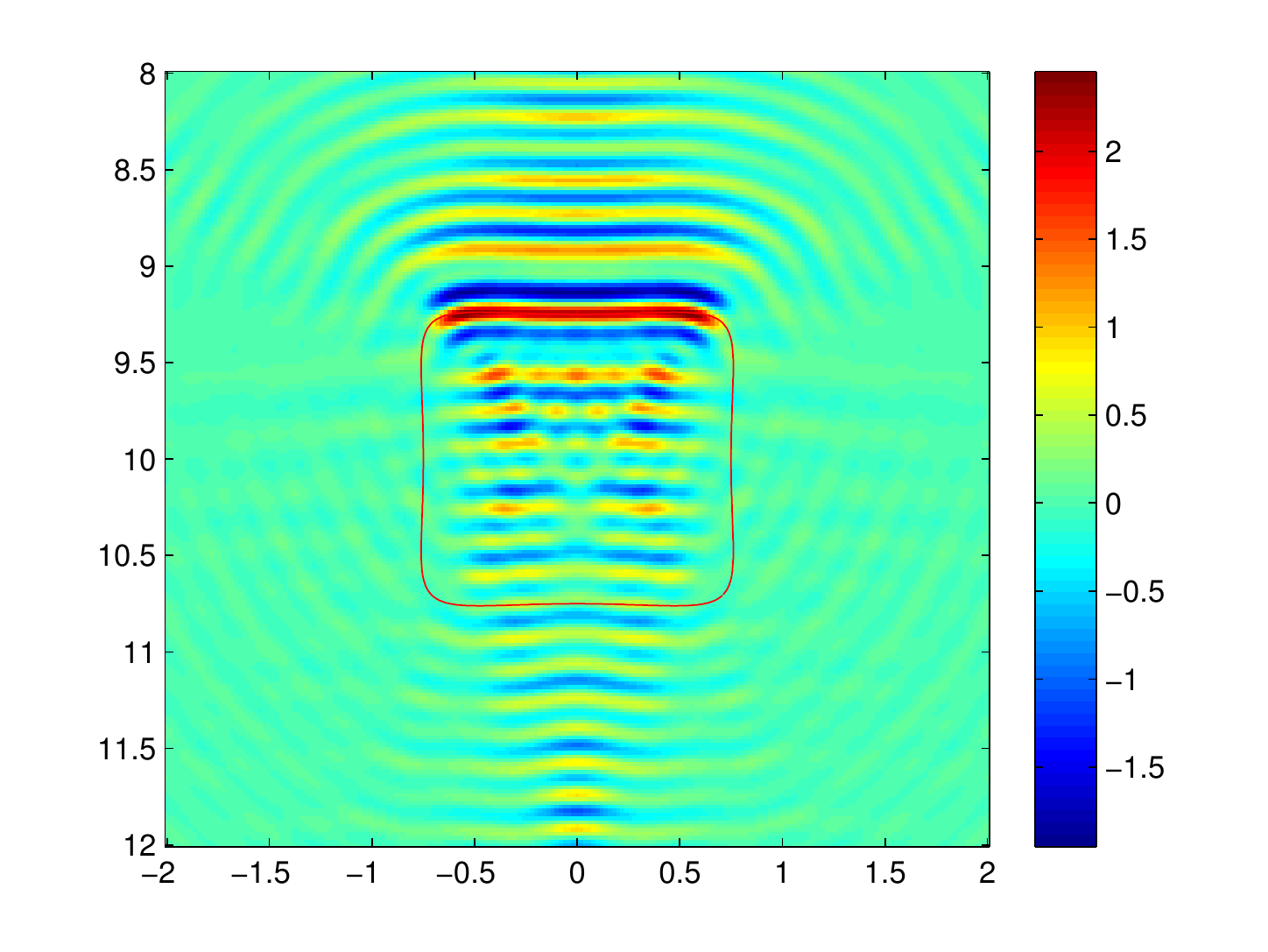}
	\includegraphics[width=0.32\textwidth]{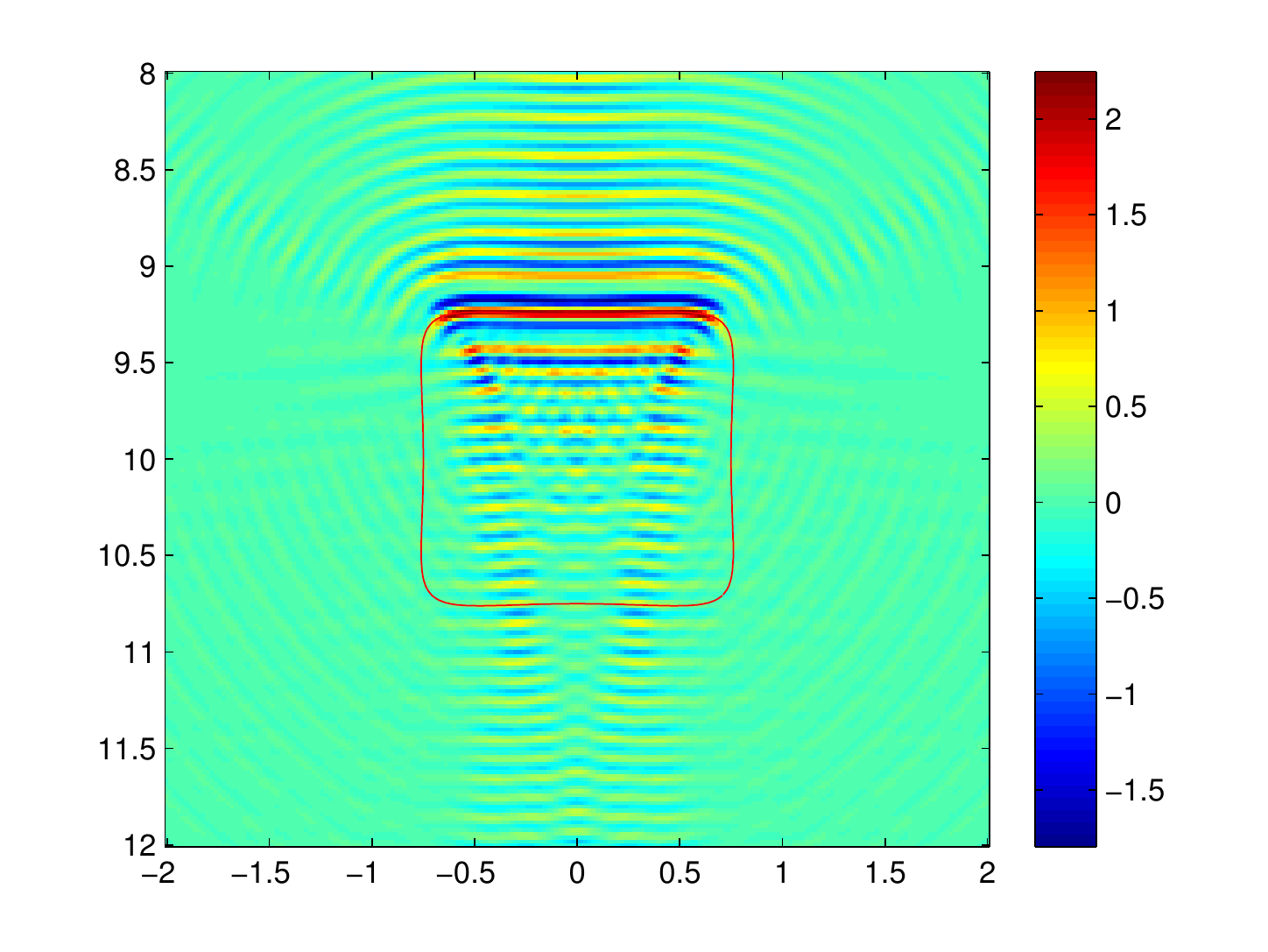}
	\includegraphics[width=0.32\textwidth]{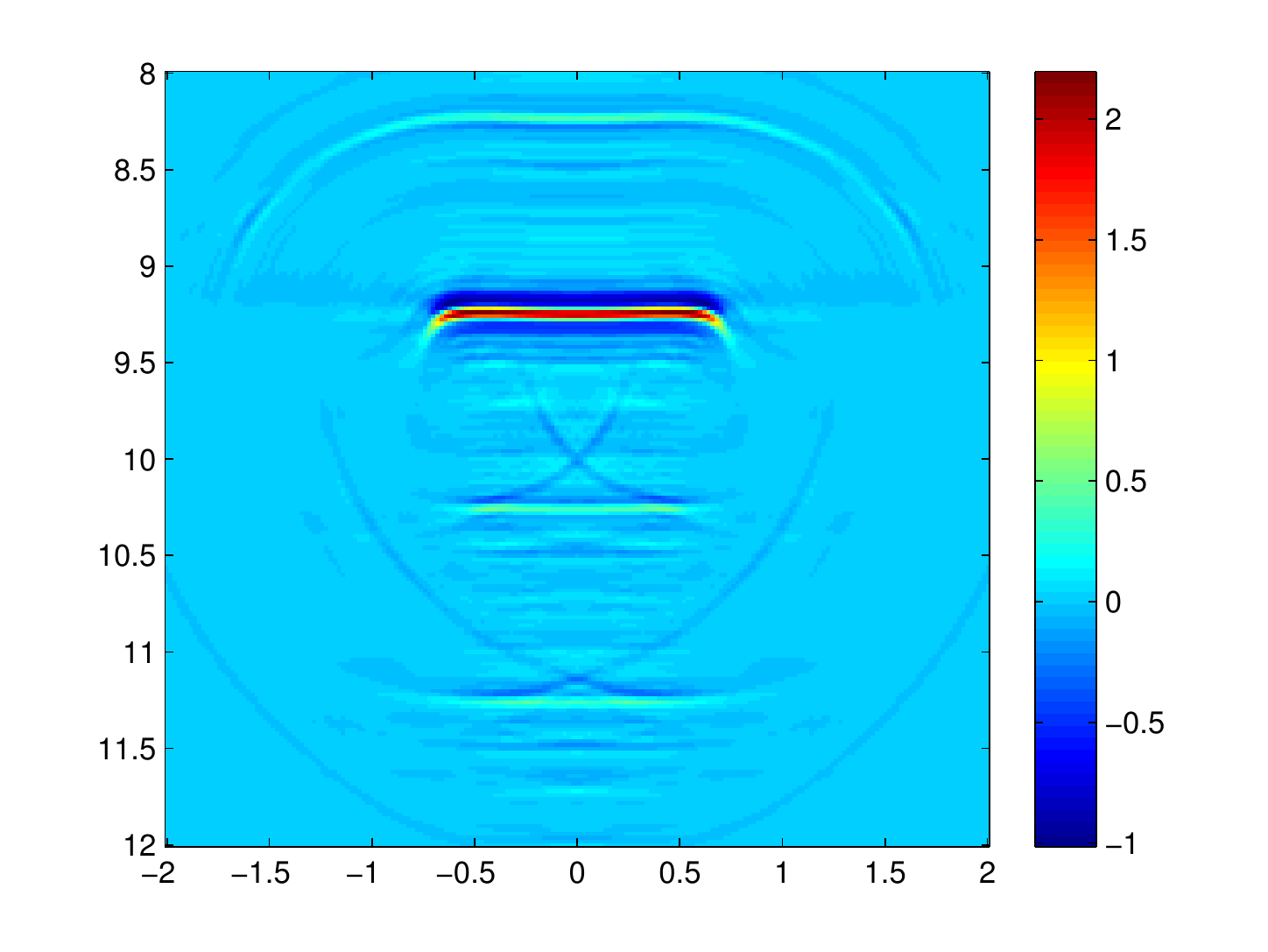}
	
	\caption{Example 2: Imaging results of Dirichlet obstacles
		with different shapes from top to below. The left column is imaged with single frequency data $\om=3\pi$, The middle column is imaged with single frequency data $\om=5\pi$ and the right column is imaged with multiple frequency data.}\label{figure_2}
\end{figure}

\bigskip
\textbf{Example 3} We consider the imaging of two Neumann obstacles. The first model
consists of two circles along horizontal direction and the second one is a circle and a
peanut along the vertical direction. The angular frequency $\om=3\pi$ for the test of the single frequency and $\om=\pi\times[2:0.5:8]$ for the test of multiple frequencies. Figure \ref{figure_31} shows the imaging result of the first model. The
imaging domain $\Om=(-4, 4) \times (8,12)$ with mesh size $401 \times 201$. We set $N_s = N_r = 301$. Figure \ref{figure_32}
 shows the imaging result of the second model. The
 imaging domain $\Om=(-4, 4)\times (8,12)$ with mesh size $401 \times 401$. We set $N_s = N_r = 301$. The multi-frequency RTM imaging results in Figure \ref{figure_31} and Figure \ref{figure_32} are obtained by adding the inmaging results from different frequencies. We observe from these two figures that imaging results can be greatly improved by stacking the multiple single frequency imaging results. 
\begin{figure}
	\centering
	\includegraphics[width=0.32\textwidth,height=0.16\textheight]{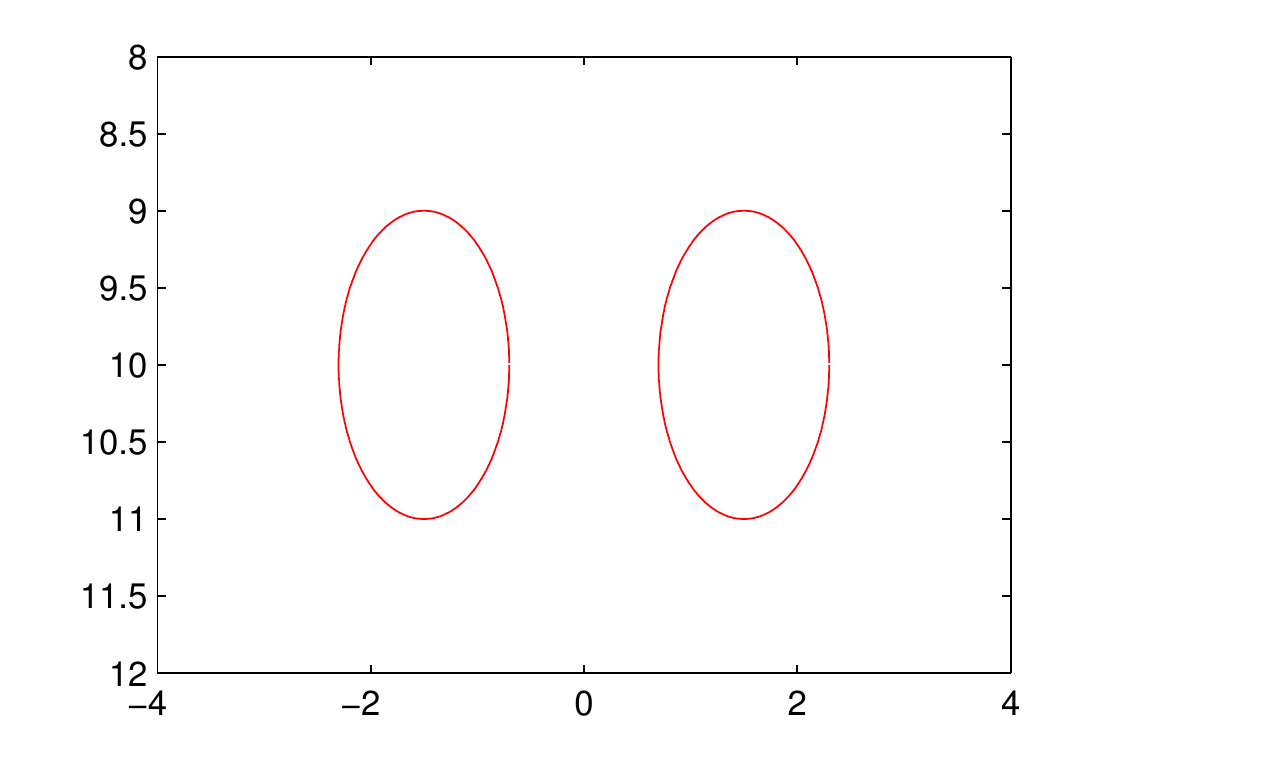}
	\includegraphics[width=0.32\textwidth]{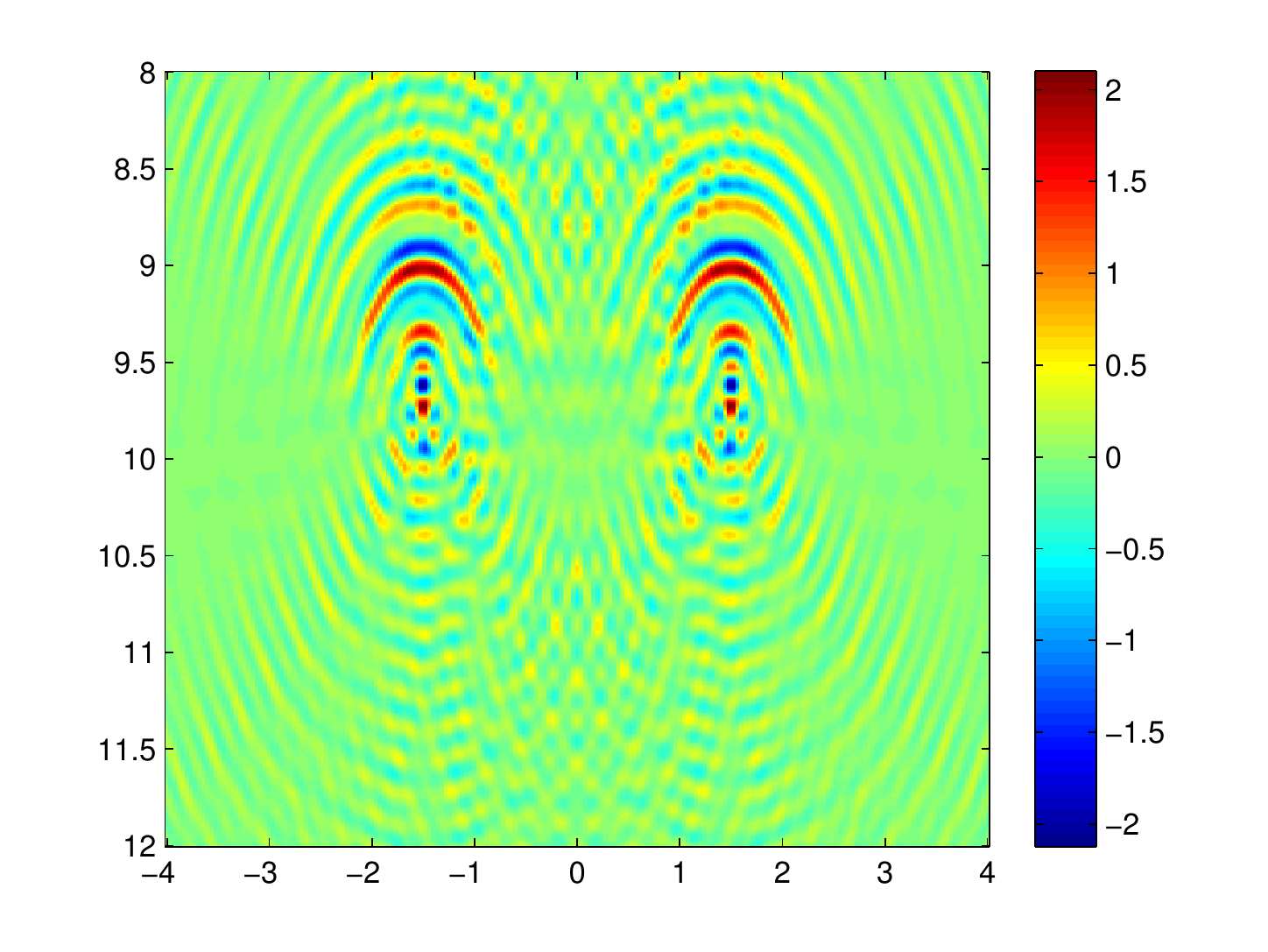}
	\includegraphics[width=0.32\textwidth]{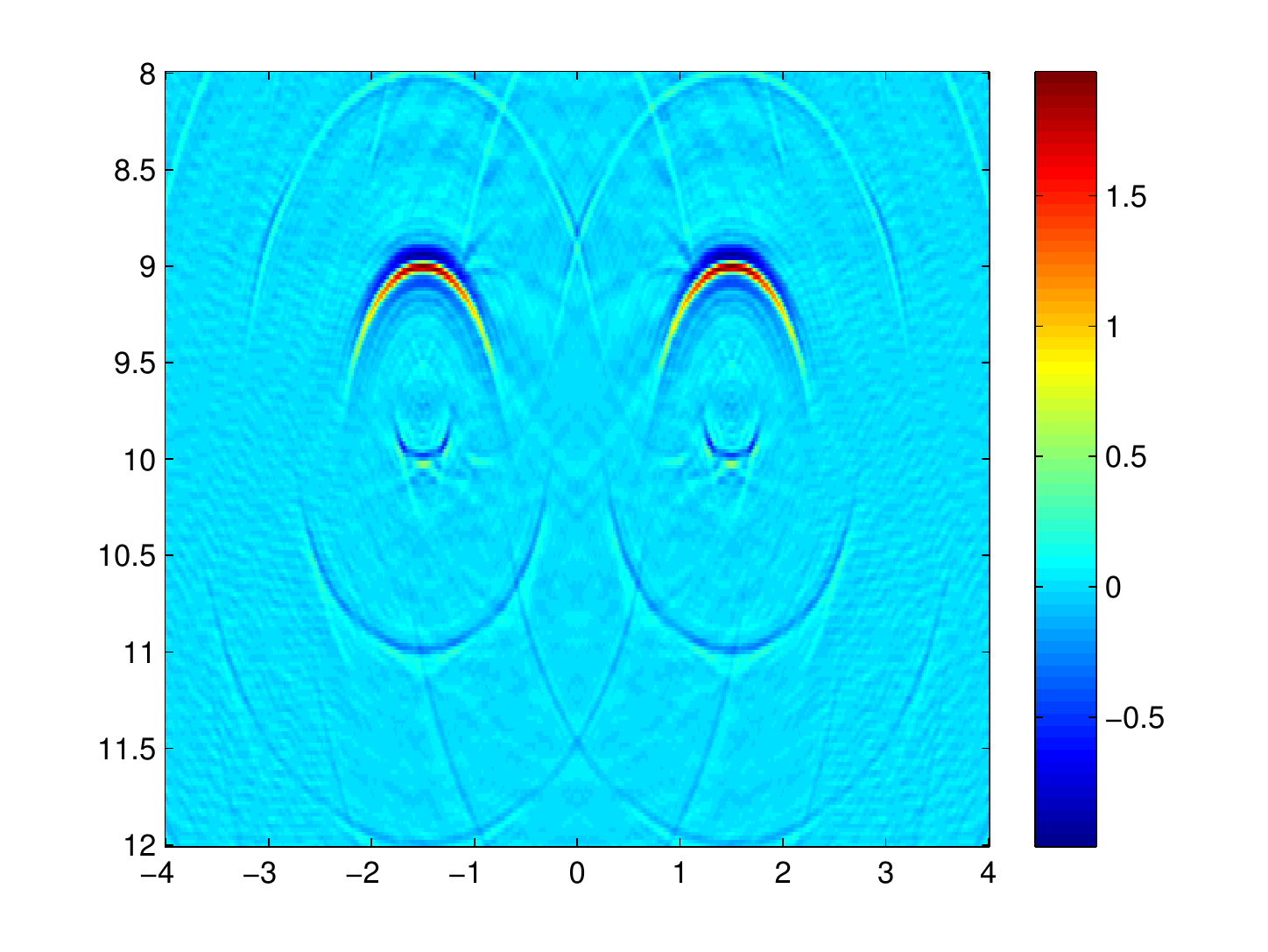}
	
	\caption{Example 3: From left to right,  true obstacle model with two circles, the imaging result
		with single frequency data $\om=3\pi$, the imaging result with multiple frequency data.}\label{figure_31}
\end{figure}

\begin{figure}
	\centering
	\includegraphics[width=0.32\textwidth,height=0.16\textheight]{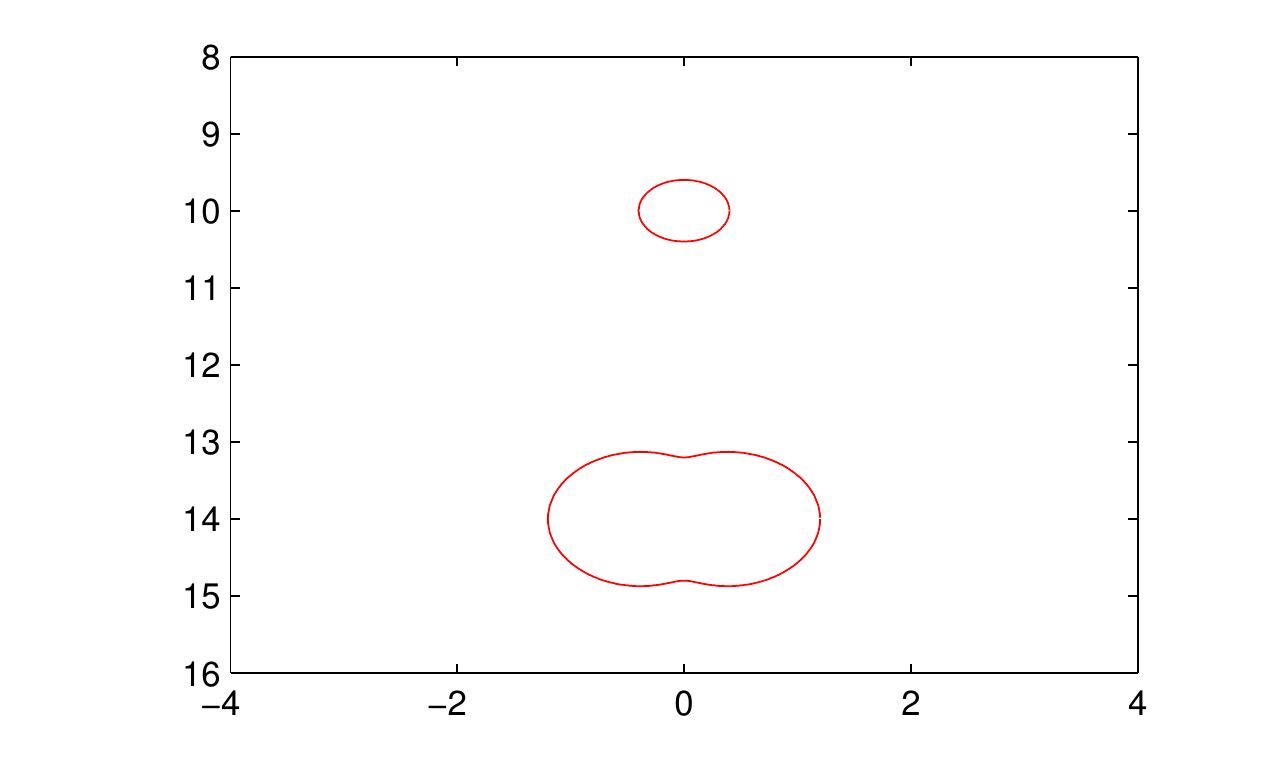}
	\includegraphics[width=0.32\textwidth]{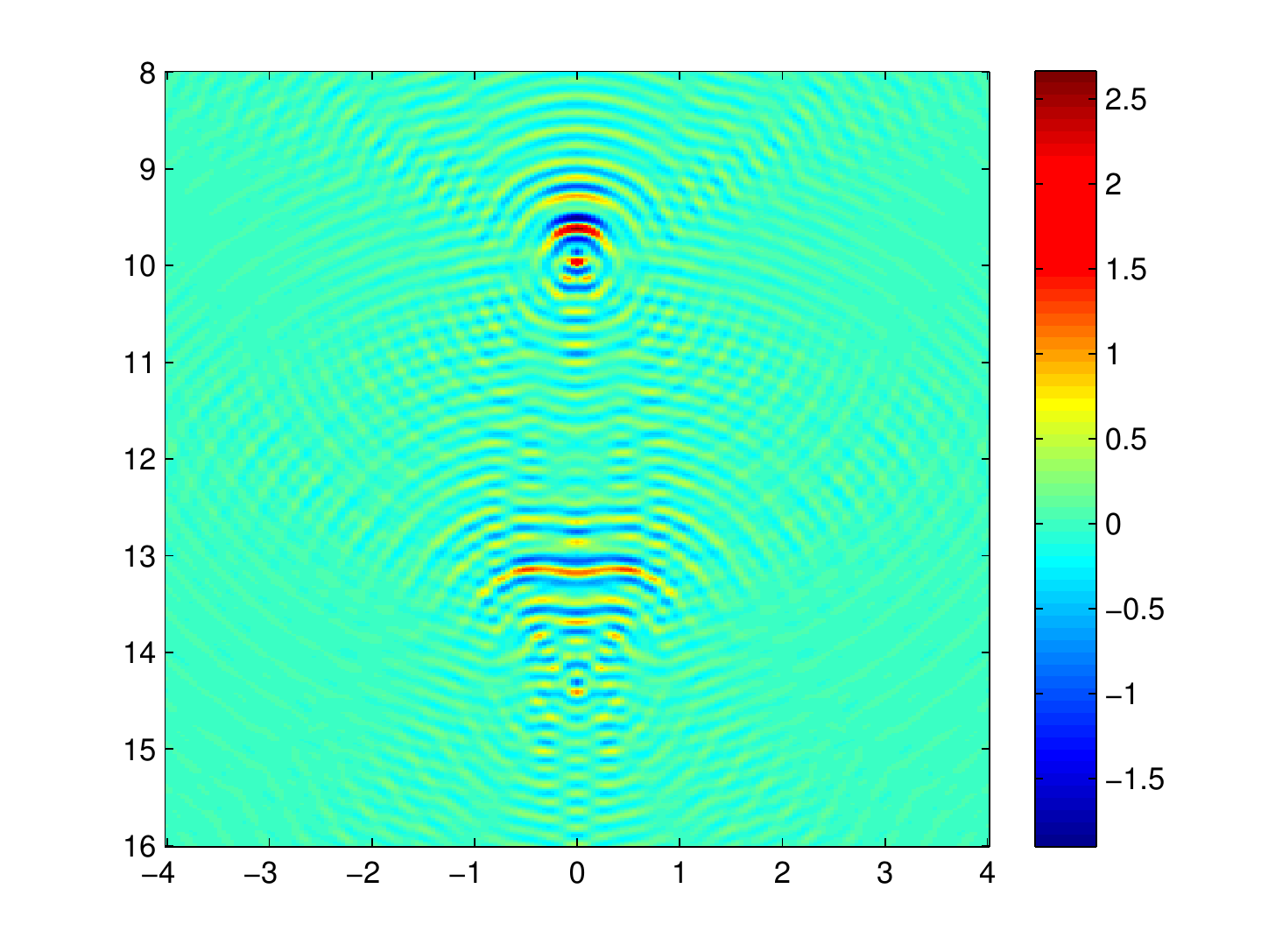}
	\includegraphics[width=0.32\textwidth]{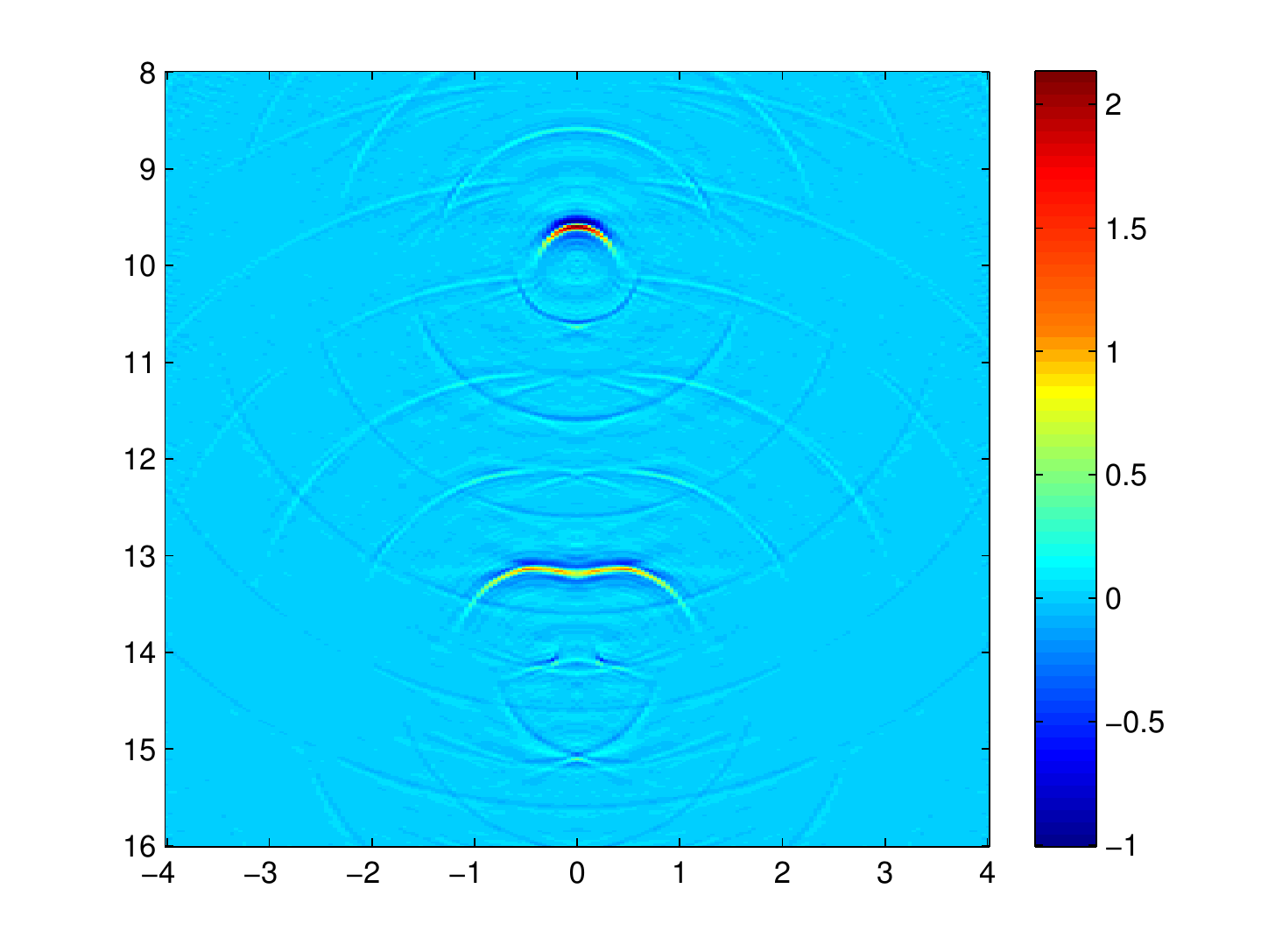}
	
	\caption{Example 3: From left to right,  true obstacle model with one circle and one peanut, the imaging result
		with single frequency data where $\om=3\pi$, the imaging result with multiple frequency data.}\label{figure_32}
\end{figure}

\bigskip
\textbf{Example 4}
We consider the stability of our half-space RTM imaging
algorithm with respect to the complex additive Gaussian random noise. We introduce
the additive Gaussian noise as $u_{\rm noise}=u_s+\nu_{\rm noise}$,
where $u_s$ is the synthesized data and $\nu_{\rm noise}$ is the Gaussian noise with mean zero and standard deviation $\sigma$ times the maximum of  the data $|u_s|$, i.e. $\nu_{\rm noise}=\frac{\sigma \max |u_s|}{\sqrt{2}}(\ep_1+\i\ep_2)$ and $\ep_i\thicksim \mathcal{N}(0,1)$.

\begin{figure}
	\centering
	\includegraphics[width=0.32\textwidth]{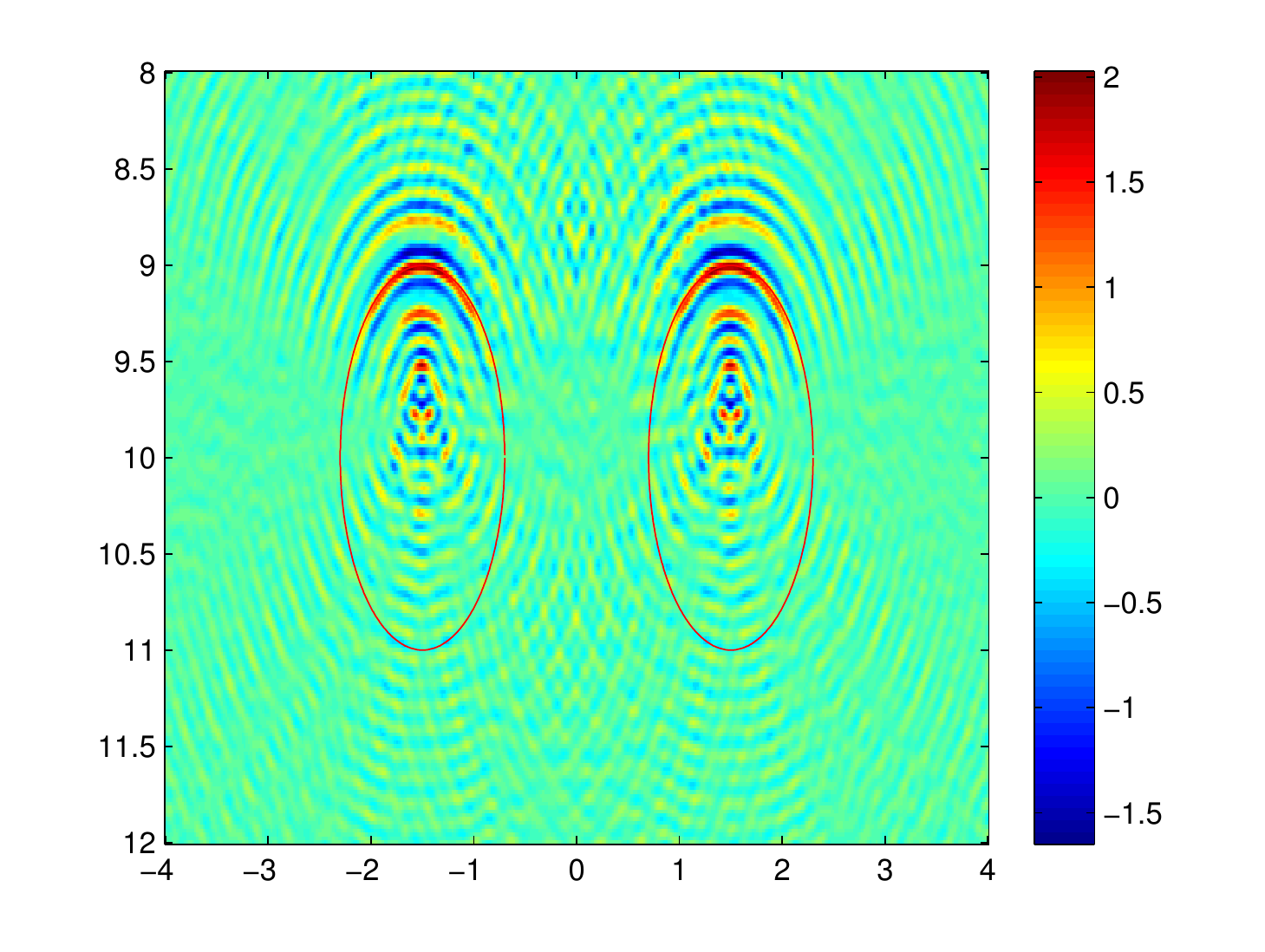}
	\includegraphics[width=0.32\textwidth]{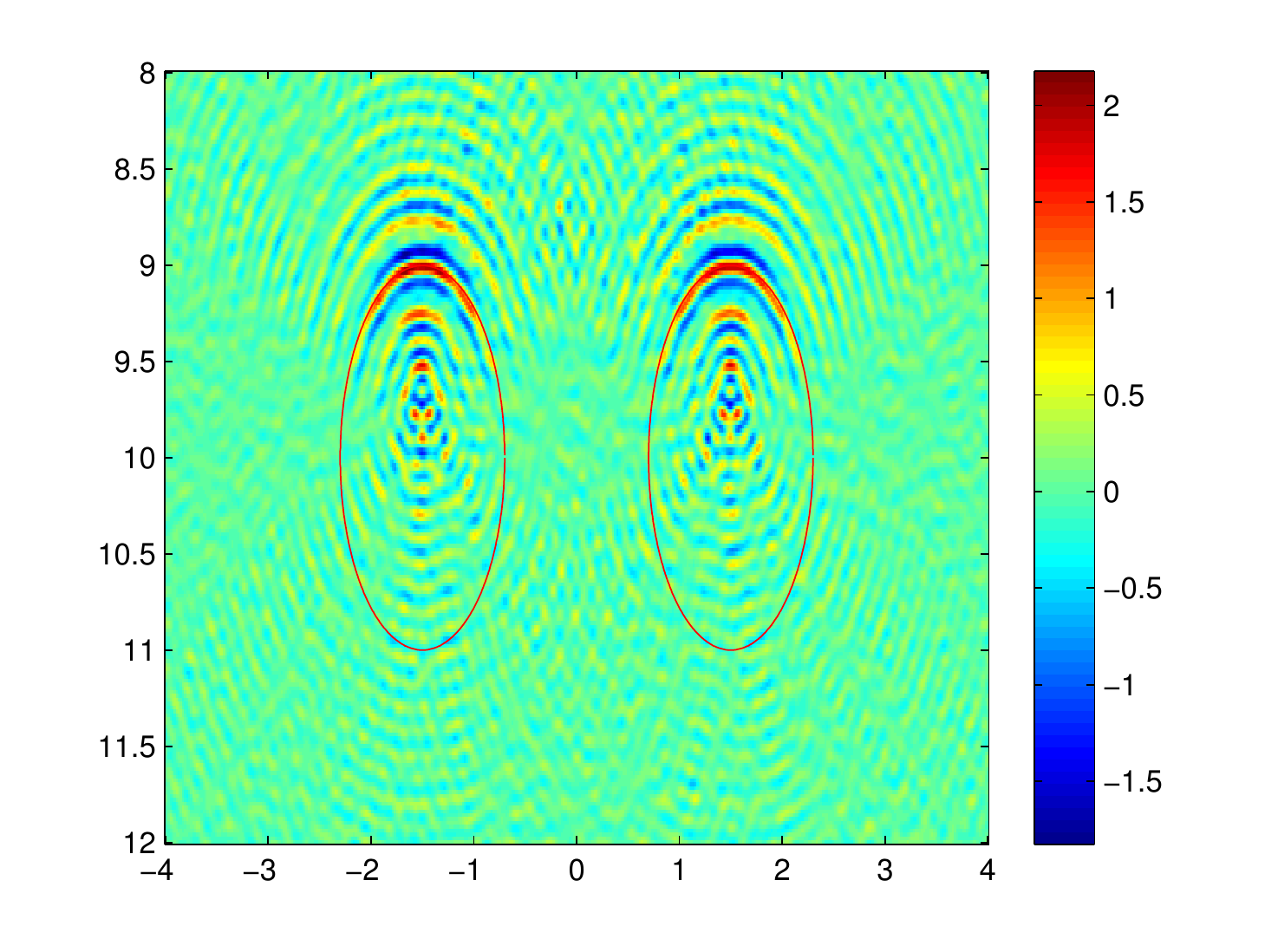}
	\includegraphics[width=0.32\textwidth]{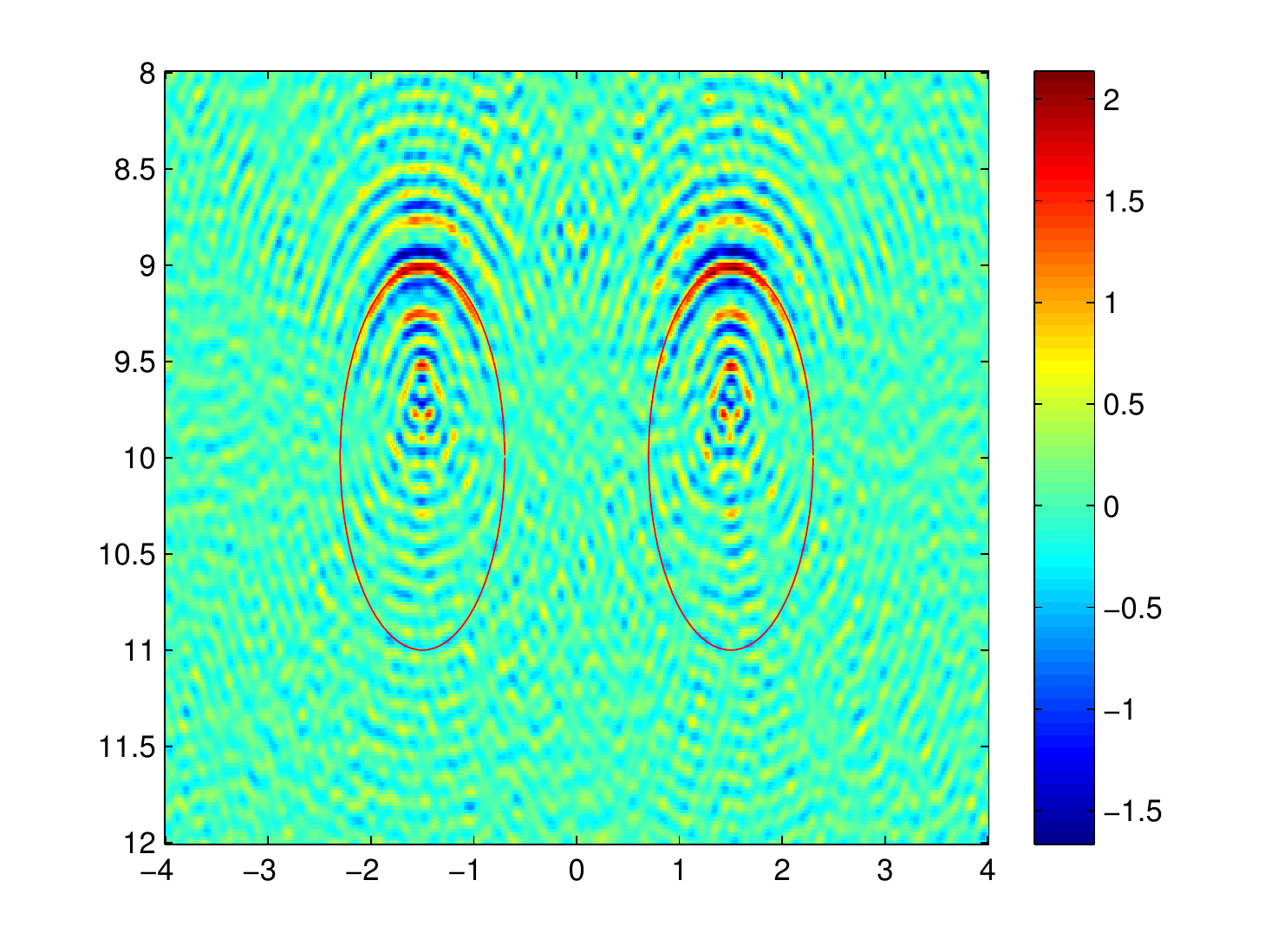}
	\includegraphics[width=0.32\textwidth]{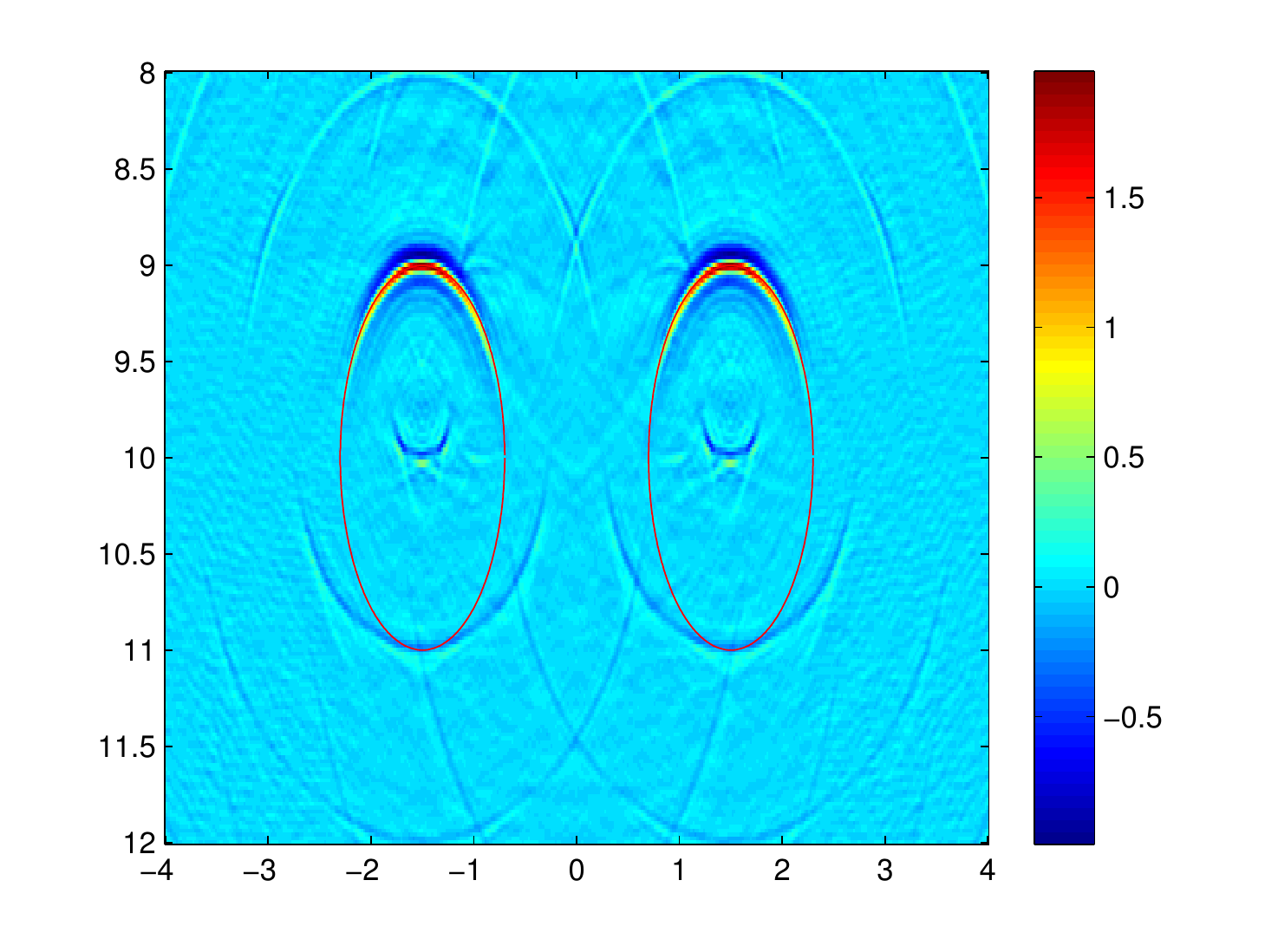}
	\includegraphics[width=0.32\textwidth]{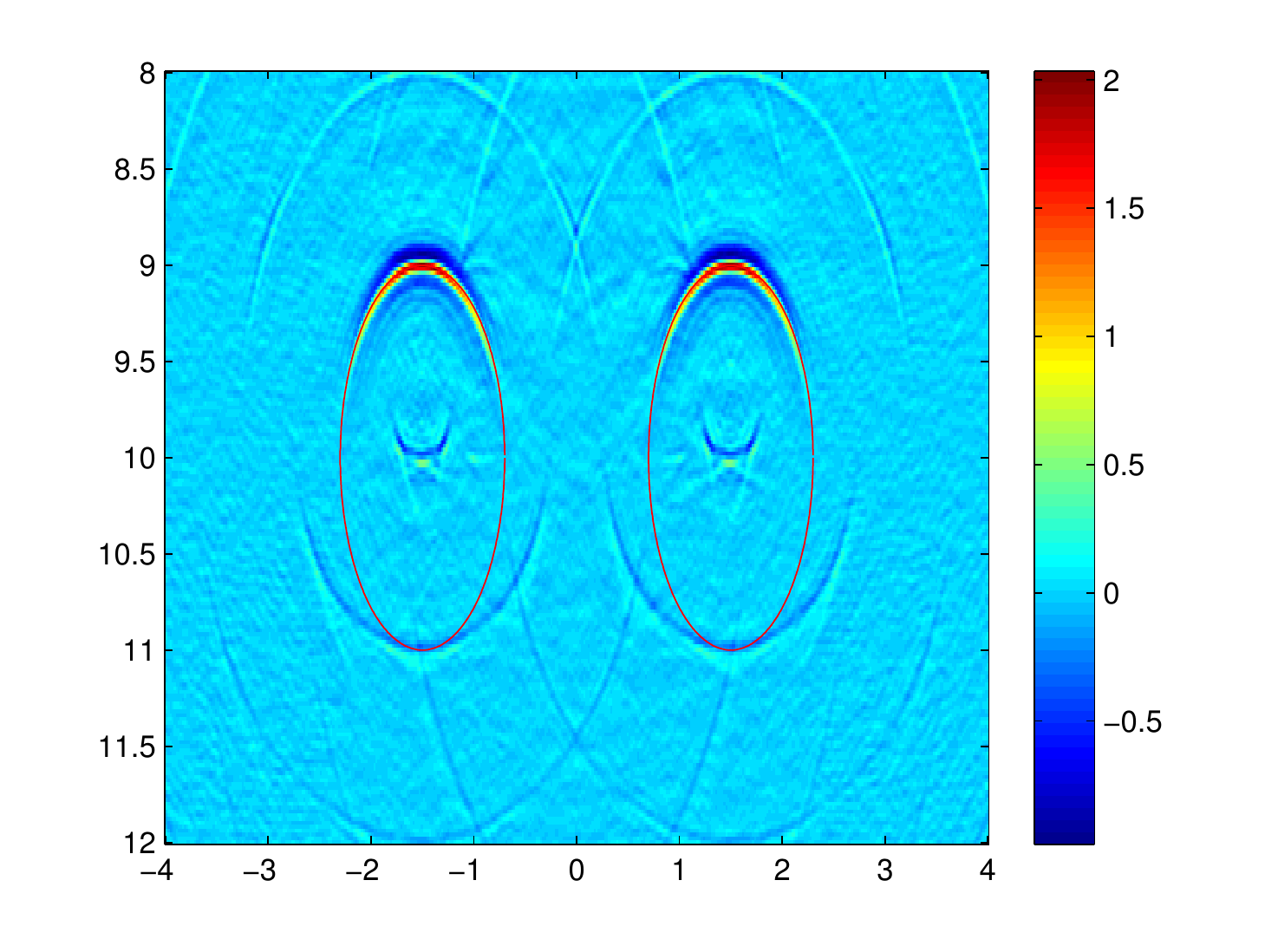}
	\includegraphics[width=0.32\textwidth]{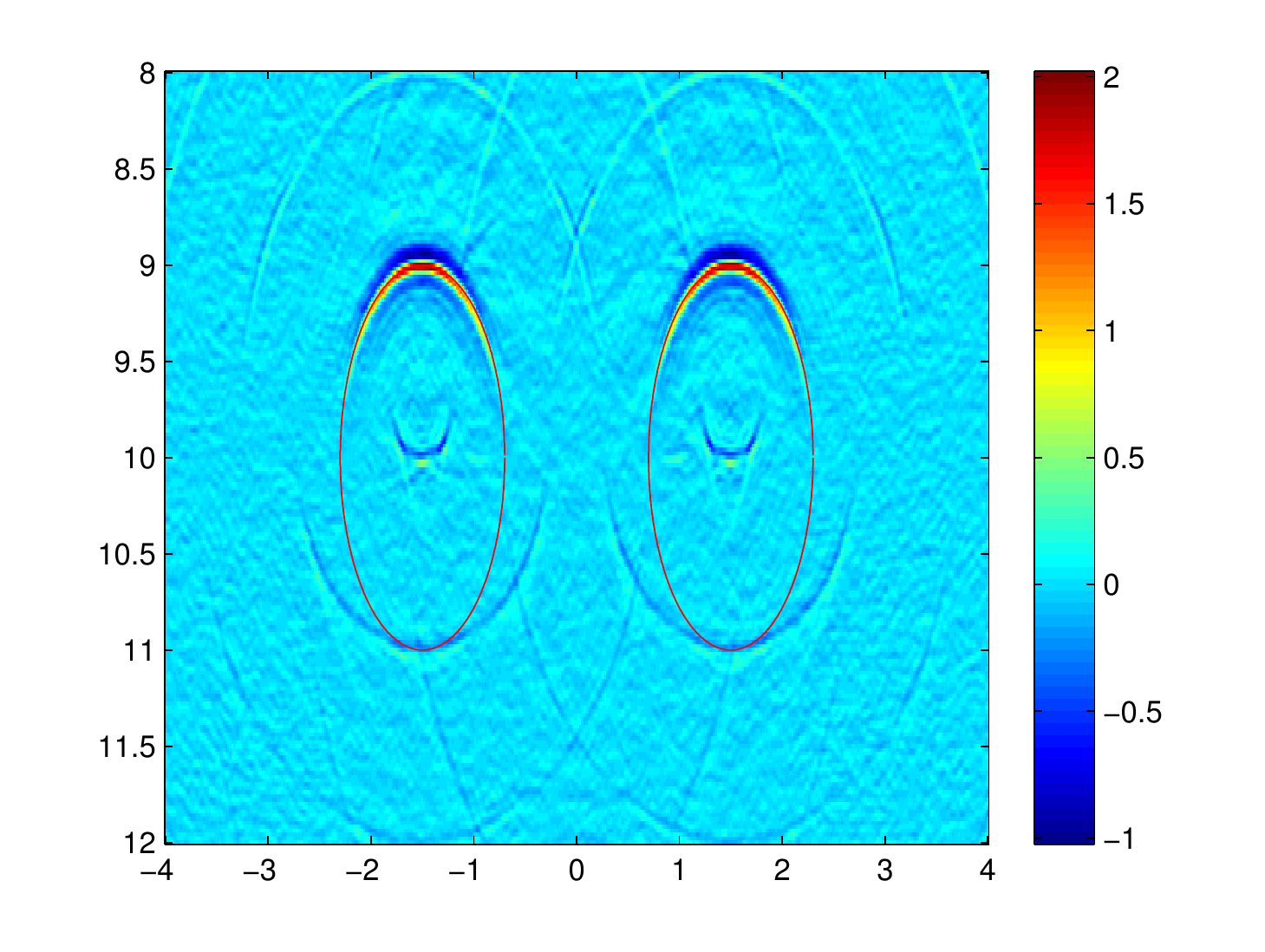}
	
	\caption{Example 4: Imaging results of a Dirichlet obstacle with noise levels $\mu =  0.2; 0.3; 0.4$ (from left to
		right). The top row is imaged with single frequency data $\om=4\pi$, and the
		bottom row is imaged with multi-frequency data.}\label{figure_4}
\end{figure}

Figure \ref{figure_4} shows the imaging results using single frequency data added with additive
Gaussian noise. The imaging quality can be greatly improved by using multi-frequency data $\omega = \pi\times [2:0.5:8]$.

\section{Appendix: Proof of Theorem \ref{thm:4.2}}

In this section we prove Theorem \ref{thm:4.2}.
Let $w(x)$ be the scattering solution of the problem:
\be\label{f2}
\Delta_e w + \omega^2 w=0 \ \ \mbox{\rm in } \R^2_+,\ \ 
\sigma(w)e_2=-\sigma(u_2)e_2 \ \ \mbox{\rm on } \Ga_0 .
\ee
Then $u_1-u_2-w$ is the scattering solution of the problem (\ref{e1}) with the boundary condition $u_1-u_2-w=-w$ on $\Gamma_D$. Thus by Theorem \ref{thm:4.1} and (\ref{q0}), we have
\be
\hskip-1cm\|\sigma(u_1-u_2)\nu\|_{H^{-1/2}(\Gamma_D)}&\leq&\|T_1(u_1-u_2-w)\|_{H^{-1/2}(\Gamma_D)}+\|\sigma(w)\nu\|_{H^{-1/2}(\Gamma_D)}\nn\\
&\leq&C (1+\|T_1\|)\max_{x\in \bar D}(|w(x)|+d_D|\nabla w(x)|),\label{f5}
\ee
where we recall that $T_1:H^{1/2}(\Ga_D)\to H^{-1/2}(\Ga_D)$ is the Dirichlet to Neumann mapping associated to the half-space elastic scattering problem (\ref{e1}) and $\|T_1\|$ denotes its operator norm.

By the integral representation formula, the scattering solution of the problem (\ref{f2}) satisfies
\be\label{f3}
w(y)\cdot e_j=\int_{\Ga_0} \sigma(u_2(x))e_2\cdot \N(x,y) e_j ds(x),\ \ \forall y\in \R^2,\ j=1,2.
\ee
On the other hand, by the integral representation formula, we have
\ben
u_2(x)\cdot e_j=\GG(u_2(\cdot),\G(\cdot,x)e_j),\ \ \forall x\in\Ga_0, \ j=1,2,
\een
where $\GG(\cdot,\cdot)$ is defined in (\ref{g1}) and $\G(\cdot,\cdot)$ is the fundamental solution tensor of the elastic wave equation introduced in section 2. For any $x\in\Ga_0, z\in\R^2$, denote by $\T(z,x)\in\C^{2\times 2}$ the traction tensor, $\T(z,x)q=\sigma(\G(z,x)q)e_2,\forall q\in\R^2$. The $(i,j)$-th element of $\T(z,x)$ is 
\ben
[\T(z,x)]_{ij}=[\sigma(\G(z,x)e_j)e_2]e_i,\ \ i,j=1,2.
\een
Simple calculation shows that
\ben
\sigma(u_2(x))e_2\cdot e_i=\GG(u_2(\cdot),\T(\cdot,x)^Te_i),\ \ \forall x\in\Ga_0, i=1,2,
\een
which yields from (\ref{f3}) that
\be
\hskip-1cmw(y)\cdot e_j&=&\GG(u_2(\cdot),\left[\int_{\Ga_0}\sum^2_{i=1}[\T(\cdot,x)^Te_i]\cdot [e_i^T\N(x,y)e_j]ds(x)\right])\nn\\
\hskip-1cm&=&\GG(u_2(\cdot),\V(\cdot,y)e_j),\label{f4}
\ee
where 
\ben
\V(z,y)=\int_{\Ga_0}\T(z,x)^T\N(x,y)ds(x),\ \ \forall y\in \R^2, z\in\Ga_D.
\een
Notice that $\|\sigma(u_2)\nu\|_{H^{-1/2}(\Gamma_D)}\leq \|T_2\| \|g\|_{H^{1/2}(\Gamma_D)}$, where $T_2:H^{1/2}(\Ga_D)\to H^{-1/2}(\Ga_D)$ is the Dirichlet to Neumann mapping associated to the scattering problem (\ref{e2}) and $\|T_2\|$ denotes its operator norm. We obtain from (\ref{f4}) and (\ref{q0}) that
\be
\hskip-1.5cm |w(y)|+d_D|\na w(y)|&\le&C (1+\|T_2\|) \|g\|_{H^{1/2}(\Gamma_D)}
\max_{z\in \Ga_D}\sum^1_{i,j=0}d_D^{i+j}|\na_z^i\na_y^j\V(z,y)|.\label{f6}
\ee
To estimate the term involving $\V(z,y)$, we use Parserval identity and Lemma \ref{lem:2.2} to obtain
\ben
\V(z,y)&=&\frac 1{2\pi}\,\pv\int_{\R}\hat{\T}(z_2;\xi,0)^T\hat{\N}(\xi,0;y_2)e^{-\i\xi(y_1+z_1)}d\xi\\
&&-\frac\i 2
\left[\hat{\T}(z_2;\xi,0)^T\hat{\N}(\xi,0;y_2)e^{-\i\xi(y_1+z_1)}\right]^{k_R}_{-k_R}.
\een
It is easy to see from (\ref{G1})-(\ref{G2}) that
\ben \hspace{-1cm}
\hat{\T}(z_2;\xi,0)&=&\frac \mu{2\om^2}\Bigg( \begin{array}{cc}
	\beta & \frac{\xi\beta}{\mu_s} \\
	2\xi\mu_s & 2\xi^2
\end{array} \Bigg)e^{\i\mu_s z_2}
+\frac{\mu}{2\om^2}
\Bigg(\begin{array}{cc}
	2\xi^2 & -2\xi\mu_p \\
	-\frac{\xi\beta}{\mu_p} & \beta 
\end{array}\Bigg)e^{\i\mu_p z_2} \\
&:=&\tilde{\mathbb{T}}_{s}(\xi)e^{\i\mu_p z_2}+\tilde{\mathbb{T}}_{p}(\xi)e^{\i\mu_s z_2}.
\een
Now by using (\ref{d2}) we have
\ben
\hskip-1cm\V(z,y)&=&\frac 1{2\pi}\sum_{\al,\beta=p,s}{\rm p.v.}\int_{\R}\frac{\tilde{\mathbb T}_{\alpha}(\xi)^T{\mathbb N}_\beta(\xi)}{\de(\xi)}e^{\i(\mu_\al z_2+\mu_\beta y_2)-\i(y_1+z_1)\xi}d\xi\\
\hskip-1cm& &-\frac\i 2\sum_{\al,\beta=p,s}\left[\frac{\tilde{\mathbb T}_\alpha(\xi)^T{\mathbb{N}}_\beta(\xi)}{\de'(\xi)}e^{\i(\mu_\al z_2+\mu_\beta y_2)-\i(y_1+z_1)\xi}\right]^{k_R}_{-k_R}:={\rm V}_1+{\rm V}_2.
\een
To estimate ${\rm V}_1$ we split the integral into two domains $(-k_s,k_s)$ and $\R\bks[-k_s,k_s]$ and use the Van der Corput lemma \ref{van} to estimate the integral in the first interval and the argument in Lemma \ref{lem:3.4}
to estimate
the integral in $\R\bks[-k_s,k_s]$. This yields $|{\rm V}_1|\le C\mu^{-1}(k_sh)^{-1/2}$. By the same argument as in
Lemma \ref{lem:3.2} we can show $|{\rm V}_2|\le C\mu^{-1}e^{-\sqrt{k_R-k_s}h}$. This shows
\ben
\max_{z\in\Ga_D}|\V(z,y)|\le \frac C\mu (k_sh)^{-1/2},\ \ \forall y\in \bar D.
\een
A similar argument shows that
\ben
\max_{z\in\Ga_D} k_s^{i+j}|\na_z^i\na^j_y\V(z,y)|\le \frac C\mu (k_sh)^{-1/2},\ \ \forall y\in \bar D,\ i,j=0,1.
\een
Substitute the above two estimates into (\ref{f6}) we obtain
\ben
\hskip-1cm\max_{y\in \bar D}(|w(y)|+d_D|\na w(y)|)\le\frac C\mu (1+\|T_2\|)(1+k_sd_D)^2(k_sh)^{-1/2} \|g\|_{H^{1/2}(\Gamma_D)}.
\een
This completes the proof of the theorem from (\ref{f5}).
\finproof

{\bf Acknowledgement.} The first author is grateful for the partial support by the China NSF under the grant 118311061. The authors are indebted to Dr. Guanghui Huang from Michigan State University for helpful discussions.

\end{document}